\documentclass[aos,preprint]{imsart}

\usepackage{graphicx}
\usepackage{amssymb}
\usepackage{epstopdf}
\RequirePackage[authoryear]{natbib}
\usepackage{amsmath, mathabx}
\usepackage{mathrsfs}
\usepackage{hyperref}
\usepackage{comment}
\usepackage{xcolor}
\usepackage{subfigure}
\usepackage{enumitem}
\usepackage{algorithm}
\usepackage{algpseudocode}
\usepackage{stmaryrd}
\usepackage[T1]{fontenc}

\startlocaldefs
\newcommand{\X}{{\mathbf{X}}}
\newcommand{\A}{{\mathbf{A}}}
\newcommand{\Y}{{\mathbf{Y}}}
\newcommand{\B}{{\mathbf{B}}}

\newcommand{\D}{{\mathbf{D}}}
\newcommand{\U}{{\mathbf{U}}}
\newcommand{\I}{{\mathbf{I}}}
\newcommand{\Z}{{\mathbf{Z}}}

\newcommand{\bcX}{{\mathbfcal{X}}}
\newcommand{\bcA}{{\mathbfcal{A}}}
\newcommand{\bcY}{{\mathbfcal{Y}}}

\newcommand{\bcZ}{{\mathbfcal{Z}}}
\newcommand{\bcB}{{\mathbfcal{B}}}

\newcommand{\bcW}{{\mathbfcal{W}}}

\newcommand{\bcM}{{\mathbfcal{M}}}

\newcommand{\bk}{\mathbf{k}}
\newcommand{\bn}{\mathbf{n}}
\newcommand{\bu}{\mathbf{u}}
\newcommand{\bv}{\mathbf{v}}
\newcommand{\bx}{\mathbf{x}}
\newcommand{\bi}{\mathbf{i}}
\newcommand{\bj}{\mathbf{j}}

\newcommand{\cG}{\mathcal{G}}
\newcommand{\cC}{\mathcal{C}}
\newcommand{\cH}{\mathcal{H}}
\newcommand{\cE}{\mathcal{E}}
\newcommand{\cS}{\mathcal{S}}

\newcommand{\cL}{\mathcal{L}}
\newcommand{\cV}{\mathcal{V}}
\newcommand{\cM}{\mathcal{M}}
\newcommand{\cP}{\mathcal{P}}
\newcommand{\cX}{\mathcal{X}}
\newcommand{\cY}{\mathcal{Y}}

\newcommand{\sX}{\mathscr{X}}

\newcommand{\bbP}{\mathbb{P}}

\newcommand{\bbE}{\mathbb{E}}
\newcommand{\bbR}{\mathbb{R}}
\newcommand{\bbZ}{\mathbb{Z}}

\newcommand{\bsigma}{\boldsymbol{\sigma}}

\newcommand{\AllAlg}{{\rm{AllAlg}}}
\newcommand{\PolyAlg}{{\rm{PolyAlg}}}

\renewcommand{\exp}{{\rm{exp}}}
\newcommand{\TV}{{\rm{TV}}}
\newcommand{\KL}{{\rm{KL}}}
\newcommand{\HPC}{{\rm{HPC}}}
\newcommand{\HPDS}{{\rm{HPDS}}}
\newcommand{\CHC}{{\rm{CHC}}}
\newcommand{\ROHC}{{\rm{ROHC}}}

\newcommand{\NORM}{{\rm{NORM}}}
\newcommand{\HS}{{\rm{HS}}}
\newcommand{\SUM}{{\rm{SUM}}}
\newcommand{\Bern}{{\rm{Bern}}}
\newcommand{\Bin}{{\rm{Bin}}}
\newcommand{\RK}{{\rm{RK}}}
\newcommand{\id}{{\rm{id}}}
\newcommand{\LR}{{\rm{LR}}}

\newcommand{\BC}{{\rm{BC}}}
\newcommand{\ROS}{{\rm{ROS}}}

\newcommand{\argmax}{\mathop{\rm arg\max}}

\newcommand{\prodk}{\prod_{i=1}^d k_i}
\newcommand{\prodn}{\prod_{i=1}^d n_i}

\newtheorem{Lemma}{Lemma}
\newtheorem{Proposition}{Proposition}
\newtheorem{Theorem}{Theorem}

\newtheorem{Conjecture}{Conjecture}

\newtheorem{Remark}{Remark}

\DeclareMathAlphabet\mathbfcal{OMS}{cmsy}{b}{n}

\endlocaldefs

\begin{document}

\begin{frontmatter}
\title{Tensor Clustering with Planted Structures: Statistical Optimality and Computational Limits}
\runtitle{Tensor Clustering with Planted Structures}

\begin{aug}
\author[A]{\fnms{Yuetian} \snm{Luo}\ead[label=e1]{yluo86@wisc.edu}}
\and
\author[A,B]{\fnms{Anru R.} \snm{Zhang}\ead[label=e2]{anru.zhang@duke.edu}}
\address[A]{Department of Statistics, University of Wisconsin-Madison, \printead{e1,e2}}
\address[B]{Department of Biostatistics \& Bioinformatics, Duke University, \printead{e2}}
\end{aug}

\begin{abstract}
This paper studies the statistical and computational limits of high-order clustering with planted structures. We focus on two clustering models, constant high-order clustering (CHC) and rank-one higher-order clustering (ROHC), and study the methods and theory for testing whether a cluster exists (detection) and identifying the support of cluster (recovery). 
    
Specifically, we identify the sharp boundaries of signal-to-noise ratio for which CHC and ROHC detection/recovery are statistically possible. We also develop the tight computational thresholds: when the signal-to-noise ratio is below these thresholds, we prove that polynomial-time algorithms cannot solve these problems under the computational hardness conjectures of hypergraphic planted clique (HPC) detection and hypergraphic planted dense subgraph (HPDS) recovery. We also propose polynomial-time tensor algorithms that achieve reliable detection and recovery when the signal-to-noise ratio is above these thresholds. Both sparsity and tensor structures yield the computational barriers in high-order tensor clustering. The interplay between them results in significant differences between high-order tensor clustering and matrix clustering in literature in aspects of statistical and computational phase transition diagrams, algorithmic approaches, hardness conjecture, and proof techniques. To our best knowledge, we are the first to give a thorough characterization of the statistical and computational trade-off for such a double computational-barrier problem. Finally, we provide evidence for the computational hardness conjectures of HPC detection (via low-degree polynomial and Metropolis methods) and HPDS recovery (via low-degree polynomial method).
\end{abstract}



\end{frontmatter}

\section{Introduction}\label{sec:intro}

The high-dimensional tensor data have been increasingly prevalent in many domains, such as genetics, social sciences, engineering. In a wide range of applications, unsupervised analysis, in particular the high-order clustering, can be applied to discover the hidden modules in these high-dimensional tensor data. For example, in microbiome studies, microbiome samples are often measured across multiple body sites from multiple subjects \citep{faust2012microbial,flores2014temporal}, resulting in the three-way tensors with subjects, body sites, and bacteria taxa as three modes. It has been reported that multiple microbial taxa can coexist within or across multiple body sites and subjects can form different subpopulations  \citep{faust2012microbial}. Similar data structures can also be found in multi-tissue multi-individual gene expression data \citep{wang2019three}. Mathematically, these patterns correspond to high-order clusters, i.e., the underlying multi-way block structures in the data tensor. We also refer readers to the recent survey \citep{henriques2019triclustering} on high-order clustering in applications. 

In the literature, a number of methods have been proposed for triclustering or high-order clustering of tensor data, such as divide and conquer \citep{li2009effective}, seed growth \citep{sim2010discovering}, stochastic approach \citep{amar2015hierarchical}, exhaustive approaches \citep{jiang2004mining}, pattern-based approach \citep{ji2006mining}, etc. However, the theoretical guarantees for those existing procedures are not well established to our best knowledge. 

This paper aims to fill the void of theory in high-order clustering. Suppose we observe an $n_1 \times \cdots \times n_d$-dimensional order-$d$ tensor $\bcY$ that satisfies
\begin{equation} \label{eq: overall model}
	\bcY = \bcX + \bcZ,
\end{equation}
where $\bcX \in \bbR^{n_1 \times \cdots \times n_d}$ is the underlying signal with planted structure and $\bcZ$ is the noise that has i.i.d standard normal distributed entries. Our goal is to detect or recover the ``planted structure" of the signal $\bcX$. The specific problems in this paper are listed below.

\subsection{Problem Formulations}\label{sec:formulation}

First, we consider the signal tensor $\bcX$ that contains a constant planted structure:
\begin{equation}\label{eq: CHC problem set up}	
\bcX \in \sX_{\CHC}(\bk, \bn, \lambda),\quad \sX_{\CHC}(\bk, \bn, \lambda) = \left\{\lambda' \mathbf{1}_{I_1} \circ \cdots \circ \mathbf{1}_{I_d}: I_i \subseteq [n_i], |I_i| = k_i, \lambda'\geq\lambda\right\}.
\end{equation}
Here, ``$\circ$" denotes the vector outer product, $\mathbf{1}_{I_i}$ is the $n_i$-dimensional indicator vector such that $(\mathbf{1}_{I_i})_j = 1$ if $j\in I_i$ and $(\mathbf{1}_{I_i})_j = 0$ if $j\notin I_i$; $\lambda$ represents the signal strength. We collectively denote $\bk = (k_1, \ldots, k_d)$ and $\bn = (n_1, \ldots, n_d)$ for convenience. The support of the planted structure of $\bcX$ is denoted as $\cS(\bcX):=(I_1, \ldots, I_d)$. We refer to this model \eqref{eq: overall model}\eqref{eq: CHC problem set up} as the \textbf{constant high-order clustering (CHC)}. The constant planted clustering model in tensor or matrix biclustering (BC) data has been considered in a number of recent literature (see, e.g., \cite{butucea2013detection,butucea2015sharp,sun2013maximal,brennan2018reducibility,chi2017convex,cai2017computational,brennan2019universality,kolar2011minimax,chen2016statistical,xia2019sup}).

We also consider a more general setting that $\bcX$ contains a rank-one planted structure: 
\begin{equation}\label{eq: ROHC problem set up}
	\bcX \in \sX_{\ROHC}(\bk, \bn, \mu), \quad \sX_{\ROHC}(\bk, \bn, \mu) = \left\{\mu' \bv_1 \circ \cdots \circ \bv_d: \bv_i \in \cV_{n_i, k_i}, \mu'\geq \mu\right\},
\end{equation}
where
\begin{equation*}
	\cV_{n, k} := \left\{ \bv\in \mathbb{S}^{n-1}: k/C_1 \leq \|\bv\|_0 \leq k \text{ and } k^{-1/2} \leq |\bv_{i}| \leq C_2k^{-1/2} \text{ for } i \in S(\bv) \right\}, \quad C_1, C_2 > 1
\end{equation*}
is the set of all $k$-sparse unit vectors with near-uniform magnitude. Here $S(\bv)$ denotes the support of the vector $\bv$ and its formal definition is given in Section \ref{sec: notations}. Throughout the paper, we refer to the model in \eqref{eq: overall model}\eqref{eq: ROHC problem set up} as the \textbf{rank-one high-order clustering (ROHC)}. Especially if $d=2$, i.e., in the matrix case, this model (rank-one submatrix (ROS)) was considered in \cite{sun2013maximal,busygin2008biclustering,madeira2004biclustering,brennan2018reducibility}. For both models, we hope to answer the following questions on detection ($\cP_D$) and recovery ($\cP_R$):
\begin{enumerate}[leftmargin=*]
    \item[$\cP_D$] {\it When we can \textbf{detect} if any high-order cluster exists and when such conclusion cannot be made.} To be specific, consider the following hypothesis tests:
    \begin{equation}\label{eq: constant subtensor detection}
    \begin{split}
        \CHC_D(\bn, \bk, \lambda):& \quad  H_0 : \bcX = \mathbf{0} \quad \text{v.s.}\quad \quad H_1: \bcX \in \sX_{\CHC}(\bk, \bn, \lambda),\\
        \ROHC_D(\bn, \bk, \mu):& \quad  H_0 : \bcX = \mathbf{0} \quad \text{v.s.}\quad \quad H_1: \bcX \in \sX_{\ROHC}(\bk, \bn, \mu),
    \end{split}
    \end{equation}
we ask when there is a sequence of algorithms that can achieve reliable detection, i.e., both type-I and II errors tend to zero.

    \item[$\cP_R$] {\it How to \textbf{recover} the support of the cluster when it exists.} Specifically, we assume $H_1$ holds and aim to develop an algorithm that recovers the support $S(\bcX)$ based on the observation of $\bcY$. Denote the CHC and ROHC recovery problems considered in this paper as $\CHC_R(\bn, \bk, \lambda)$ and $\ROHC_R(\bn, \bk, \mu)$, respectively. We would like to know when there exists a sequence of algorithms that can achieve reliable recovery, i.e., the probability of correctly recovering $S(\bcX)$ tends to one. 
\end{enumerate}

We study the performance of both {\bf unconstrained-time algorithms} and {\bf polynomial-time algorithms} for both detection $\cP_D$ and recovery $\cP_R$. The class of unconstrained algorithms includes all procedures with unlimited computational resources, while an algorithm that runs in polynomial-time has access to poly($n$) independent random bits and must finish in poly$(n)$ time, where $n$ is the size of input. For convenience of exposition, we assume the explicit expressions can be exactly computed and $N(0,1)$ random variable can be sampled in $O(1)$ time.

\subsection{Main Results}\label{sec:main-result}

In this paper, we give a comprehensive characterization of the statistical and computational limits of the detection and recovery for both $\CHC$ and $\ROHC$ models. Denote $n:= \max_i n_i, k:= \max_i k_i$, and assume $d$ is fixed. For technical convenience, our discussions are based on two asymptotic regimes:
\begin{equation}\tag{A1}\label{assum: asymptotic assumption}
\forall i \in [d], \quad n_i \to \infty, \quad k_i \to \infty \quad \text{and} \quad k_i/n_i \to 0;
\end{equation}
\begin{equation}\tag{A2}\label{eq: comp lower bound asymptotic regime}
\begin{split}
     \text{or} \qquad & \text{for fixed $0\leq \alpha\leq 1, \beta \in \mathbb{R}$ }, n\to \infty, \quad n_1 = \cdots = n_d = \tilde{\Theta}(n), \\ 
     & k = k_1 = \cdots = k_d= \tilde{\Theta}(n^{\alpha}), \quad \lambda = \tilde{\Theta}(n^{-\beta}), \quad \mu/\sqrt{k^d} =\tilde{\Theta}(n^{-\beta}).
\end{split}  
\end{equation} 
In \eqref{eq: comp lower bound asymptotic regime}, $\alpha$ and $\beta$ represent the sparsity level and the signal strength of the cluster, respectively. The cluster becomes sparser as $\alpha$ decreases and the signal becomes stronger as $\beta$ decreases. A rescaling of $\mu$ in \eqref{eq: comp lower bound asymptotic regime} is to make the magnitude of normalized entries in cluster of $\ROHC$ to be approximately one, which enables a valid comparison between the computational hardness of $\CHC$ and $\ROHC$. 

The following informal statements summarize the main results of this paper. 
\begin{Theorem}[Informal: Phase Transitions in CHC]\label{thm:informal-CHC} Define
\begin{equation}\label{eq: CHC phase transition threshold}
\begin{split}
    &\beta^s_{\CHC_D} := \left(d\alpha  - d/2\right) \lor (d-1)\alpha/2, \quad  \beta^s_{\CHC_R}:=(d-1)\alpha/2,\\
    &\beta^c_{\CHC_D} :=  \left(d\alpha - d/2\right) \lor 0, \quad \beta^c_{\CHC_R} := \left((d-1)\alpha -(d-1)/2\right) \lor 0.
\end{split}
\end{equation}
Under the asymptotic regime \eqref{eq: comp lower bound asymptotic regime}, the statistical and computational limits of $\CHC_D(\bk, \bn, \lambda)$ and $\CHC_R(\bk, \bn, \lambda)$ exhibit the following phase transitions:
\begin{itemize}[leftmargin=*]
    \item CHC Detection:
    \begin{enumerate}[label=(\roman*)]
        \item $\beta>\beta^s_{\CHC_D}$: reliable detection is information-theoretically impossible.
        \item $\beta^c_{\CHC_D} <  \beta < \beta^s_{\CHC_D} $: the computational inefficient test $\psi_{\CHC_D}^s$ in Section \ref{sec: CHC and ROHC detection information limit} succeeds, but polynomial-time reliable detection is impossible based on the hypergraphic planted clique (HPC) conjecture (Conjecture \ref{conj: hardness of tensor clique detection}).
        \item $\beta<\beta^c_{\CHC_D}$: the polynomial-time test $\psi_{\CHC_D}^c$ in Section \ref{sec: CHC and ROHC detection poly alg limit} based on combination of sum and max statistics succeeds.
    \end{enumerate}
    \item CHC Recovery:
    \begin{enumerate}[label=(\roman*)]
        \item $\beta > \beta^s_{\CHC_R}$: reliable recovery is information-theoretically impossible.
        \item $\beta^c_{\CHC_D} < \beta < \beta^s_{\CHC_R}$: the exhaustive search (Algorithm \ref{alg: CHC recovery combina search}) succeeds, but polynomial-time reliable recovery is impossible based on HPC conjecture (Conjecture \ref{conj: hardness of tensor clique detection}) and hypergraphic planted dense subgraph (HPDS) recovery conjecture (Conjecture \ref{conj: HPDS recovery conjecture}).
        \item $\beta < \beta^c_{\CHC_D}$: the combination of polynomial-time Algorithms \ref{alg: CHC and ROHC recovery simple thresholding} and \ref{alg: CHC recovery aggregated SVD} succeeds. 
    \end{enumerate}
\end{itemize}
\end{Theorem}

\begin{Theorem}[Informal: Phase Transitions in ROHC]\label{thm:informal-ROHC} Define
\begin{equation}\label{eq: ROHC phase transition threshold}
\begin{split}
  &\beta^s_{\ROHC} = \beta^s_{\ROHC_D}= \beta^s_{\ROHC_R}:= (d-1)\alpha/2,\\ &\beta^c_{\ROHC} =\beta^c_{\ROHC_D} = \beta^c_{\ROHC_R} := \left(\alpha d/2 - d/4 \right) \lor 0.  
\end{split}
\end{equation}
Under the asymptotic regime \eqref{eq: comp lower bound asymptotic regime}, the statistical and computational limits of $\ROHC_D(\bk, \bn, \mu)$ and $\ROHC_R(\bk, \bn, \mu)$ exhibit the following phase transitions:
    \begin{enumerate}[label=(\roman*)]
        \item $\beta>\beta^s_{\ROHC}$: reliable detection and recovery are information-theoretically impossible.
        \item $\beta^c_{\ROHC} <  \beta < \beta^s_{\ROHC} $: the computational inefficient test $\psi_{\ROHC_D}^s$ in Section \ref{sec: CHC and ROHC detection information limit} succeeds in detection and the search Algorithm \ref{alg: ROHC recovery combinatorial search} succeeds in recovery, but polynomial-time reliable detection and recovery are impossible based on the HPC conjecture (Conjecture \ref{conj: hardness of tensor clique detection}).
        \item $\beta<\beta^c_{\ROHC}$: the polynomial-time test $\psi_{\ROHC_D}^c$ in Section \ref{sec: CHC and ROHC detection poly alg limit} succeeds in detection and the combination of polynomial-time Algorithms \ref{alg: CHC and ROHC recovery simple thresholding} and \ref{alg: ROHC and CHC recover via HOOI} succeeds in recovery.
    \end{enumerate}
\end{Theorem}
In Table \ref{tab: phase transition table}, we summarize the statistical and computational limits in Theorems \ref{thm:informal-CHC} and \ref{thm:informal-ROHC} in terms of the original parameters $k,n,\lambda, \mu$ and provide the corresponding algorithms that achieve these limits.
\begin{table}
	\centering
	\begin{tabular}{c | c | c |c}
	\hline
	 & $\CHC_D$ & $\CHC_R$ & $\ROHC_D$ \& $\ROHC_R$\\
	 \hline
	 \hline
	 Impossible & $\lambda^2 \ll \frac{n^d}{k^{2d}} \wedge \frac{1}{k^{d-1}}$ & $\lambda^2 \ll \frac{1}{k^{d-1}}$ & $\frac{\mu^2}{k^d} \ll \frac{1}{k^{d-1}}$\\
	 \hline
	 \hline
	 Hard &  $\frac{n^d}{k^{2d}} \wedge \frac{1}{k^{d-1}} \lesssim \lambda^2 \ll \frac{n^d}{k^{2d}} \wedge 1$ & $ \frac{1}{k^{d-1}} \lesssim \lambda^2 \ll \frac{n^{d-1}}{k^{2(d-1)}} \wedge 1$ & $ \frac{1}{k^{d-1}} \lesssim \frac{\mu^2}{k^d} \ll \frac{n^{d/2}}{k^d} \wedge 1 $\\
	 \hline
	 Algorithms & $\psi_{\CHC_D}^s$ & Alg \ref{alg: CHC recovery combina search} & $\psi_{\ROHC_D}^s$\& Alg \ref{alg: ROHC recovery combinatorial search} \\
	 \hline
	 \hline
	 Easy & $\lambda^2 \gtrsim \frac{n^d}{k^{2d}} \wedge 1$ & $\lambda^2 \gtrsim \frac{n^{d-1}}{k^{2(d-1)}} \wedge 1$ & $\frac{\mu^2}{k^d} \gtrsim \frac{n^{d/2}}{k^d} \wedge 1 $\\
	 \hline
	 Algorithms & $\psi_{\CHC_D}^c$ & Algs \ref{alg: CHC and ROHC recovery simple thresholding} and \ref{alg: CHC recovery aggregated SVD} & $\psi_{\ROHC_D}^c$ \& Algs \ref{alg: CHC and ROHC recovery simple thresholding} and \ref{alg: ROHC and CHC recover via HOOI}\\
	 \hline
	 \hline
	\end{tabular}
	\caption{Phase transition and algorithms for detection and recovery in CHC and ROHC under the asymptotic regime \eqref{eq: comp lower bound asymptotic regime}. Here, easy, hard, and impossible mean polynomial-time solvable, unconstrained-time solvable but polynomial-time unsolvable, and unconstrained-time unsolvable, respectively.} \label{tab: phase transition table}
\end{table}

\begin{figure}[h!]
	\centering
	\subfigure[Constant high-order clustering ($\CHC$)]{
		\includegraphics[height = .4\textwidth]{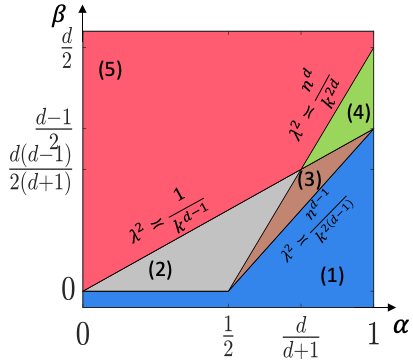}}
	\subfigure[Matrix biclustering ($\BC$)]{
		\includegraphics[height = .4\textwidth]{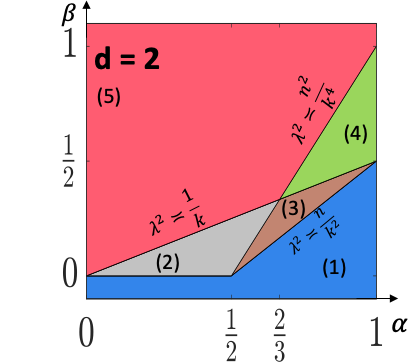}}
	\subfigure[Rank-one high-order clustering ($\ROHC$)]{
		\includegraphics[height = .4\textwidth]{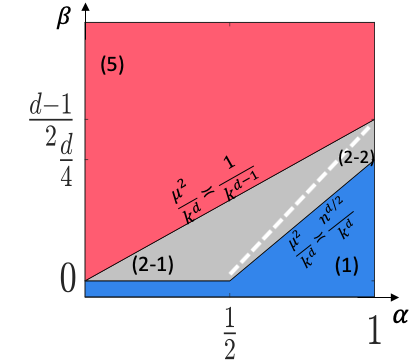}}
	\subfigure[Rank-one submatrix clustering ($\ROS$)]{
		\includegraphics[height = .4\textwidth]{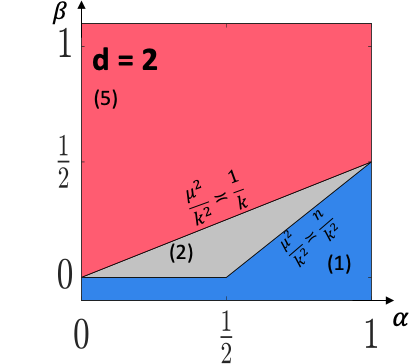}}
	\caption{Statistical and computational phase transition diagrams for constant high-order and rank-one high-order ($d\geq 3$) clustering models (CHC and ROHC) (left two panels) and constant biclustering and rank-one submatrix ($d=2$) clustering models (BC and ROS) (right two panels) under asymptotic regime \eqref{eq: comp lower bound asymptotic regime}. Meaning of each region: (1) all problems detection and recovery both easy; (2),(2-1),(2-2) all problems detection hard and recovery hard; (3) CHC and BC detection easy and recovery hard; (4) CHC and BC detection easy and recovery impossible; (5) all problems detection and recovery impossible.}\label{fig: detection and recovery phase transition}
    \end{figure}
    We also illustrate the phase transition diagrams for both $\CHC, \ROHC$ $(d\geq 3)$ in Figure \ref{fig: detection and recovery phase transition}, Panels (a) and (c). When $d = 2$, the phase transition diagrams in Panels (a) and (c) of Figure \ref{fig: detection and recovery phase transition} reduce to constant biclustering ($\BC$) diagram \citep{ma2015computational,cai2017computational,brennan2018reducibility,chen2016statistical} and rank-one submatrix ($\ROS$) diagram \citep{brennan2018reducibility} in Panels (b) and (d) of Figure \ref{fig: detection and recovery phase transition}. 
    
    \subsection{Comparison with Matrix Clustering and Our Contributions}\label{sec:comparison-matrix}
    The high-order ($d\geq 3$) clustering problems show many distinct aspects from their matrix counterparts ($d=2$). We summarize the differences and highlight our contributions in the aspects of \textit{phase transition diagrams}, \textit{algorithms}, \textit{hardness conjecture}, and \textit{proof techniques} below.

({\bf Phase transition diagrams}) We can see the order-$d$ ($d\geq 3$) tensor clustering has an additional regime: (2-2) in Figure \ref{fig: detection and recovery phase transition} Panel (c). Specifically if $d=2$, $\CHC_R, \ROHC_R$ become $\BC_R, \ROS_R$ that share the same computational limit and there is no gap between the statistical limit and computational efficiency for $\alpha=1$ in $\ROS_R$ (see Panels (b) and (d), Figure \ref{fig: detection and recovery phase transition}). If $d\geq 3$, we need a strictly stronger signal-to-noise ratio to solve $\ROHC_R$ than $\CHC_R$ and there is always a gap between the statistical optimality and computational efficiency for $\ROHC_R$. This difference roots in two level computation barriers, sparsity and tensor structure, in high-order ($d\geq 3$) clustering. To our best knowledge, we are the first to characterize such double computational barriers. 

\begin{figure}[h!]
	\centering
		\includegraphics[width=0.9\textwidth]{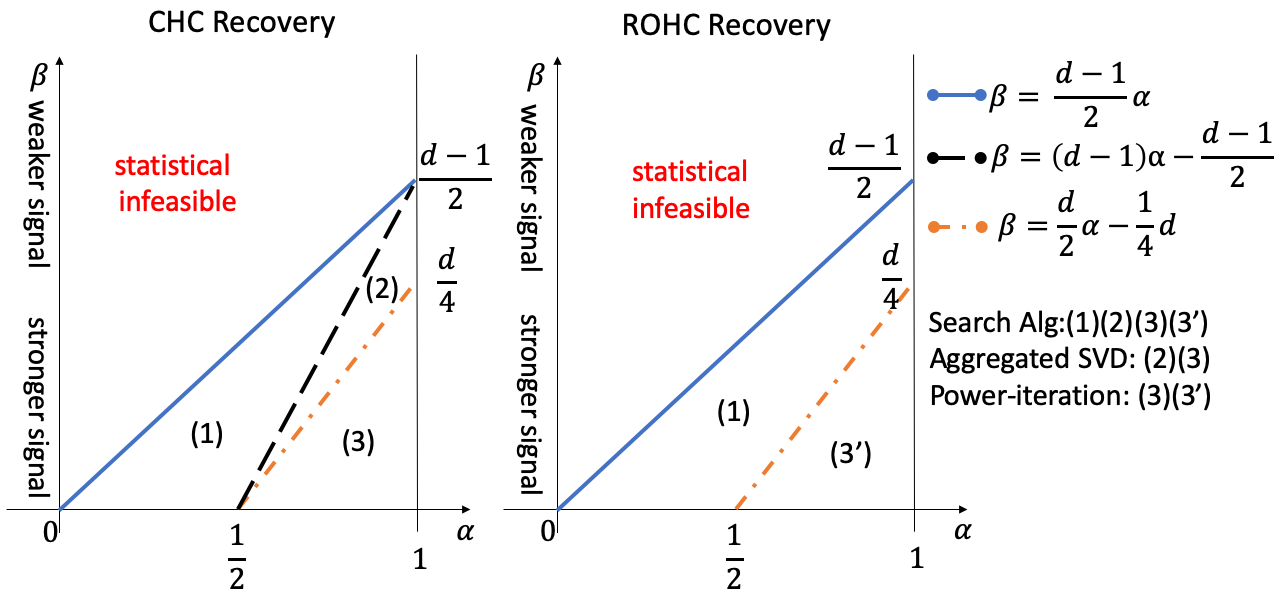}
	\caption{$\CHC_R$ and $\ROHC_R$ diagrams for Exhaustive search, Aggregated-SVD and Power-iteration algorithms under asymptotic regime \eqref{eq: comp lower bound asymptotic regime}. In the right bottom corner, we provide the feasible signal-to-noise ratio regimes for each algorithm. }\label{fig: recovery phase transition for diff algorithm}
\end{figure}

{\bf (Algorithms)} In addition, we develop new algorithms for high-order clustering. For $\CHC_R$ and $\ROHC_R$, we introduce polynomial-time algorithms \textit{Power-iteration} (Algorithm \ref{alg: ROHC and CHC recover via HOOI}), \textit{Aggregated-SVD} (Algorithm \ref{alg: CHC recovery aggregated SVD}), both of which can be viewed as high-order analogues of the matrix spectral clustering. Also, see Section \ref{sec: related work} for a comparison with the methods in the literature. We compare these algorithms and the \textit{exhaustive search} (Algorithms \ref{alg: CHC recovery combina search} and \ref{alg: ROHC recovery combinatorial search}) under the asymptotic regime \eqref{eq: comp lower bound asymptotic regime} in Figure \ref{fig: recovery phase transition for diff algorithm}. Compared to matrix clustering recovery diagram, i.e. Figure \ref{fig: detection and recovery phase transition}(d), a new Regime (2) appears in the high-order ($d\geq 3$) clustering diagram. Different from the matrix clustering, where the polynomial-time spectral method reaches the computational limits for both $\BC_R$ and $\ROS_R$ when $\frac{1}{2}\leq \alpha \leq 1$, the optimal polynomial-time algorithms for $\CHC_R$ and $\ROHC_R$ are distinct: Power-iteration is optimal for $\ROHC_R$ but is suboptimal for $\CHC_R$; the Aggregated-SVD is optimal for $\CHC_R$ but does not apply for $\ROHC_R$. This difference stems from the unique tensor algebraic structure in CHC. 

{\bf (Hardness conjecture)} We adopt the average-case reduction approach to establish the computational lower bounds. It would be ideal to do average-case reduction from the commonly raised conjectures, such as the planted clique (PC) detection or Boolean satisfiability (SAT), so that all of the hardness results of these well-studied conjectures can be inherited to the target problem. However, this route is complicated by the multiway structure in the high-order clustering. Instead, we apply a new average-case reduction scheme from hypergraphic planted clique (HPC) and the hypergraphic planted dense subgraph (HPDS) since HPC and HPDS have a more natural tensor structure that enables a more straightforward average-case reduction. Despite the widely studied planted clique (PC) and planted dense subgraph (PDS) in literature, the HPC and HPDS are far less understood and so are their computational hardness conjectures. The relationship between the computational hardness of PC and HPC remains an open problem \citep{luo2020open}. This paper is among the first to explore the average computational complexities of HPC and HPDS. To provide evidence for the computation hardness conjecture, we show a class of powerful algorithms, including the polynomial-time low-degree polynomials and Metropolis algorithms, are not able to solve HPC detection unless the clique size is sufficiently large. Also, we show low-degree polynomial method only succeeds in HPDS recovery in a restricted parameter regime. These results on HPC and HPDS may be of independent interests in analyzing average-case computational complexity, given the steadily increasing popularity on tensor data analysis recently and the commonly observed statistical-computational gaps therein \citep{richard2014statistical,barak2016noisy,zhang2018tensor,perry2020statistical,hopkins2015tensor,lesieur2017statistical,wein2019kikuchi,dudeja2020statistical}.

{\bf (Proof techniques)} The theoretical analysis in high-order clustering incorporates sparsity, low-rankness, and tensor algebra simultaneously, which is significantly more challenging than its counterpart in matrix clustering. Specifically, to prove the statistical lower bound of $\ROHC_D$, we introduce the new Lemma \ref{lm: MGF bound of power of a sum random walk}, which gives an upper bound for the moment generating function of any power of a symmetric random walk on $\bbZ$ stopped after a hypergeometric distributed number of steps. This lemma is proved by utilizing Hoeffding's inequality and the tail bound integration, which is different from the literature and can be of independent interest. To prove the statistical lower bound of $\CHC_D$, we introduce a new  technique that ``sequentially decompose event" to bound the second moment of the truncated likelihood ratio (see Lemma \ref{lm: bound on the second moment of likelihood ratio}). To prove the computational lower bounds, we introduce new average-case reduction schemes from $\HPC$ and $\HPDS$, including a new reduction technique of {\it tensor reflection cloning} (Algorithm \ref{alg: Tensor Reflecting Clone}). This technique spreads the signal in the planted high-order cluster along each mode evenly, maintains the independence of entries in the tensor, and only mildly reduces the signal magnitude.

\subsection{Related Literature} \label{sec: related work}
This work is related to a wide range of literature on biclustering, tensor decomposition, tensor SVD, and theory of computation. When the order of the observation $d = 2$, the problem \eqref{eq: overall model} reduces to the matrix clustering \citep{ames2011nuclear,butucea2015sharp,chi2017convex,mankad2014biclustering,tanay2002discovering,busygin2008biclustering}. The statistical and computational limits of matrix clustering have been extensively studied in the literature \citep{balakrishnan2011statistical,kolar2011minimax,butucea2013detection,ma2015computational,chen2016statistical,cai2017computational,brennan2018reducibility,brennan2019universality,schramm2020computational}. As discussed in Section \ref{sec:comparison-matrix}, the high-order ($d\geq 3$) tensor clustering exhibits significant differences from the matrix problems in various aspects. 

Another related topic is on tensor decomposition and best low-rank tensor approximation. Although the best low-rank matrix approximation can be efficiently solved by the matrix singular value decomposition (Eckart$-$Young$-$Mirsky Theorem), the best low-rank tensor approximation is NP-hard to calculate in general \citep{hillar2013most}. Various polynomial-time algorithms, which can be seen as the polynomial-time relaxations of the best low-rank tensor approximation, have been proposed in the literature, including the Newton method \citep{zhang2001rank}, alternating minimization \citep{zhang2001rank,richard2014statistical}, high-order singular value decomposition \citep{de2000multilinear}, high-order orthogonal iteration \citep{de2000best}, $k$-means power iteration \citep{anandkumar2014guaranteed,sun2017provable}, sparse high-order singular value decomposition \citep{zhang2019optimal-statsvd}, regularized gradient descent \citep{han2020optimal}, etc. The readers are referred to surveys \cite{kolda2009tensor,cichocki2015tensor}. Departing from most of these previous results, the high-order clustering considered this paper involves both sparsity and low-rankness structures, which requires new methods and theoretical analysis as discussed in Section \ref{sec:comparison-matrix}. 

Our work is also related to a line of literature on average-case computational hardness and the statistical and computational trade-offs. The average-case reduction approach has been commonly used to show computational lower bounds for many recent high-dimensional problems, such as testing $k$-wise independence \citep{alon2007testing}, biclustering \citep{ma2015computational, cai2017computational, cai2020statistical}, community detection \citep{hajek2015computational}, RIP certification \citep{wang2016average, koiran2014hidden}, matrix completion \citep{chen2015incoherence}, sparse PCA \citep{berthet2013complexity,berthet2013optimal,brennan2018reducibility, gao2017sparse,wang2016statistical,brennan2019optimal}, universal submatrix detection \citep{brennan2019universality}, sparse mixture and robust estimation \citep{brennan2019average}, a financial model with asymmetry information \citep{arora2011computational}, finding dense common subgraphs \citep{charikar2018finding}, graph logistic regression \citep{berthet2020statistical}, online local learning \citep{awasthi2015label}. See also a web of average-case reduction to a number of problems in \cite{brennan2018reducibility,brennan2020reducibility} and a recent survey \citep{wu2018statistical}. The average-case reduction is delicate, requiring that a distribution over instances in a conjecturally hard problem be mapped precisely to the target distribution. For this reason, many recent literature turn to show computational hardness results under the restricted models of computation, such as sum of squares \citep{ma2015sum,hopkins2017power,barak2019nearly}, statistical query \citep{feldman2017statistical,diakonikolas2017statistical,diakonikolas2019efficient,feldman2018complexity,wang2015sharp,fan2018curse,kannan2017hidden}, class of circuit \citep{rossman2008constant,rossman2014monotone}, convex relaxation \citep{chandrasekaran2013computational}, local algorithms \citep{gamarnik2014limits}, meta-algorithms based on low-degree polynomials \citep{hopkins2017bayesian,kunisky2019notes} and others. As discussed in Section \ref{sec:comparison-matrix}, this paper is among the first to investigate the hypergraphic planted clique (HPC) and hypergraphic planted dense subgraph (HPDS) problems and their computational hardness. We perform new average-case reduction scheme from these conjectures and develop the computational lower bounds for CHC and ROHC.

\subsection{Organization}\label{sec:organization}

The rest of this article is organized as follows. After a brief introduction of notation and preliminaries in Section \ref{sec: notations}, the statistical limits of high-order cluster recovery and detection are given in Sections \ref{sec: CHC and ROHC recovery} and \ref{sec: CHC and ROHC detection}, respectively. In Section \ref{sec: Comp lower bound for CHC and ROHC}, we establish the computational limits of high-order clustering, along with the hypergraphic planted clique (HPC) and hypergraphic planted dense subgraph (HPDS) models, computational hardness conjectures, and evidence. Discussion and future work are given in Section \ref{sec: discussion}. The technical proofs are collected in supplementary materials \cite{luo2020tensorsupp}. 

\section{Notation and Definitions}\label{sec: notations}
The following notation will be used throughout this article. For any non-negative integer $n$, let $[n] = \{1, \ldots, n\}$. The lowercase letters (e.g., $a, b$), lowercase boldface letters (e.g., $\bu, \bv$), uppercase boldface letters (e.g., $\A, \U$), and boldface calligraphic letters (e.g., $\bcA, \bcX$) are used to denote scalars, vectors, matrices, and order-3-or-higher tensors respectively. For any two series of numbers, say $\{a_n\}$ and $\{b_n\}$, denote $a \asymp b$ if there exist uniform constants $c, C>0$ such that $ca_n \leq b_n\leq Ca_n, \forall n$; and $a  = \Omega (b)$ if there exists uniform constant $c > 0$ such that $a_n \geq c b_n, \forall n$. The notation $a = \tilde{\Theta}(b)$ and $a \gg b$ mean $\lim_{n \to \infty} a_n / n = \lim_{n \to \infty} b_n / n$ and $\lim_{n \to \infty} \log(a_n) / \log n > \lim_{n \to \infty} \log(b_n) / \log n$, respectively. $a \lesssim b$ means $a \leq b$ up to polylogarithmic factors in $n$. We use bracket subscripts to denote sub-vectors, sub-matrices, and sub-tensors. For example, $\bv_{[2:r]}$ is the vector with the $2$nd to $r$th entries of $\bv$; $\D_{[(r+1):n_1, :]}$ contains the $(r+1)$-th to the $n_1$-th rows of $\D$; $\bcA_{[1:s_1, 1:s_2, 1:s_3]}$ is the $s_1$-by-$s_2$-by-$s_3$ sub-tensor of $\bcA$ with index set $\{(i_1, i_2, i_3): 1 \leq i_1\leq s_1, 1\leq i_2 \leq s_2, 1\leq i_3\leq s_3\}$. For any vector $\bv \in \mathbb{R}^{n_1}$, define its $\ell_2$ norm as $\|\bv \|_2 = \left(\sum_i |\bv_{i}|^2\right)^{1/2}$ and $\|\bv\|_0$ is defined to be the number of non-zero entries in $\bv$. Given vectors $
\{\bv_i\}_{i=1}^d \in \bbR^{n_i}$, the outer product $\bcA \in \bbR^{n_1 \times \cdots \times n_d} = \bv_1 \circ \cdots \circ \bv_d$ is defined such that $\bcA_{[i_1, \ldots, i_d]} = (\bv_1)_{i_1} \cdots (\bv_d)_{i_d}$. For any event $A$, let $\bbP(A)$ be the probability that $A$ occurs. 

For any order-$d$ tensor $\bcA \in \mathbb{R}^{n_1\times \cdots\times n_d}$. The matricization $\mathcal{M}(\cdot)$ is the operation that unfolds or flattens the order-$d$ tensor $\bcA\in\mathbb{R}^{n_1\times \cdots \times n_d}$ into the matrix $\mathcal{M}_z(\bcA)\in \mathbb{R}^{n_z\times \prod_{j\neq z}n_j}$ for $z=1,\ldots, d$. Specifically, the mode-$z$ matricization of $\bcA$ is formally defined as 
\begin{equation*}
     \bcA_{[i_1,\ldots, i_d]} = \left(\mathcal{M}_z(\bcA)\right)_{\left[i_z, j\right]}, \quad j = 1 + \sum_{\substack{l=1\\l\neq z}}^d\left\{(i_l-1)\prod_{\substack{m=1\\m\neq z}}^{l-1}n_m\right\}
\end{equation*}
for any $1\leq i_l \leq n_l, l=1,\ldots, d$. Also see \cite[Section 2.4]{kolda2009tensor} for more discussions on tensor matricizations. The mode-$z$ product of $\bcA \in \mathbb{R}^{n_1 \times \cdots \times n_d}$ with a matrix $\U\in \mathbb{R}^{k_z\times n_z}$ is denoted by $\bcA \times_z \U$ and is of size $n_1 \times \cdots \times n_{z-1}\times k_z \times n_{z+1}\times \cdots \times n_d$, such that 
$$(\bcA \times_z \U)_{[i_1, \ldots, i_{z-1}, j, i_{z+1}, \ldots, i_d]} = \sum_{i_z=1}^{n_z} \bcA_{[i_1, i_2, \ldots, i_d]} \U_{[j,i_z]}.$$ For any two distinct $k_1, k_2 \in [d]\, (k_1 < k_2)$ and $j_1 \in [n_{k_1}]$ and $j_2 \in [n_{k_2}]$, we denote 
$$\bcA^{(k_1, k_2)}_{j_1, j_2} = \bcA_{\Big[: , ~ \ldots ~ , : , \underbrace{j_1}_{k_1\text{th index}}, : , ~\ldots, ~ :, \underbrace{j_2}_{k_2\text{th index}},  : ~ \ldots, ~ : \Big]}\in \bbR^{n_1 \times \cdots \times n_{k_1-1} \times n_{k_1 + 1} \times \cdots \times n_{k_2-1} \times n_{k_2+1} \times \cdots \times n_d }$$ 
as a subtensor of $\bcA$. The support of an order-$d$ tensor $\bcX \in \bbR^{n_1 \times \cdots \times n_d}$ is denoted by $S(\bcX) := (I_1,\cdots,I_d)$ where $I_j \in \bbR^{n_j}$ and $(I_j)_i$ equals to zero when $\cM_j(\bcX)_{[i,:]}$ is a zero vector and equals to one otherwise. In particular, when the tensor order is $1$, we simply have the support of a vector $\bv$ is $S(\bv) =  \{j: \bv_{j} \neq 0\}$.

Given a distribution $Q$, let $Q^{\otimes n}$ be the distribution of $(X_1, \ldots, X_n)$ if $\{X_i\}_{i=1}^n$ are i.i.d. copies of $Q$. Similarly, let $Q^{\otimes m \times n}$ and $Q^{\otimes (n^{\otimes d})}$ denote the distribution on $\bbR^{m \times n}$ and $\bbR^{n^{\otimes d}}$ with i.i.d. entries distributed as $Q$. Here $n^{\otimes d} := n \times n \times \cdots \times n$ denotes the order-$d$ Cartesian product. In addition, we use $C, C_1, C_2, c$ and other variations to represent the large and small constants, whose actual values may vary from line to line.


Next, we formally define the statistical and computational risks to quantify the fundamental limits of high-order clustering. First, we define the risk of testing procedure $\phi_D(\bcY)\in\{0,1\}$ as the sum of Type-I and Type-II errors  for detection problems $\CHC_D$ and $\ROHC_D$:
\begin{equation*}
\cE_{\cP_D}(\phi_D) = \bbP_0 (\phi_D (\bcY) = 1) + \sup_{\substack{\bcX \in \sX_{\CHC}(\bk, \bn, \lambda)\\ (\text{or } \bcX \in \sX_{\ROHC}(\bk, \bn, \mu))}} \bbP_{\bcX}(\phi_D(\bcY) = 0),
\end{equation*}
where $\bbP_0$ is the probability under $H_0$ and $\bbP_{\bcX}$ is the probability under $H_1$ with the signal tensor $\bcX$. We say $\{\phi_D\}_n$ reliably detect in $\cP_D$ if $\lim_{n\to\infty}\cE_{\mathcal{P}_D}(\phi_D) = 0$. Second, for recovery problem $\CHC_R$ and $\ROHC_R$, define the recovery error for any recovery algorithm $\phi_R(\bcY) \in \{(I_1, \ldots, I_d): I_i \subseteq \{1,\ldots, n_i\}\}$ as
\begin{equation*}
	\cE_{\cP_R}(\phi_R) = \sup_{\bcX \in \sX_{\CHC}(\bk, \bn, \lambda)} \bbP_{\bcX} (\phi_R (\bcY) \neq S(\bcX)) \quad \text{or}\quad \sup_{\bcX \in \sX_{\ROHC}(\bk, \bn, \mu)} \bbP_{\bcX} (\phi_R (\bcY) \neq S(\bcX)).
\end{equation*} 
We say $\{\phi_R\}_n$ reliably recover in $\cP_R$ if $\lim_{n\to\infty}\cE_{\cP_R}(\phi_R)= 0$. Third, denote $\AllAlg^{D}$, $\AllAlg^{R}$,  $\PolyAlg^{D}$, $\PolyAlg^{R}$ as the collections of unconstrained-time algorithms and polynomial-time algorithms for detection and recover problems, respectively. Then we can define four different statistical and computational risks as follows,
\begin{equation*}
	\cE^s_{\cP_D} := \inf_{\phi_D \in \AllAlg^{D}} \cE_{\cP_D}(\phi_D), \quad  \cE^c_{\cP_D} := \inf_{\phi_D \in \PolyAlg^{D}} \cE_{\cP_D}(\phi_D),
\end{equation*}
\begin{equation*}
	\cE^s_{\cP_R} := \inf_{\phi_R \in \AllAlg^{R}} \cE_{\cP_R}(\phi_R), \quad \cE^c_{\cP_R} := \inf_{\phi_R \in \PolyAlg^{R}} \cE_{\cP_R}(\phi_R).
\end{equation*}

\section{High-order Cluster Recovery: Statistical Limits and Polynomial-time Algorithms} \label{sec: CHC and ROHC recovery}

This section studies the statistical limits of high-order cluster recovery. We first present the statistical lower bounds of $\lambda$ and $\mu$ that guarantee reliable recovery, then we give unconstrained-time algorithms that achieves these lower bounds. We also propose computationally efficient algorithms, Thresholding Algorithm, Power-iteration, and Aggregated-SVD, with theoretical guarantees.

\subsection{$\CHC_R$ and $\ROHC_R$: Statistical Limits} \label{sec: CHC and ROHC recovery stat limit}

Recall \eqref{eq: CHC phase transition threshold} and \eqref{eq: ROHC phase transition threshold}, we first present the statistical lower bounds for reliable recovery of $\CHC_R$ and $\ROHC_R$.
\begin{Theorem}[Statistical Lower Bounds for $\CHC_R$ and $\ROHC_R$]\label{thm: CHC and ROHC recovery stat lower bound}
	Consider $\CHC_{R}(\bk, \bn, \lambda)$ and $\ROHC_R(\bk, \bn, \mu)$. Let $0 < \eta < \frac{1}{8}$ be fixed. Under the asymptotic regime \eqref{assum: asymptotic assumption}, if 
	\begin{equation*}
		\lambda \leq \max\left( \left\{ \sqrt{\frac{\eta \log(n_i - k_i)}{\prod_{z=1, z\neq i}^d k_z}} \right\}_{i=1}^d \right) \quad  \left(\text{ or }\quad \frac{\mu}{\sqrt{\prodk}} \leq \max\left( \left\{ \sqrt{\frac{\eta \log(n_i - k_i)}{\prod_{z=1, z\neq i}^d k_z}} \right\}_{i=1}^d \right) \right),
	\end{equation*}
	we have
	\begin{equation*}
	\cE_{\CHC_R}^s (\text{ or } \cE^s_{\ROHC_R}) \geq \frac{\sqrt{M}}{1+ \sqrt{M}} \left(1- 2\eta - \frac{2\eta}{\log M}\right) \to 1 - 2\eta,
	\end{equation*}
	where $M = \max( \{n_i - k_i\}_{i=1}^d )$. Moreover, under the asymptotic regime \eqref{eq: comp lower bound asymptotic regime}, if $\beta > \beta^s_{\CHC_R}$ (or $\beta^s_{\ROHC_R}$), we have $\cE^s_{\CHC_R} (\text{ or } \cE^s_{\ROHC_R}) \to 1-2\eta$.
\end{Theorem}

We further propose the $\CHC_R$ Search (Algorithm \ref{alg: CHC recovery combina search}) and $\ROHC_R$ Search (Algorithm \ref{alg: ROHC recovery combinatorial search}) with the following theoretical guarantees. These algorithms exhaustively search all possible cluster positions and find one that best matches the data. In particular, Algorithm \ref{alg: CHC recovery combina search} is exactly the maximum likelihood estimator. It is note worthy in Algorithm \ref{alg: ROHC recovery combinatorial search}, we generate $\bcZ_1$ with i.i.d. standard Gaussian entries and construct $\bcA = \frac{\bcY + \bcZ_1}{\sqrt{2}}$ and $\bcB = \frac{\bcY - \bcZ_1}{\sqrt{2}}$. In that case, $\bcA$ and $\bcB$ becomes two independent sample tensors, which facilitate the theoretical analysis. Such a scheme is mainly for technical convenience and not necessary in practice.

\begin{Theorem}[Guarantee of $\CHC_R$ Search] \label{thm: CHC recovery stat upper bound}
	Consider $\CHC_R(\bk, \bn, \lambda)$ under the asymptotic regime \eqref{assum: asymptotic assumption}. There exists $C_0 > 0$ such that when $\lambda \geq C_0  \sqrt{\frac{\sum_{i=1}^d \log (n_i - k_i) }{\underset{1\leq i \leq d}\min \{ \prod_{z=1, z\neq i}^d k_z \} }}$,
	Algorithm \ref{alg: CHC recovery combina search} identifies the true support of $\bcX$ with probability at least $1 - C\sum_{i=1}^d (n_i - k_i)^{-c}$ for some $c, C > 0$. Moreover, under the asymptotic regime \eqref{eq: comp lower bound asymptotic regime}, Algorithm \ref{alg: CHC recovery combina search} achieves the reliable recovery of $\CHC_R$ when $\beta < \beta^s_{\CHC_R}$.
\end{Theorem}

\begin{Theorem}[Guarantee of $\ROHC_R$ Search] \label{thm: ROHC recovery stat upper bound}
	Consider $\ROHC_R(\bk, \bn, \mu)$ under the asymptotic regime \eqref{assum: asymptotic assumption}. There is an absolute constant $C_0 > 0$ such that if $\mu \geq C_0 \sqrt{k  \log n}$, then Algorithm \ref{alg: ROHC recovery combinatorial search} identifies the true support of $\bcX$ with probability at least $1 - C\sum_{i=1}^d (n_i - k_i)^{-1}$ for some constant $C > 0$. Moreover, under the asymptotic regime \eqref{eq: comp lower bound asymptotic regime}, Algorithm \ref{alg: ROHC recovery combinatorial search} achieves the reliable recovery of $\ROHC_R$ when $\beta < \beta^s_{\ROHC_R}$.
\end{Theorem}

\begin{algorithm}[h] \caption{$\CHC_R$ Search} \label{alg: CHC recovery combina search}
\begin{algorithmic}[1]
  \State \textbf{Input}: $\bcY \in \bbR^{n_1 \times \cdots \times n_d}$, sparsity level $\bk = (k_1, \ldots, k_d)$.
\State \textbf{Output}: \begin{equation*}
			(\hat{I}_1, \ldots, \hat{I}_d) = \arg \max_{\substack{I_i \subseteq [n_i], |I_i| = k_i\\ i = 1, \ldots, d}} \sum_{i_1 \in I_1} \ldots \sum_{i_d \in I_d} \bcY_{[i_1, \ldots, i_d]}.
		\end{equation*}  
\end{algorithmic}

\end{algorithm}

\begin{algorithm}[h] \caption{$\ROHC_R$ Search} \label{alg: ROHC recovery combinatorial search}
	\begin{algorithmic}[1]
	\State \textbf{Input:} $\bcY\in \bbR^{n_1 \times \cdots \times n_d}$, sparsity upper bound $\bk$.
		\State Sample $\bcZ_1 \sim N(0,1)^{\otimes n_1 \times \cdots \times n_d}$ and construct $\bcA = \frac{\bcY + \bcZ_1}{\sqrt{2}}$ and $\bcB = \frac{\bcY - \bcZ_1}{\sqrt{2}}$.
		\State For each $\{\bar{k}_i \}_{i=1}^d$ satisfying $\bar{k}_i \in [1, k_i]\, (1 \leq i \leq d)$ do:
		\begin{enumerate}[label=(\alph*),leftmargin=*]
		    \item  Compute
		    \begin{equation*}
		        (\hat{\bu}_1, \ldots, \hat{\bu}_d) = \arg \max_{(\bu_1, \ldots, \bu_d) \in S_{\bar{k}_1}^{n_1} \times \cdots \times S_{\bar{k}_d}^{n_d}} \bcA \times_1 \bu_1^\top \times \cdots \times_d \bu_d^\top.
		    \end{equation*}
		    Here, $S_{t}^n$ is the set of vectors $\bu \in \{-1, 1, 0\}^n$ with exactly $t$ nonzero entries.
		    \item For each tuple $(\hat{\bu}_1, \ldots, \hat{\bu}_d)$ computed from Step (a), mark it if it satisfies \begin{equation*}
		        \left\{j: \left(\bcB \times_1 \hat{\bu}_1^\top \times \cdots \times_{i-1} \hat{\bu}_{i-1}^\top \times_{i+1} \hat{\bu}_{i+1}^\top \times \cdots \times_d \hat{\bu}_d^\top \right)_{j} (\hat{\bu}_i)_{j} \geq \frac{1}{2\sqrt{2}} \frac{\mu}{\sqrt{\prod_{i=1}^d k_i}  } \prod_{z\neq i} \bar{k}_z\right\}
		    \end{equation*} is exactly the support of $\hat{\bu}_i$, $S(\hat{\bu}_i)$ for all $1 \leq i \leq d$.
		\end{enumerate}
		\State Among all marked tuples $(\hat{\bu}_1, \ldots, \hat{\bu}_d)$, we find the one, say $(\tilde{\bu}_1, \ldots, \tilde{\bu}_d)$, that maximizes $\sum_{i=1}^d |S(\hat{\bu}_i)|$. 
		\State {\bf Output:} $\hat{I}_i = S(\tilde{\bu}_i)\, (1 \leq i \leq d)$.
	\end{algorithmic}
\end{algorithm}

\begin{algorithm}[h] \caption{$\CHC_R$ and $\ROHC_R$ Thresholding Algorithm} \label{alg: CHC and ROHC recovery simple thresholding}
\begin{algorithmic}[1]
   \State \textbf{Input}: $\bcY \in \bbR^{n_1 \times \cdots \times n_d}$.
\State \textbf{Output}: \begin{equation*}
			(\hat{I}_1, \ldots , \hat{I}_d) = \{(i_1, \ldots, i_d): |\bcY_{[i_1, \ldots, i_d]}| \geq \sqrt{2(d+1)\log n} \}.
		\end{equation*} 
\end{algorithmic}
\end{algorithm}

Combining Theorems \ref{thm: CHC and ROHC recovery stat lower bound}, \ref{thm: CHC recovery stat upper bound}, and \ref{thm: ROHC recovery stat upper bound}, we can see if $k_1 \asymp k_2 \asymp \cdots \asymp k_d$, Algorithms \ref{alg: CHC recovery combina search}, \ref{alg: ROHC recovery combinatorial search} achieve the minimax statistical lower bounds for $\CHC_R$, $\ROHC_R$. On the other hand, Algorithms \ref{alg: CHC recovery combina search} and \ref{alg: ROHC recovery combinatorial search} are based on computationally inefficient exhaustive search. Next, we introduce the polynomial-time algorithms.

\subsection{$\CHC_R$ and $\ROHC_R$: Polynomial-time Algorithms}

The polynomial-time algorithms for solving $\CHC_R$ and $\ROHC_R$ rely on the sparsity level $k_i\, (1 \leq i \leq d)$. First, when $k \lesssim \sqrt{n}$ (sparse regime), we propose Thresholding Algorithm (Algorithm \ref{alg: CHC and ROHC recovery simple thresholding}) that selects the high-order cluster based on the largest entry in absolute value from each tensor slice. The theoretical guarantee of this algorithm is given in Theorem \ref{thm: CHC and ROHC recovery thresholding}. 

\begin{Theorem}[Guarantee of Thresholding Algorithm for $\CHC_R$ and $\ROHC_R$]\label{thm: CHC and ROHC recovery thresholding} 
Consider $\CHC_R(\bk, \bn, \lambda)$ and $\ROHC_R(\bk, \bn, \mu)$. If $\lambda \geq 2\sqrt{2(d+1) \log n}$ (or $\mu/\sqrt{\prodk} \geq 2\sqrt{2(d+1) \log n}$),
	Algorithm \ref{alg: CHC and ROHC recovery simple thresholding} exactly recovers the true support of $\bcX$ with probability at least $1 - O(n ^{-1})$. Moreover, under the asymptotic regime \eqref{eq: comp lower bound asymptotic regime}, Algorithm \ref{alg: CHC and ROHC recovery simple thresholding} achieves the reliable recovery of $\CHC_R$ and $\ROHC_R$ when $\beta < 0$.
\end{Theorem}

Second, when $k \gtrsim \sqrt{n}$ (dense regime), we consider the Power-iteration given in Algorithm \ref{alg: ROHC and CHC recover via HOOI}, which is a modification of the tensor PCA methods in the literature \citep{richard2014statistical,anandkumar2014guaranteed,zhang2018tensor} and can be seen as an tensor analogue of the matrix spectral clustering method.

\begin{algorithm}[!h] \caption{Power-iteration for $\CHC_R$ and $\ROHC_R$} \label{alg: ROHC and CHC recover via HOOI}
	\begin{algorithmic}[1]
	\State \textbf{Input:} $\bcY\in \bbR^{n_1 \times \cdots \times n_d}$, maximum number of iterations $t_{max}$.
	\State Sample $\bcZ_1 \sim N(0,1)^{\otimes n_1 \times \cdots \times n_d}$ and construct $\bcA = (\bcY + \bcZ_1)/\sqrt{2}$ and $\bcB = (\bcY - \bcZ_1)/\sqrt{2}$.
		\State (Initiation) Set $t = 0$. For $ i = 1: d$, compute the top left singular vector of $\cM_i(\bcA)$ and denote it as $\hat{\bu}_i^{(0)}$.
		\State For $t=1, \ldots, t_{\max}$, do
        \begin{enumerate}[label=(\alph*)]
        \item For $i = 1$ to $d$, update
		\begin{equation*}
			\hat{\bu}_i^{(t)} = \NORM( \bcA \times_1 ( \hat{\bu}_1^{(t)})^\top \times \cdots \times_{i-1} (\hat{\bu}_{i-1}^{(t)})^\top \times_{i+1} (\hat{\bu}_{i+1}^{(t-1)})^\top \times \cdots \times_d (\hat{\bu}_d^{(t-1)})^\top).
		\end{equation*}
		Here, $\NORM(\bv) = \bv/\|\bv\|_2$ is the normalization of vector $\bv$.
		\end{enumerate}
		\State Let $(\hat{\bu}_1, \ldots, \hat{\bu}_d) := ( \hat{\bu}_1^{(t_{\max})}, \ldots, \hat{\bu}_d^{(t_{\max})} )$. For $i = 1,\ldots, d$, calculate 
	\begin{equation}
		\bv_i := \bcB \times_1 \hat{\bu}_1^\top \times \cdots \times_{i-1} \hat{\bu}_{i-1}^\top \times_{i+1} \hat{\bu}_{i+1}^\top  \times \cdots \times_{d} \hat{\bu}_d^\top  \in \mathbb{R}^{n_i}.
	\end{equation}
	\begin{itemize}
	    \item If the problem is $\CHC_R$, the component values of $\bv_i$ form two clusters. Sort $\{ (\bv_i)_{j} \}_{j=1}^{n_i}$ and cut the sequence at the largest gap between the consecutive values. Let the index subsets of two parts be $\hat{I}_i$ and $[n_i] \setminus \hat{I}_i$. 
	    Output: $\hat{I}_i$
	    \item If the problem is $\ROHC_R$, the component values of $\bv_i$ form three clusters. Sort the sequence $\{ (\bv_i)_{j} \}_{j=1}^{n_i}$, cut at the two largest gaps between the consecutive values, and form three parts. Among the three parts, pick the two smaller-sized ones, and let the index subsets of these two parts be $\hat{I}^1_i, \hat{I}^2_i$. Output: $\hat{I}_i = \hat{I}^1_i \bigcup \hat{I}^2_i$
	\end{itemize}
	\State \textbf{Output:} $ \{\hat{I}_i\}_{i=1}^d$.
	\end{algorithmic}
\end{algorithm}

We also propose another polynomial-time algorithm, \emph{Aggregated SVD}, in Algorithm \ref{alg: CHC recovery aggregated SVD} for the dense regime of $\CHC_R$. As its name suggests, the central idea is to first transform the tensor $\bcY$ into a matrix by taking average, then apply matrix SVD. Aggregated-SVD is in a similar vein of the hypergraph adjacency matrix construction in the hypergraph community recovery literature \citep{ghoshdastidar2017consistency,kim2017community}. 

\begin{algorithm}[h] \caption{Aggregated-SVD for $\CHC_R$} \label{alg: CHC recovery aggregated SVD}
	\begin{algorithmic}[1]
	\State \textbf{Input:} $\bcY\in \bbR^{n_1 \times \cdots \times n_d}$.
		\State \textbf{For} $i = 1, 2, \ldots, d$, do:
		\begin{enumerate}[label=(\alph*)]
		    \item Find $i^* = \arg \min_{j \neq i} n_j$ and calculate $\Y^{(i, i^*)} \in \bbR^{n_{i} \times n_{i^*}}$ where $\Y^{(i, i^*)}_{[k_1, k_2]} := \SUM (\bcY^{(i, i^*)}_{k_1, k_2})/ \sqrt{\prod_{j=1, j\neq i, i^*}^d n_j}$ for $ 1 \leq k_1 \leq n_{i}, 1 \leq k_2 \leq n_{i^*}$.
		Here $\SUM(\bcA) := \sum_{i_1} \cdots \sum_{i_d} \bcA_{[i_1, \ldots, i_d]}$ 
		and $\bcY_{k_1,k_2}^{(i,i^*)}$ is the subtensor of $\bcY$ defined in Section \ref{sec: notations}.
		\item Sample $\Z_1 \sim N(0,1)^{\otimes n_{i} \times n_{i^*} }$ and form $\A^{(i, i^*)} = (\Y^{(i, i^*)} + \Z_1)/\sqrt{2}$ and $\B^{(i, i^*)} = (\Y^{(i, i^*)} - \Z_1)/\sqrt{2}$. Compute the top right singular vector of $\A^{(i, i^*)}$, denote it as $\bv$.
		\item To compute $\hat{I}_{i}$, calculate $\left(\B^{(i, i^*)}_{[j,:]} \cdot \bv \right)$ for $1 \leq j \leq n_{j}$. These values form two data driven clusters and a cut at the largest gap at the ordered value of $\left\{  \B^{(i, i^*)}_{[j,:]} \cdot \bv \right\}_{j=1}^{n_{i}}$ returns the set $\hat{I}_{i}$ and $[n_{i}] \setminus \hat{I}_{i}$.
		\end{enumerate}
		\State \textbf{Output:} $ \{\hat{I}_i\}_{i=1}^d$.
	\end{algorithmic}
\end{algorithm}

We give guarantees of Power-iteration and Aggregated-SVD for high-order cluster recovery. In particular, Aggregated SVD achieves strictly better performance than Power-iteration in $\CHC_R$, but does not apply for $\ROHC_R$. 

\begin{Theorem}[Guarantee of Power-iteration for $\CHC_R$ and $\ROHC_R$]\label{thm: CHC and ROHC recovery HOOI upper bound} Consider $\CHC_R(\bk, \bn, \lambda)$ and $\ROHC_R(\bk, \bn, \mu)$. Assume $n_i \geq c_0 n\,(1 \leq i \leq d)$ for constant $c_0>0$ where $n := \max_i n_i$. Under the asymptotic regime \eqref{assum: asymptotic assumption}, there exists a uniform constant $C_0>0$ such that if $\lambda\sqrt{\prodk} \geq C_0 n^{\frac{d}{4}}$ (or $\mu \geq C_0 n^{\frac{d}{4}})$ 
$$\text{and}\quad t_{\max} \geq C\log\left(\frac{n}{\lambda\sqrt{\prodk}}\right) \lor C \quad  \left(\text{ or } t_{\max} \geq C(\log(n/\mu) \lor 1)\right),$$
	Algorithm \ref{alg: ROHC and CHC recover via HOOI} identifies the true support of $\bcX$ with probability at least $1 - \sum_{i=1}^d n_i^{-c} - C \exp(-c n)$ for constants $c, C > 0$. Moreover, under the asymptotic regime \eqref{eq: comp lower bound asymptotic regime}, Algorithm \ref{alg: ROHC and CHC recover via HOOI} achieves the reliable recovery of $\CHC_R$ and $\ROHC_R$ when $\beta < (\alpha - 1/2)d/2$. 
\end{Theorem}

\begin{Theorem}[Guarantee of Aggregated-SVD for $\CHC_R$] \label{thm: CHC recovery aggregated svd} Consider $\CHC_R(\bk, \bn, \lambda)$ and Algorithm \ref{alg: CHC recovery aggregated SVD}. There exists a uniform constant $C_0 > 0$ such that if
	\begin{equation} \label{cond: CHC recovery aggregated SVD condition}
\lambda  \geq C_0  \frac{k\sqrt{\prodn}}{\sqrt{n_{\min}}\prodk} \left( 1  +  \sqrt{\frac{ k \log n }{ n_{\min} }} \right),
	\end{equation} 
	the support recovery algorithm based on Aggregated-SVD identifies the true support of $\bcX$ with probability at least $1 - \sum_{i=1}^d n_i^{-c} - C \exp(-c n_{\min})$. Here, $n_{\min} = \min(n_1, \ldots, n_d)$. Moreover, under the asymptotic regime \eqref{eq: comp lower bound asymptotic regime}, Aggregated-SVD achieves reliable recovery of $\CHC_R$ when $\beta < (\alpha -1/2)(d-1)$. 
\end{Theorem}

Combining Theorems \ref{thm: CHC and ROHC recovery thresholding}--\ref{thm: CHC recovery aggregated svd}, we can see the reliable recovery of $\CHC_R$ and $\ROHC_R$ is polynomial-time possible if $\beta < \beta^c_{\CHC_R}: = (\alpha -1/2)(d-1) \vee 0$ and $\beta < \beta^c_{\ROHC_R}:= (\alpha - 1/2)d/2 \vee 0.$ Since $\beta^c_{\CHC_R} < \beta^s_{\CHC_R}$ and $\beta^c_{\ROHC_R} < \beta^s_{\ROHC_R}$, the proposed polynomial-time algorithms (Algorithms \ref{alg: CHC and ROHC recovery simple thresholding}, \ref{alg: ROHC and CHC recover via HOOI} and \ref{alg: CHC recovery aggregated SVD}) require a strictly stronger signal-to-noise ratio than the proposed unconstrained-time ones (Algorithms \ref{alg: CHC recovery combina search} and \ref{alg: ROHC recovery combinatorial search}) which leaves a significant gap between statistical optimality and computational efficiency to be discussed in Section \ref{sec: Comp lower bound for CHC and ROHC}.

\section{High-order Cluster Detection: Statistical Limits and Polynomial-time Algorithms}\label{sec: CHC and ROHC detection}

In this section, we investigate the statistical limits of both $\CHC_D$ and $\ROHC_D$. For each model, we first present the statistical lower bounds of signal strength that guarantees reliable detection, then we propose the algorithms, though being computationally intense, that provably achieve the statistical lower bounds. Finally, we introduce the computationally efficient algorithms and provide the theoretical guarantees under the stronger signal-to-noise ratio.

\subsection{$\CHC_D$ and $\ROHC_D$: Statistical Limits} \label{sec: CHC and ROHC detection information limit}

Recall \eqref{eq: CHC phase transition threshold} and \eqref{eq: ROHC phase transition threshold}, Theorems \ref{thm: CHC testing statistical lower bound} and \ref{thm: ROHC detection lower bound} below give the statistical lower bounds that guarantee reliable detection for $\CHC_D$ and $\ROHC_D$, respectively.
\begin{Theorem}[Statistical Lower Bound of $\CHC_D$] \label{thm: CHC testing statistical lower bound}
Consider $\CHC_D(\bk, \bn, \lambda)$ under the asymptotic regime \eqref{assum: asymptotic assumption} and assume
\begin{equation}\label{assum: CHC detection lower bound assumption}
 \frac{\log (n_j/k_j)}{k_i} \to 0,\quad  \frac{\log \log (n_i/k_i)  }{\log (n_j/k_j)} \to 0, \quad \text{ and }\quad k_i \log \frac{n_i}{k_i} \asymp k_j \log \frac{n_j}{k_j}
\end{equation}
for all $i, j \in [d], i \neq j$. Then if 
\begin{equation}\label{cond: CHC detection stat lower bound condition}
\frac{\lambda \prodk }{\sqrt{ \prodn }} \to 0 \quad \text{ and }\quad  \limsup_{n \to \infty} \frac{\lambda  (\prodk)^{\frac{1}{2}  }}{ \sqrt{2 (\sum_{i=1}^d k_i \log (n_i/k_i))}  }  < 1,
\end{equation}
we have $\cE_{\CHC_D}^s \to 1$. Moreover, under \eqref{eq: comp lower bound asymptotic regime}, if $\beta > \beta^s_{\CHC_D}$, $\cE^s_{\CHC_D} \to 1$.
\end{Theorem}
 
\begin{Theorem}[Statistical Lower Bound of $\ROHC_D$] \label{thm: ROHC detection lower bound} Consider $\ROHC_D(\bk, \bn, \mu)$. Under the asymptotic regime \eqref{assum: asymptotic assumption}, if $\frac{\mu}{\sqrt{ k \log(en/k)}} \to 0$, then $\cE^s_{\ROHC_D} \to 1$. Under the asymptotic regime \eqref{eq: comp lower bound asymptotic regime}, if $\beta > \beta^s_{\ROHC_D}$, $\cE^s_{\ROHC_D} \to 1$.
\end{Theorem}

Next, we present the hypothesis tests $\psi_{\CHC_D}^s$ and $\psi_{\ROHC_D}^s$ that achieve reliable detection on the statistical limits in Theorems \ref{thm: CHC testing statistical lower bound} and \ref{thm: ROHC detection lower bound}. For $\CHC_D$, define $\psi_{\CHC_D}^s := \psi_{sum} \lor \psi_{scan}$. Here, $\psi_{sum}$ and $\psi_{scan}$ are respectively the sum and scan tests:
\begin{equation}\label{eq:psi_lin}
\psi_{sum} = \mathbf{1}\left(\sum_{i_1 = 1}^{n_1} \cdots \sum_{i_d = 1}^{n_d} \bcY_{[i_1, \ldots, i_d]}/ \sqrt{n_1 \cdots n_d} > W\right) 
\end{equation}
for some to-be-specified $W>0$ and
\begin{equation}
\psi_{scan} = \mathbf{1}\left(T_{scan} > \sqrt{2\log (G_{\bk}^{\bn})}\right), \quad T_{scan} = \max_{C \in \cS_{\bk, \bn}}\frac{\sum_{(i_1, \ldots, i_d) \in C} \bcY_{[i_1, \ldots, i_d]} }{ \sqrt{k_1 \cdots k_d} },
\end{equation}
where $G_{\bk}^{\bn} = \binom{n_1}{k_1} \binom{n_2}{k_2} \cdots \binom{n_d}{k_d}$ and $\cS_{\bk, \bn}$ represents the set of all possible supports of planted signal: 
\begin{equation}\label{eq: def of S-k-n}
    \cS_{\bk, \bn} = \left\{(I_1 \times I_2\times \cdots \times I_d):   I_1 \subseteq [n_1], I_2 \subseteq [n_2], \ldots, I_d \subseteq [n_d]\text{ and } |I_i| = k_i, 1\leq i \leq d \right\}.
\end{equation}
The following Theorem \ref{thm: CHC testing stat upped bound} provides the statistical guarantee for $\psi_{\CHC_D}^s$.

\begin{Theorem}[Guarantee for $\psi_{\CHC_D}^s$]\label{thm: CHC testing stat upped bound}Consider $\CHC_D(\bk, \bn, \lambda)$. Under the asymptotic regime \eqref{assum: asymptotic assumption}, when \begin{equation}\label{cond: sum test condition}
\frac{\lambda \prodk }{\sqrt{ \prodn}} \to \infty,\quad W\to \infty, \quad W \leq c\lambda \frac{ \prodk }{\sqrt{\prodn}} (0< c < 1),
\end{equation} or when \begin{equation}\label{cond: scan test condition}
\liminf_{n \to \infty} \frac{\lambda  (\prodk)^{\frac{1}{2}  }}{ \sqrt{2 (\sum_{i=1}^d k_i \log (\frac{n_i}{k_i}  ))}  } > 1,
\end{equation}
we have $\cE_{\CHC_D}(\psi_{\CHC_D}^s) \to 0$. Under the asymptotic regime \eqref{eq: comp lower bound asymptotic regime}, $\psi^s_{\CHC_D}$ succeeds in reliable detection when $\beta < \beta^s_{\CHC_D}$.
\end{Theorem}

The test for $\ROHC_D$ is built upon the $\ROHC$ Search (Algorithm \ref{alg: ROHC recovery combinatorial search} in Section \ref{sec: CHC and ROHC recovery}) designed for $\ROHC_R$. To be specific, generate $\bcZ_1$ with i.i.d. standard Gaussian entries and calculate $\bcA = \frac{\bcY + \bcZ_1}{\sqrt{2}}$ and $\bcB = \frac{\bcY - \bcZ_1}{\sqrt{2}}$. Then $\bcA$ and $\bcB$ becomes two independent samples. Apply Algorithm \ref{alg: ROHC recovery combinatorial search} on $\bcA$ and let $(\bu_1, \ldots, \bu_d)$ be the output of Step 4 of Algorithm \ref{alg: ROHC recovery combinatorial search}. Define the test statistic as
\begin{equation*}
    \psi_{\ROHC_D}^s = \mathbf{1} \left(\bcB \times_1 \bu_1^\top/\sqrt{k_1} \times \cdots \times_d \bu_d^\top/\sqrt{k_d} \geq C \sqrt{k} \right),
\end{equation*}
where $C > 0$ is a fixed constant. We have the following theoretical guarantee for $\psi_{\ROHC_D}^s$.
\begin{Theorem}[Guarantee for $\psi_{\ROHC_D}^s$]\label{thm: ROHC detection stat upper bound}
Consider $\ROHC_D(\bk, \bn, \mu)$ under the asymptotic regime \eqref{assum: asymptotic assumption}. There exists some constant $C > 0$ such that when $\mu \geq C \sqrt{k \log n}$, $\cE_{\ROHC_D}(\psi_{\ROHC_D}^s) \to 0$. Moreover, under the asymptotic regime \eqref{eq: comp lower bound asymptotic regime}, $\psi^s_{\ROHC_D}$ succeeds in reliable detection when $\beta < \beta^s_{\ROHC_D}$.
\end{Theorem}
 Combining Theorems \ref{thm: CHC testing statistical lower bound} and \ref{thm: CHC testing stat upped bound}, we have shown $\psi_{\CHC_D}^s$ achieves sharply minimax lower bound of $\lambda$ for reliable detection of $\CHC_D$. From Theorems \ref{thm: ROHC detection lower bound} and \ref{thm: ROHC detection stat upper bound}, we see $\psi_{\ROHC}^s$ achieves the minimax optimal rate of $\mu$ for reliable detection of $\ROHC_D$. However, both $\psi_{\CHC_D}^s$ and $\psi_{\ROHC_D}^s$ are computationally inefficient.

\begin{Remark}
The proposed $\psi_{\CHC_D}^s$ and $\psi_{ROHC_D}^s$ share similar spirits with the matrix clustering algorithms in the literature \citep{butucea2013detection,brennan2018reducibility}, though the tensor structure here causes extra layer of difficulty. Particularly when $d = 2$, the lower and upper bounds results in Theorem \ref{thm: CHC testing statistical lower bound}-\ref{thm: ROHC detection stat upper bound} match the ones in \cite{butucea2013detection,brennan2018reducibility}, although the proof for high-order clustering is much more complicated.
\end{Remark}

\subsection{$\CHC_D$ and $\ROHC_D$: Polynomial-time Algorithms} \label{sec: CHC and ROHC detection poly alg limit}

Next, we introduce polynomial-time algorithms for high-order cluster detection. For $\CHC_D$, define $\psi_{\CHC_D}^c:= \psi_{sum} \lor \psi_{max}$, where $\psi_{sum}$ is defined in \eqref{eq:psi_lin} and $\psi_{max}$ is defined below based on max statistic,
\begin{equation}\label{eq:psi_max}
\psi_{max} = \mathbf{1}\left( \max_{\substack{1 \leq i_j \leq n_j\\j=1,\ldots, d}}\bcY_{[i_1, \ldots, i_d]} > \sqrt{2\sum_{i=1}^d \log n_i} \right).
\end{equation}

\begin{Theorem}[Theoretical Guarantee for $\psi_{\CHC_D}^c$] \label{thm: CHC detection upper bound for poly test} 
Consider $\CHC_D(\bk, \bn, \lambda)$. Under the asymptotic regime \eqref{assum: asymptotic assumption}, if condition \eqref{cond: sum test condition} holds or
\begin{equation}\label{cond: poly time max test condition}
\liminf_{n \to \infty} \frac{\lambda }{ \sqrt{2\sum_{i=1}^d \log n_i}} > 1,
\end{equation}
holds, then $\cE_{\CHC_D}(\psi_{\CHC_D}^c) \to 0$. Moreover, under the asymptotic regime \eqref{eq: comp lower bound asymptotic regime}, $\psi^c_{\CHC_D}$ succeeds in reliable detection when $\beta < \beta^c_{\CHC_D}$.
\end{Theorem}  

We also propose a polynomial-time algorithm for $\ROHC_D$ based on a high-order analogue of the largest matrix singular value in tensor.
Following the procedure of $\psi_{\ROHC_D}^s$, we construct $\bcA$ and $\bcB$ as two independent copies. Apply Algorithm \ref{alg: ROHC and CHC recover via HOOI} in Section \ref{sec: CHC and ROHC recovery} on $\bcA$ and let $(\bu_1, \ldots, \bu_d)$ to be the output of Step 4 of Algorithm \ref{alg: ROHC and CHC recover via HOOI}. We define 
\begin{equation}\label{eq: ROHC detection poly test}
    \psi_{\ROHC_D}^c = \psi_{sing} \lor \psi_{max}, \quad \psi_{sing} = \mathbf{1}\left(\bcB \times_1 \bu_1^\top \times \cdots \times_d \bu_d^\top \geq C \sqrt{k}\right),
\end{equation}
where $\psi_{max}$ is defined in \eqref{eq:psi_max} and $C>0$ is a fixed constant.

\begin{Theorem}[Theoretical Guarantee for $\psi_{\ROHC_D}^c$]\label{thm: ROHC detection poly test upper bound}
Consider $\ROHC_D(\bk, \bn, \mu)$ under the asymptotic regime \eqref{assum: asymptotic assumption}. There exists a constant $C > 0$ such that when 
\begin{equation} \label{cond: ROHC_D poly time test condition}
\mu \geq C n^{\frac{d}{4}} \quad \text{ or } \quad \liminf_{n \to \infty} \frac{\mu }{ \sqrt{2(\prodk) (\sum_{i=1}^d \log n_i)}  } > 1,
\end{equation} 
we have $\cE_{\ROHC_D}(\psi_{\ROHC_D}^c) \to 0$. Moreover, under the asymptotic regime \eqref{eq: comp lower bound asymptotic regime}, $\psi^c_{\ROHC_D}$ succeeds in reliable detection when $\beta < \beta^c_{\ROHC_D}$.
\end{Theorem}

Since $\beta^c_{\CHC_D} < \beta^s_{\CHC_D}$ and $\beta^c_{\ROHC_D} < \beta^s_{\ROHC_D}$, the proposed polynomial-time algorithms $\psi_{\CHC_D}^c$ and $\psi_{\ROHC_D}^c$ require a strictly stronger signal-noise-ratio than the unrestricted-time algorithms.

\section{Computational Lower Bounds}\label{sec: Comp lower bound for CHC and ROHC}

To provide the computational lower bounds for high-order clustering, it suffices to focus on the asymptotic regime \eqref{eq: comp lower bound asymptotic regime} as it also implies the computational lower bounds in the general parameterization regime \eqref{assum: asymptotic assumption}. We first consider the detection of $\CHC$. Theorem \ref{thm: CHC_D comp lower bound} below and Theorem \ref{thm: CHC detection upper bound for poly test} in Section \ref{sec: CHC and ROHC detection poly alg limit} together yield a tight computational lower bound for $\CHC_D$.
\begin{Theorem}[Computational Lower Bound of $\CHC_D$] \label{thm: CHC_D comp lower bound} Consider $\CHC_D(\bk, \bn, \lambda)$ under the asymptotic regime \eqref{eq: comp lower bound asymptotic regime}. If $\beta > \beta^c_{\CHC_D}$, then $ \liminf_{n \to \infty} \cE^c_{\CHC_D} \geq 1/2$ under the HPC detection conjecture \ref{conj: hardness of tensor clique detection}.
\end{Theorem}

Next, 
Theorems \ref{thm: CHC and ROHC recovery thresholding}, \ref{thm: CHC recovery aggregated svd}, and Theorem \ref{thm: CHC_R computation lower bound} below together give a tight computational lower bound for $\CHC_R$.

\begin{Theorem}[Computational Lower Bound of $\CHC_R$]\label{thm: CHC_R computation lower bound}
	Consider $\CHC_R(\bk, \bn, \lambda)$ under the asymptotic regime \eqref{eq: comp lower bound asymptotic regime}.
	 If $\alpha \geq 1/2$ and $\beta > (d-1)\alpha - (d-1)/2$, then $\liminf_{n \to \infty} \cE^c_{\CHC_R} \geq 1/2$ under the HPDS recovery conjecture (Conjecture \ref{conj: HPDS recovery conjecture}).
	 If $0<\alpha < 1/2, \beta > 0$, then $\liminf_{n \to \infty} \cE^c_{\CHC_R} \geq 1/2$ under the HPC detection conjecture (Conjecture \ref{conj: hardness of tensor clique detection}).
	Combined together, we have if $\beta > \beta^c_{\CHC_R}$, then $\liminf_{n \to \infty} \cE^c_{\CHC_R} \geq 1/2$ under Conjectures \ref{conj: hardness of tensor clique detection} and \ref{conj: HPDS recovery conjecture}.
\end{Theorem}

Then, we consider rank-one high-order cluster detection and recovery. By Lemma 10 in \cite{luo2020tensorsupp} Section B, we can show that the computational lower bound of $\ROHC_R$ is implied by $\ROHC_D$. We specifically have the following theorem.
\begin{Theorem}[Computational Lower Bounds of $\ROHC_D$ and $\ROHC_R$] \label{thm: ROHC_D and ROHC_R comp lower bound}
Consider $\ROHC_D(\bk, \bn, \mu)$ and $\ROHC_R(\bk, \bn, \mu)$ under the asymptotic regime \eqref{eq: comp lower bound asymptotic regime} and the HPC detection Conjecture \ref{conj: hardness of tensor clique detection}. If $\beta > \beta^c_{\ROHC_D}$, then $\liminf_{n \to \infty} \cE^c_{\ROHC_D} \geq 1/2$, $\liminf_{n \to \infty} \cE^c_{\ROHC_R} \geq 1/2.$
\end{Theorem} 
Combining Theorems \ref{thm: CHC and ROHC recovery thresholding}, \ref{thm: CHC and ROHC recovery HOOI upper bound}, \ref{thm: ROHC_D and ROHC_R comp lower bound}, and \ref{thm: ROHC detection poly test upper bound} (provided in Section \ref{sec: CHC and ROHC detection poly alg limit}), we have obtained the tight computational lower bounds for $\ROHC_D$ and $\ROHC_R$. Furthermore, since ROHC is a special case of sparse tensor PCA/SVD studied in literature \citep{zhang2019optimal-statsvd,sun2017provable}, Theorem \ref{thm: ROHC_D and ROHC_R comp lower bound} also provides a computational lower bound for the signal-to-noise ratio requirement for sparse tensor PCA/SVD.

\begin{Remark}
The computational lower bounds in Theorems \ref{thm: CHC_D comp lower bound}, \ref{thm: CHC_R computation lower bound} and \ref{thm: ROHC_D and ROHC_R comp lower bound} are for asymmetric tensor clustering under the CHC and ROHC models. To establish the computational lower bounds for a symmetric version of the CHC or ROHC models that both the planted signal and the noise tensors are symmetric, a new proof scheme is required as the same sparsity across all modes must be ensured while constructing instance tensors in performing the average-case reduction.
\end{Remark}

Theorems \ref{thm: CHC_D comp lower bound} -- \ref{thm: ROHC_D and ROHC_R comp lower bound} above are based on the HPC and HPDS conjectures. Next, we will elaborate the HPC, HPDS conjectures in Sections \ref{sec: PC in hypergraph}, \ref{sec: HPDS in hypergraph}, and discuss the evidence in Section \ref{sec: evidence for HPC}. Then in Section \ref{sec: ROHC_D and ROHC_R comp lower bound}, we provide the high level ideas on the average-case reduction from HPC and HPDS to high-order clustering, and prove these computational lower bounds. 

\subsection{Hypergraphic Planted Clique Detection} \label{sec: PC in hypergraph}

A $d$-hypergraph can be seen as an order-$d$ extension of regular graph. In a $d$-hypergraph $G = (V(G), E(G))$, each hyper-edge $e \in E$ includes an unordered group of $d$ vertices in $V$. 
Define $\mathcal{G}_d(N, 1/2)$ as Erd\H{o}s-R{\'e}nyi random $d$-hypergraph with $N$ vertices, where each hyper-edge $(i_1, \ldots, i_d)$ is independently included in $E$ with probability $\frac{1}{2}$. Given a $d$-hypergraph $G = (V(G),E(G))$, define its adjacency tensor $\bcA:= \bcA(G) \in (\{ 0,1 \}^N)^{\otimes d}$ as
\begin{equation*}
 	\bcA_{[i_1, \ldots, i_d]} = 
 	\begin{cases}
 	1, & \text{ if } (i_1, \ldots, i_d) \in E;\\
 	0, & \text{ otherwise}.	
 	\end{cases}
\end{equation*}

We define $\mathcal{G}_d(N, \frac{1}{2}, \kappa)$ as the hypergraphic planted clique (HPC) model with clique size $\kappa$. To generate $G \sim \mathcal{G}_d(N,\frac{1}{2}, \kappa)$, we sample a random hypergraph from $\mathcal{G}_d(N,\frac{1}{2})$, pick $\kappa$ vertices uniformly at random from $[N]$, denote them as $K$, and connect all hyper-edges $e$ if all vertices of $e$ are in $K$. The focus of this section is on the {\bf hypergraphic planted clique detection (HPC)} problem:
\begin{equation}\label{eq: HPC detection problem}
H^G_0: G \sim \mathcal{G}_d \left(N, 1/2\right)\quad \text{v.s.} \quad H^G_1: G \sim \mathcal{G}_d \left(N, 1/2, \kappa \right).
\end{equation}
Given the hypergraph $G$ and its adjacency tensor $\bcA$, the risk of test $\phi$ for \eqref{eq: HPC detection problem} is defined as the sum of Type-I and II errors $\cE_{\HPC_D}(\phi) = \bbP_{H_0^G} \left(\phi(\bcA) = 1\right) + \bbP_{H_1^G} \left(\phi(\bcA) = 0\right).$ Our aim is to find out the consistent test $\phi = \{\phi_N\}$ such that $\lim_{N\to \infty}\cE_{\HPC_D}(\phi_N) = 0$.

When $d=2$, HPC detection \eqref{eq: HPC detection problem} reduces to the planted clique (PC) detection studied in literature. It is helpful to have a quick review of existing results for PC before addressing HPC. Since the size of the largest clique in Erd\H{o}s R{\'e}nyi graph $G \sim \mathcal{G}_2(N,\frac{1}{2})$ converges to $2\log_2 N$ asymptotically, reliable PC detection can be achieved by exhaustive search whenever $\kappa \geq (2+\epsilon)\log_2 N$ for any $\epsilon > 0$ \citep{bollobas1976cliques}. When $\kappa = \Omega(\sqrt{N})$, many computational-efficient algorithms, including the spectral method, approximate message passing, semidefinite programming, nuclear norm minimization, and combinatorial approaches \citep{alon1998finding,ames2011nuclear, feige2000finding, ron2010finding, mcsherry2001spectral, dekel2014finding, deshpande2015finding, chen2016statistical}, have been developed for PC detection. Despite enormous previous efforts, no polynomial-time algorithm has been found for reliable detection of PC when $\kappa = o(N^{1/2})$ and it has been widely conjectured that no polynomial-time algorithm can achieve so. The hardness conjecture of PC detection was strengthened by several pieces of evidence, including the failure of Metropolis process methods \citep{jerrum1992large}, low-degree polynomial methods \citep{hopkins2018statistical, brennan2019average}, statistical query model \citep{feldman2017statistical}, Sum-of-Squares \citep{barak2019nearly, deshpande2015improved, meka2015sum}, landscape of optimization \citep{gamarnik2019landscape}, etc.

When moving to HPC detection \eqref{eq: HPC detection problem} with $d\geq 3$, the computational hardness remains little studied. \cite{bollobas1976cliques} proved that $\frac{K^d_N}{\left(d!\log_2(N)\right)^{1/(d-1)}} \overset{a.s.} \to 1$ if $K_N^d$ is the largest clique in $G \sim \mathcal{G}_d(N, \frac{1}{2})$. So HPC detection problem \eqref{eq: HPC detection problem} is statistical possible by exhaustive search when $\kappa \geq \left((d!+\epsilon) \log_2(N)\right)^{1/(d-1)}$ for any $\epsilon > 0$. However, \cite{zhang2018tensor} observed that the spectral algorithm solves HPC detection if $\kappa = \Omega(\sqrt{N})$ but fails when $\kappa = N^{\frac{1}{2} - \epsilon}$ for any $\epsilon > 0$.
We present the following hardness conjecture for HPC detection.
\begin{Conjecture}[HPC Detection Conjecture]\label{conj: hardness of tensor clique detection}
Consider the HPC detection problem \eqref{eq: HPC detection problem} and suppose $d\geq 2$ is a fixed integer. If
\begin{equation}\label{eq:HPC-condition}
\limsup_{N \to \infty} \log \kappa / \log \sqrt{N} \leq 1 - \tau \quad \text{for any }\tau > 0,
\end{equation} 
for any polynomial-time test sequence $\{\phi \}_N: \bcA \to \{0,1\}$, $\liminf_{N \to \infty} \cE_{\HPC_D}(\phi(\bcA)) \geq \frac{1}{2}$. 
\end{Conjecture}

\begin{Remark}[Choice of Type-I, II Error Lower Bound]
We set the lower bound for the sum of Type-I, II errors to be $1/2$ in the HPC Detection Conjecture above (i.e., $\{\phi \}_N: \bcA \to \{0,1\}$, $\liminf_{N \to \infty} \cE_{\HPC_D}(\phi(\bcA)) \geq 1/2$).   
In the literature, there is no universal choice of this constant. For example, \cite{berthet2013complexity} considers PC detection conjecture with the sum of type I and type II errors to be some constant close to $1$; \cite{ma2015computational} uses the PC detection conjecture with the error constant $2/3$; \cite{brennan2018reducibility,brennan2019optimal,brennan2020reducibility,hajek2015computational} choose this constant to be $1$. 
\end{Remark}

In Section \ref{sec: evidence for HPC}, we provide two pieces of evidence for HPC detection conjecture: a general class of Monte Carlo Markov Chain process methods \citep{jerrum1992large} and a general class of low-degree polynomial tests \citep{hopkins2017bayesian, hopkins2018statistical, kunisky2019notes, brennan2019average} fail to solve HPC detection under the asymptotic condition \eqref{eq:HPC-condition}. Also, see a recent note \cite{luo2020open} for several open questions on HPC detection conjecture, in particular, whether HPC detection is equivalently hard as PC detection.

\subsection{Hypergraphic Planted Dense Subgraph}\label{sec: HPDS in hypergraph}
We consider the {\bf hypergraphic planted dense subgraph (HPDS)}, a hypergraph model with denser connections within a community and sparser connections outside, in this section. Let $\cG_d$ be a $d$-hypergraph. To generate a HPDS $G = (V(G), E(G)) \sim \mathcal{G}_d(N,\kappa, q_1, q_2 )$ with $q_1 > q_2$, we first select a size-$\kappa$ subset $K$ from $[N]$ uniformly at random, then for each hyper-edge $e = (i_1, \ldots, i_d)$,
\begin{equation*}
\bbP\left(e\in E(G)\right) = \left\{\begin{array}{ll}
q_1, & i_1,\ldots, i_d \in K\\
q_2, & \text{otherwise}.
\end{array}\right.
\end{equation*}
The aim of HPDS detection is to test
\begin{equation}\label{eq: HPDS detection}
    H_0: G \sim \cG_d(N,q_2) \quad \text{versus} \quad H_1: G \sim \cG_d(N, \kappa, q_1, q_2);
\end{equation}
the aim of HPDS recovery is to locate the planted support $K$ given $G \sim \mathcal{G}_d(N, \kappa, q_1, q_2)$. 

When $d=2$, HPDS reduces to the planted dense subgraph (PDS) considered in literature. Various statistical limits of PDS have been studied \citep{chen2016statistical,hajek2015computational,arias2014community,verzelen2015community,brennan2018reducibility,feldman2017statistical} and generalizations of PDS recovery has also been considered in \cite{hajek2016information,montanari2015finding,candogan2018finding}. In \cite{hajek2015computational,brennan2018reducibility}, a reduction from PC has shown the statistical and computational phase transition for PDS detection problem for all $q_1 > q_2$ with $q_1 -q_2 = O(q_2)$ where $q_2 = \tilde{\Theta}(N^{-\beta})$. For PDS recovery problem, \cite{chen2016statistical,brennan2018reducibility,hajek2015computational} observed that PDS appears to have a detection-recovery gap in the regime when $\kappa \gg \sqrt{N}$. 

When moving to HPDS detection, if $q_1 = \omega(q_2)$, the computational barrier for this problem is conjectured to be the log-density threshold $\kappa = \tilde{\Theta}(N^{\log_{q_2} q_1})$ when $\kappa \ll \sqrt{N}$ \citep{chlamtac2012everywhere,chlamtavc2017minimizing}. Recently, \cite{chlamtac2018sherali} showed that $\tilde{\Omega}(\log N)$ rounds of the Sherali-Adams hierarchy cannot solve the HPDS detection problem below the log-density threshold in the regime $q_1 = \omega(q_2)$. The HPDS recovery, to the best of our knowledge, remains unstudied in the literature. 

In the following Proposition \ref{prop: spectral method for HPDS}, we show that a variant of Aggregated-SVD (presented in Algorithm \ref{alg: HPDS recovery aggregated SVD}) requires a restricted condition on $\kappa, q_1, q_2, N$ for reliable recovery in $\HPDS$ in the regime $\kappa \gg \sqrt{N}$.

\begin{algorithm}[h] \caption{Support Recovery of HPDS via Aggregated-SVD} \label{alg: HPDS recovery aggregated SVD}
	\begin{algorithmic}[1]
	\State \textbf{Input:} $\bcA$.
	    \State Let $\widecheck{\bcA} = \bcA_{\left[1:  \lfloor\frac{N}{d} \rfloor, \lfloor\frac{N}{d} \rfloor + 1: 2\lfloor\frac{N}{d} \rfloor,\ldots, (d-1) \lfloor \frac{N}{d} \rfloor+1: N \rfloor  \right]}$.
		\State Let $\widetilde{\bcA}_{[i_1, \ldots, i_d]} = \frac{\widecheck{\bcA}_{[i_1, \ldots, i_d]} - q_2}{\sqrt{q_2(1-q_2) }}$ for all $1\leq i_1 \leq \lfloor \frac{N}{d} \rfloor, \ldots, (d-1) \lfloor \frac{N}{d} \rfloor+1 \leq i_d \leq N $. Then apply Algorithm \ref{alg: CHC recovery aggregated SVD} with input $\widetilde{\bcA}$ and denote the estimated support for each mode of $\widetilde{\bcA}$ as $\hat{K}_i$.
		\State Compute $\hat{K} = \bigcup_{i=1}^d \hat{K}_i$.
		\State    \textbf{Output:} $ \hat{K}$.
	\end{algorithmic}

\end{algorithm}
 
\begin{Proposition}\label{prop: spectral method for HPDS}
Suppose $G \sim \cG_d(N, \kappa, q_1, q_2)$ with $q_1 > q_2$. Let $\bcA$ be the adjacency tensor of $G$. When $\liminf_{N \to \infty} \log_N \kappa \geq 1/2$, and
\begin{equation} \label{cond: spectral HPDS lambda requirement}
    \limsup_{N \to \infty} \log_N \left(\frac{\kappa^{d-1} (q_1-q_2)}{\sqrt{q_2(1-q_2)}}\right) \geq \frac{d}{2} -\frac{1}{2},
\end{equation} 
Algorithm \ref{alg: HPDS recovery aggregated SVD} recovers the support of the planted dense subgraph with probability at least $1 - d \left(N/d\right)^{-c} - C \exp\left(-c N/d \right)$ for some $c, C > 0$.
\end{Proposition}

On the other hand, the theoretical analysis in Proposition \ref{prop: spectral method for HPDS} breaks down when condition \eqref{cond: spectral HPDS lambda requirement} does not hold. We conjecture that the signal-to-noise ratio requirement in \eqref{cond: spectral HPDS lambda requirement} is essential for HPDS recovery and propose the following computational hardness conjecture. 
\begin{Conjecture}[HPDS Recovery Conjecture] \label{conj: HPDS recovery conjecture}
	Suppose $G \sim \cG_d(N, \kappa, q_1, q_2)$ with $1-\Omega(1) > q_1 > q_2$. Denote its adjacency tensor as $\bcA$. If
	\begin{equation} \label{eq: HPDS conjecture}
		\liminf_{N \to \infty} \log_N \kappa \geq \frac{1}{2} \quad \text{and}\quad  \limsup_{N \to \infty} \log_N \left(\frac{\kappa^{d-1} (q_1-q_2)}{\sqrt{
		q_2(1-q_2)}}\right) < \frac{d}{2} -\frac{1}{2},
	\end{equation}
	then for any randomized polynomial-time algorithm $\{\phi \}_N$, $\liminf_{N \to \infty} \cE_{\HPDS_R} (\phi (\bcA)) \geq \frac{1}{2}$.
\end{Conjecture}

In the proof of Proposition \ref{prop: spectral method for HPDS}, we provide evidence for Conjecture \ref{conj: HPDS recovery conjecture} by showing a variant of Aggregated-SVD fails to solve HPDS recovery under the PC detection conjecture. A stronger piece of evidence for Conjecture \ref{conj: HPDS recovery conjecture} via low-degree polynomial method is given in Section \ref{sec: low-degree-evidence-HPDS}. 

\subsection{Evidence for HPC Detection Conjecture}\label{sec: evidence for HPC}

In this section, we provide two pieces of evidence for HPC conjecture \ref{conj: hardness of tensor clique detection} via Monte Carlo Markov Chain process and low-degree polynomial test and one piece of evidence for HPDS recovery conjecture \ref{conj: HPDS recovery conjecture} via low-degree polynomial method. 

\subsubsection{Evidence of HPC Conjecture \ref{conj: hardness of tensor clique detection} via Metropolis process}

We first show a general class of Metropolis processes are not able to detect or recover the large planted clique in hypergraph.
Motivated by \cite{alon2007testing}, in Lemma \ref{lm: HPC recovery hard imply detection hard} we first prove that if it is computationally hard to recover a planted clique in HPC, it is also computationally hard to detect.

\begin{Lemma}\label{lm: HPC recovery hard imply detection hard}
Assume $\kappa > \Omega(\log N)$. Consider the $\HPC_D(N, 1/2, \kappa)$ problem: test
\begin{equation*}
H_0: G\sim \cG_d(N, 1/2)\quad \text{versus} \quad H_1: G\sim \cG_d(N, 1/2, \kappa)
\end{equation*}
and $\HPC_R(N, 1/2, \kappa)$ problem:  recover the exact support of the planted clique if $H_1$ holds. If there is no polynomial time recovery algorithm can output the right clique of $\HPC_R(N,1/2,\kappa)$ with success probability at least $1-1/N$, then there is no polynomial time detection algorithm can output the right hypothesis for $\HPC_D(N,1/2,\kappa/3)$ with probability $1-1/(4N^d)$.
\end{Lemma}
By Lemma \ref{lm: HPC recovery hard imply detection hard}, to show the computational hardness of HPC detection, we only need to show the HPC recovery. 

Motivated by the seminal work of \cite{jerrum1992large}, we consider the following simulated annealing method for planted clique recovery in hypergraph. Given a hypergraph $G = (V, E) \sim \mathcal{G}_d(N, 1/2,\kappa)$ on the vertex set $V =\{0, \ldots, N-1\}$ and a real number $\theta\geq 1$, we consider a Metropolis process on the state space of the collection $\Gamma \subseteq 2^{V}$ of all cliques in $G$, i.e., all subsets of $V$ which induces the complete subgraph in $G$. A transition from state $K$ to state $K'$ is allowed if $|K \oplus K'| \leq 1$ (Here, $K\oplus K' = \{i: i\in K, i\notin K'\}\bigcup \{i: i\in K', i\notin K\}$ is the set symmetric difference). 

For all distinct states $K,K' \in \Gamma$, the transition probability from $K$ to $K'$ is
\begin{equation}\label{eq: transition prob of MCMC}
P(K, K') = \left\{ \begin{array}{l l}
\frac{1}{N \theta}, & \text{if } K \oplus K' = 1, K \supset K' ;\\
\frac{1}{N}, & \text{if } K \oplus K' = 1, K \subset K';\\
0, & \text{if } |K\oplus K'| \geq 2.
\end{array} \right.
\end{equation} 
The loop probability $P(K, K) = 1 - \sum_{K' \neq K} P(K , K')$ are defined by complementation. The transition probability can be interpreted by the following random process. Suppose the current state is $K$. Pick a vertex $v$ uniformly at random from $V$. 
\begin{enumerate}[leftmargin=*]
\item If $v\notin K$ and $K\bigcup\{v\}$ is a clique, then let $K'=K \bigcup \{v\}$; 
\item If $v\notin K$ and $K\bigcup\{v\}$ is not a clique, then let $K' = K$;
\item If $v \in K$, with probability $\frac{1}{\theta}$, set $K'=K\setminus\{v\}$, else set $K' = K$. 
\end{enumerate}

When $\theta > 1$, the Metropolis process defined above is aperiodic and then has a unique statitionary distribution. Let $\pi: \Gamma \to [0, 1]$ be defined as $$\pi(K) = \frac{\theta^{|K|}}{\sum_{K \in \Gamma} \theta^{|K|}}.$$
We can check that $\pi$ satisfies the following detailed balance property: 
\begin{equation}
\theta^{|K|} P(K, K') = \theta^{|K'|} P(K', K), \text{ for all } K, K' \in \Gamma.
\end{equation}
This means $\pi$ is indeed the stationary distribution of this Markov chain. The following theorem shows that it takes superpolynomial time to locate a clique in $G$ of size $\Omega \left((\log_2 N)^{1/(d-1)}\right)$ by described Metropolis process.
\begin{Theorem}[Hardness of Finding Large Clique in $\mathcal{G}_d(N, 1/2, N^\beta), 0< \beta < \frac{1}{2}$] \label{thm: Metropolis process is hard to find large clique}
Suppose $\epsilon > 0$ and $0 < \beta < \frac{1}{2}$. For almost every $G \in \mathcal{G}_d (N, 1/2, N^\beta)$ and every $\theta > 1$, there exists an initial state from which the expected time for the Metropolis process to reach a clique of size at least $m$ exceeds $N^{\Omega \left((\log_2 N)^{1/(d-1)} \right)}$. Here, 
$$m = 2\left\lceil \left(\left(1+ \frac{2}{3} \epsilon\right) \frac{d!}{2} \log_2 N\right)^{\frac{1}{d-1}}\right\rceil - \left\lceil \left(\left(1+ \frac{2}{3} \epsilon\right) (d-1)! \log_2 N\right)^{\frac{1}{d-1}}\right\rceil \asymp_d (\log_2 N)^{\frac{1}{d-1}}.$$
\end{Theorem}

\subsubsection{Evidence of HPC Conjecture \ref{conj: hardness of tensor clique detection} via Low-degree Polynomial Test} \label{sec: low-degree-evid-HPC-conjecture}

We also consider the low-degree polynomial tests to establish the computational hardness for hypergraphic planted clique detection. The idea of using low-degree polynomial to predict the statistical and computational gap is recently developed in a line of papers \citep{hopkins2017bayesian, hopkins2017power,hopkins2018statistical,barak2019nearly}. Many state-of-the-art algorithms, such as spectral algorithm, approximate message passing  \citep{donoho2009message} can be represented as low-degree polynomial functions as the input, where ``low" means logarithmic in the dimension. In comparison to sum-of-squares (SOS) computational lower bounds, the low-degree method is simpler to carry out and appears to always yields the same results for natural average-case problems, such as the planted clique detection \citep{hopkins2018statistical, barak2019nearly}, community detection in stochastic block model \citep{hopkins2017bayesian, hopkins2018statistical}, the spiked tensor model \citep{hopkins2017power, hopkins2018statistical, kunisky2019notes}, the spiked Wishart model \citep{bandeira2020computational}, sparse PCA \citep{ding2019subexponential}, spiked Wigner model \citep{kunisky2019notes}, sparse clustering \citep{loffler2020computationally}, certifying RIP \citep{ding2020average} and a variant of planted clique and planted dense subgraph models \citep{brennan2019average}. It is gradually believed that the low-degree polynomial method is able to capture the essence of what makes SOS succeed or fail \citep{hopkins2017bayesian, hopkins2017power, hopkins2018statistical, kunisky2019notes, raghavendra2018high}. Therefore, we apply this method to give the evidence for the computational hardness of HPC detection \eqref{eq: HPC detection problem}. Specifically, we have the following Theorem \ref{prop: fail of low-degree polynomial test} for low degree polynomial tests in HPC.

\begin{Theorem}[Failure of Low-degree Polynomial Tests for HPC] \label{prop: fail of low-degree polynomial test}
	Consider the HPC detection problem \eqref{eq: HPC detection problem} for $\kappa = N^\beta \,  (0<\beta < \frac{1}{2})$. Suppose $\bcA$ is the adjacency tensor of $G$ and $f(\bcA)$ is a polynomial test such that $\bbE_{H_0^G} f(\bcA) = 0$, $\bbE_{H_0^G} (f^2(\bcA)) = 1$, and the degree of $f$ is at most $D$ with $D \leq C \log N$ for constant $C>0$. Then we have $\bbE_{H_1^G} f(\bcA) = O(1)$.
\end{Theorem}

It has been widely conjectured in the literature that for a broad class of hypothesis testing problems: $H_0$ versus $H_1$, there is a test with runtime $n^{\tilde{O}(D)}$ and Type I + II error tending to zero if and only if there is a successful $D$-simple statistic, i.e., a polynomial $f$ of degree at most $D$, such that $\bbE_{H_0} f(X) = 0$, $\bbE_{H_0} (f^2(X)) = 1$, and $\bbE_{H_1} f(X) \to \infty$ \citep{hopkins2018statistical,kunisky2019notes,brennan2019average,ding2019subexponential}. Thus, Theorem \ref{prop: fail of low-degree polynomial test} provides the firm evidence that there is no polynomial-time test algorithm that can reliably distinguish between $\mathcal{G}_d(N, 1/2)$ and $\mathcal{G}_d(N, 1/2, N^\beta)$ for $0<\beta < 1/2$.

\subsubsection{Evidence of HPDS Recovery Conjecture \ref{conj: HPDS recovery conjecture} via Low-degree Polynomial Method} \label{sec: low-degree-evidence-HPDS}
Compared to the hardness evidence for the hypothesis testing problems, it is much less explored in the literature to establish hardness evidence for the {\it estimation or recovery problems}. Recently, \cite{schramm2020computational} provides the first sharp computational lower bounds for recovery in biclustering and planted dense subgraph via the low-degree polynomial method and resolve the ``detection-recovery gap'' open problem mentioned in \cite{ma2015computational,chen2016statistical,brennan2018reducibility,hajek2015computational}. In this work, we leverage the results in \cite{schramm2020computational} and provide the firm evidence for HPDS recovery conjecture \ref{conj: HPDS recovery conjecture} via the low-degree polynomial method.

Recall the HPDS recovery problem in Section \ref{sec: HPDS in hypergraph}. Let $G \sim \cG_d(N,\kappa,q_1,q_2)$ with $q_1 > q_2 $ and planted subset $K$. Denote $v_1\in \{0,1\}$ as the membership of vertex $1$ such that $v_1=1$ if the first vertex is in $K$ and $v_1=0$ otherwise. The following theorem shows that it is impossible to estimate $v_1$ well in the conjectured hard regime via low-degree polynomials, which implies the computational difficulty of recovering $K$ in general. 
\begin{Theorem}[Failure of Low-degree Polynomials for HPDS Recovery] \label{thm: low-degree-HPDS-evidence}
	Suppose $G \sim \cG_d(N,\kappa,q_1,q_2)$ with $q_1 > q_2$ and $\bcA$ is the adjacency tensor of $G$. For any $0<r <1$ and  $D \geq 1 $, if 
	\begin{equation} \label{ineq: HPDS-low-degree-SNR-lower-bound}
		\frac{q_1 - q_2}{\sqrt{q_2(1-q_1)}} \leq \frac{\sqrt{r}}{D+1} \min \left((D(d-1) + 1)^{-d/2}, \frac{N^{(d-1)/2} }{\sqrt{Dd(d-1)} \kappa^{d-1} } \right),
	\end{equation}
	then for any $f: \bcA \to \bbR$ with degree at most $D$, we have $\bbE(f(\bcA) - v_1)^2 \geq \frac{\kappa}{N} - (\frac{\kappa}{N})^2 (1+ \frac{r}{(1-r)^2})$.
	
	In particular, suppose $q_2 < q_1 < 1- \Omega(1)$. Consider the asymptotic regime of Conjecture \ref{conj: HPDS recovery conjecture} that
	$$ \liminf_{N \to \infty} \log_N \kappa \geq \frac{1}{2} \quad \text{and}\quad  \limsup_{N \to \infty} \log_N \left(\frac{\kappa^{d-1} (q_1-q_2)}{\sqrt{q_2(1-q_2)}}\right) < \frac{d}{2} -\frac{1}{2}.$$
Let $f_0$ be the trivial constant estimator of $v_1$: $f_0(\bcA) = \kappa/N$. Then for any polynomial $f:\bcA \to \mathbb{R}$ with degree at most $D$ with $D \leq \text{polylog}(N)$, we have 
	$$\liminf_{N\to \infty}\frac{\mathbb{E}(f(\bcA) - v_1)^2}{\mathbb{E}(f_0(\bcA) - v_1)^2} \geq 1.$$
\end{Theorem}
Theorem \ref{thm: low-degree-HPDS-evidence} shows that under the conjectured hard regime of HPDS \eqref{eq: HPDS conjecture} and $q_2 < q_1 < 1- \Omega(1)$, the mean square error of any $f$ with degree equal or less than polylog$(N)$ is no better than the trivial estimator $f_0$. This gives strong evidence for the HPDS recovery conjecture \ref{conj: HPDS recovery conjecture}.

\subsection{Proofs of Computational Lower Bounds} \label{sec: ROHC_D and ROHC_R comp lower bound}

Now, we are in position to prove the computational lower bounds. Before the detailed analysis, we first outline the high-level idea. 

Consider a hypothesis testing problem $B$: $H_0$ versus $H_1$. To establish a computational lower bound for $B$, we can construct a randomized polynomial-time reduction $\varphi$ from the conjecturally hard problem $A$ to $B$ such that the total variation distance between $\varphi(A)$ and $B$ converges to zero under both $H_0$ and $H_1$. If such a $\varphi$ can be found, whenever there exists a polynomial-time algorithm $\phi$ for solving $B$, we can also solve $A$ using $\phi \circ \varphi$ in polynomial-time. Since $A$ is conjecturally hard, we can conclude that $B$ must also be polynomial-time hard by the contradiction argument. To establish the computational lower bound for a recovery problem, we can either follow the same idea above or establish a reduction from recovery to an established detection lower bound. A key challenge of average-case reduction is often how to construct an appropriate randomized polynomial-time map $\varphi$.

We summarize the procedure of constructing randomized polynomial-time maps for the high-order clustering computational lower bounds as follows.
\begin{itemize}
    \item Input: Hypergraph $G$ and its adjacency tensor $\bcA$
    \item {\bf Step 1}: Apply the {\it rejection kernel} technique, which was proposed by \cite{ma2015computational} and formalized by \cite{brennan2018reducibility}, to simultaneously map $\Bern(p)$  distribution to $N(\xi,1)$ and $\Bern(q)$ distribution to $N(0,1)$ approximately. 
    \item {\bf Step 2}: Simultaneously change the magnitude and sparsity of the planted signal guided by the target problem. In this step, we develop several new techniques and apply several ones in the literature. In $\CHC_D$ (Algorithm \ref{alg: CHC detection reduction}), we use the average-trick idea in \cite{ma2015computational}; in $\CHC_R$ (Algorithm \ref{alg: CHC recovery reduction}), we use the invariant property of Gaussian to handle the multiway-symmetricity of hypergraph; to achieve a sharper scaling of signal strength and sparsity in $\ROHC_D, \ROHC_R$ (Algorithm \ref{alg: ROHC reduction}), the {\it tensor reflection cloning}, a generalization of reflection cloning \citep{brennan2018reducibility}, is introduced that spreads the signal in the planted high-order cluster along each mode evenly, maintains the independence of entries in the tensor, and only mildly reduces the signal magnitude.
    \item {\bf Step 3}: Randomly permute indices of different modes to transform the symmetric planted signal tensor to an asymmetric one (Lemmas 14 and 16 in \cite{luo2020tensorsupp}) that maps to the high-order clustering problem. 
\end{itemize}

Then, we give a detailed proof of Theorem \ref{thm: ROHC_D and ROHC_R comp lower bound}, i.e., computational lower bounds for $\ROHC_D$ and $\ROHC_R$. The proofs for the computational limits of $\CHC_D$ and $\CHC_R$ are similar and postponed to the supplementary materials \citep{luo2020tensorsupp}.

We first introduce the rejection kernel scheme given in Algorithm 9 in \cite{luo2020tensorsupp} Section C, which simultaneously maps $\Bern(p)$ to distribution $f_X$ and $\Bern(q)$ to distribution $g_X$ approximately. In our high-order clustering problem, $f_X$ and $g_X$ are $N(\xi, 1)$ and $N(0, 1)$, i.e., the distribution of the entries inside and outside the planted cluster, respectively. Here, $\xi$ is to be specified later. Denote $\RK(p \to f_X, q \to g_X, T)$ as the rejection kernel map, where $T$ is the number of iterations in the rejection kernel algorithm.

We then propose a new tensor reflection cloning technique in Algorithm \ref{alg: Tensor Reflecting Clone}. Note that the input tensor $\bcW_0$ to Algorithm \ref{alg: Tensor Reflecting Clone} often has independent entries and a sparse planted cluster, we multiply $\bcW^{\sigma^{\otimes d}}$, a random permutation of $\bcW_0$, by $\frac{\A + \B}{\sqrt{2}}$ in each mode to ``spread" the signal of the planted cluster along all modes while keep the entries independent. We prove some properties related to tensor reflection cloning in Lemma 16 of \cite{luo2020tensorsupp} Section C.
\begin{algorithm}[h]\caption{Tensor Reflecting Cloning}\label{alg: Tensor Reflecting Clone} 
\begin{algorithmic}[1]
\State {\bf Input:} Tensor $\bcW_0 \in \bbR^{n^{\otimes d}}$ ($n$ is an even number), number of iterations $\ell$.
	\State Initialize $\bcW = \bcW_0$.
	\State For $i = 1, \ldots, \ell$, do:
	\begin{enumerate}[label=(\alph*)]
	    \item Generate a permutation $\sigma$ of $[n]$ uniformly at random.
	    \item Calculate
	    \begin{equation*}
	        \bcW' = \bcW^{\sigma^{\otimes d}} \times_1 \frac{\A + \B}{\sqrt{2}} \times \cdots \times_d \frac{\A + \B}{\sqrt{2}},
	    \end{equation*} where $\bcW^{\sigma^{\otimes d}}$ means permuting each mode indices of $\bcW$ by $\sigma$ and $\B$ is a $n \times n$ matrix with ones on its anti-diagonal and zeros elsewhere and $\A$ is given by 
\begin{equation}\label{eq: A matrix}
    \A = \left[ \begin{array}{cc}
         \I_{\frac{n}{2}}& 0  \\
         0 & - \I_{\frac{n}{2}}
    \end{array} \right],
\end{equation} 
where $\I_{n/2}$ is a $n/2 \times n/2$ identity matrix.
	    \item Set $\bcW = \bcW'$.
	\end{enumerate}
	\State {\bf Output:} $\bcW$.
\end{algorithmic}
\end{algorithm}

We construct the randomized polynomial-time reduction from $\HPC$ to $\ROHC$ in Algorithm \ref{alg: ROHC reduction}.
\begin{algorithm}[h]\caption{Randomized Polynomial-time Reduction from HPC to $\ROHC$}\label{alg: ROHC reduction}
\begin{algorithmic}[1]
\State {\bf Input:} Hypergraph $G \sim \cG_d(n)$, number of iterations $\ell$.
	\State Let $\RK_G = \RK(1 \to N(\xi,1), \frac{1}{2} \to N(0,1), T)$ where $T = \lceil 2(d+1)\log_2 n\rceil$ and $\xi = \frac{\log 2}{ 2 \sqrt{ 2(d+1)\log n + 2 \log 2 }}$ and  compute the symmetric tensor $\bcW \in \bbR^{n^{\otimes d}}$ with $\bcW_{[i_1, \ldots, i_d]} = \RK_G(\mathbf{1}((i_1, \ldots, i_d) \in E(G) ))$. Let the diagonal entries of $\bcW_{[i, \ldots, i]}$ to be i.i.d. $N(0,1)$. 
		\State Generate $(d!-1)$ i.i.d. symmetric random tensor $\bcB^{(1)}, \ldots, \bcB^{(d!-1)}$ in the following way: their diagonal values are $0$ and non-diagonal values are i.i.d. $N(0,1)$. Given any non-diagonal index $\bi = (i_1, \ldots, i_d)$ ($i_1 \leq i_2 \leq \ldots \leq i_d$), suppose it has $D\, (D \leq d!)$ unique permutations and denote them as $\bi_{(0)}: = \bi, \bi_{(1)}, \ldots, \bi_{(D-1)}$, then we transform $\bcW$ in the following way
	\begin{equation*}
	    \left( \begin{array}{c}
	         \bcW_{\bi_{(0)}}  \\
	         \bcW_{\bi_{(1)}}\\
	         \vdots \\
	         \bcW_{\bi_{(D-1)}}
	    \end{array}  \right) = \left( \frac{\mathbf{1}}{\sqrt{d!}}, \left(\frac{\mathbf{1}}{\sqrt{d!}} \right)_{\perp}  \right)_{[1:D,:]} \times \left(  
	    \begin{array}{c}
	         \bcW_{[i_1, \ldots, i_d]}  \\
	         \bcB^{(1)}_{[i_1, \ldots, i_d]} \\
	         \vdots\\
	         \bcB^{(d!-1)}_{[i_1, \ldots, i_d]} 
	    \end{array} \right).
	\end{equation*} 
	Here $\frac{\mathbf{1}}{\sqrt{d!}}$ is a $\bbR^{d!}$ vector with all entries to be $\frac{1}{\sqrt{d!}}$ and $\left(\frac{\mathbf{1}}{\sqrt{d!}} \right)_{\perp} \in \bbR^{d! \times (d! - 1)}$ is an orthogonal complement of $\frac{\mathbf{1}}{\sqrt{d!}}$. 
	\State Generate independent permutations $\sigma_1, \ldots, \sigma_{d-1}$ of $[n]$ uniformly at random and let $\bcW = \bcW^{\id, \sigma_1, \ldots, \sigma_{d-1}}$.
	\State Apply Tensor Reflecting Cloning to $\bcW$ with $\ell$ iterations.
	\State {\bf Output:} $\bcW$.
\end{algorithmic}	
\end{algorithm}
The next lemma shows that the randomized polynomial-time mapping we construct in Algorithm \ref{alg: ROHC reduction} maps $\HPC$ to $\ROHC$ asymptotically. 

\begin{Lemma}\label{lm: ROHS reduction guarantee}
Suppose that $n$ is even and sufficiently large and $2^{\ell} \kappa < n/C$ for some large $C > 1$. Let $\xi = \frac{\log 2}{ 2 \sqrt{ 2(d+1)\log n + 2 \log 2 }}$. Then the randomized polynomial-time map $\varphi: \cG_d(n) \to \bbR^{n^{\otimes d}}$ in Algorithm \ref{alg: ROHC reduction} satisfies if $G \sim \cG_d(n, \frac{1}{2})$, it holds that
\begin{equation*}
    \TV \left( \cL(\varphi(G)), N(0,1)^{\otimes (n^{\otimes d})} \right) = O(1/n),
\end{equation*} and if $G \sim \cG_d(n, \frac{1}{2}, \kappa)$, there is a prior $\pi$ on unit vectors in $\cV_{n, 2^{\ell}\kappa}$ such that 
\begin{equation*}
    \TV\left(  \cL(\varphi(G)), \int \cL \left( \frac{\xi \kappa^{\frac{d}{2}}}{\sqrt{d!}} \bu_1 \circ \cdots \circ \bu_d + N(0,1)^{\otimes (n^{\otimes d})} \right) d \pi(\bu_1, \ldots, \bu_d) \right) = O(1/\sqrt{\log n}).
\end{equation*} 
Here $\TV$ denotes the total variation distance and $\cL(X)$ denotes the distribution of random variable $X$.
\end{Lemma} 
Lemma \ref{lm: ROHS reduction guarantee} specifically implies that if $k = 2^\ell \kappa, \mu = \frac{\xi \kappa^{\frac{d}{2}}}{\sqrt{d!}}$ with $\xi = \frac{\log 2}{ 2 \sqrt{ 2(d+1)\log n + 2 \log 2 }}$, the reduction map $\varphi(G)$ we constructed from Algorithm \ref{alg: ROHC reduction} satisfies $\TV(\varphi(\HPC_D(n,\frac{1}{2},\kappa)), \ROHC_D(\bn,\bk,\mu))\to 0$ under both $H_0$ and $H_1$.

Next, we prove the computational lower bound of $\ROHC_D$ under the asymptotic regime \eqref{eq: comp lower bound asymptotic regime} by a contradiction argument. 
\begin{itemize}[leftmargin=*]
\item If $\alpha \geq \frac{1}{2}$ ($\alpha$ is defined in \eqref{eq: comp lower bound asymptotic regime}), i.e., in the dense cluster case, let $\ell = \lceil \frac{2}{d} \beta \log_2 n\rceil$ and $\varphi$ be this mapping from Algorithm \ref{alg: ROHC reduction}. Suppose $\kappa = \lceil n^\gamma \rceil$ in $\HPC_D(n,\frac{1}{2},\kappa)$, then after mapping, the sparsity and signal strength in \eqref{eq: comp lower bound asymptotic regime} of $\ROHC(\bn,\bk,\mu)$ model satisfies 
\begin{equation*}
\begin{split}
    &\lim_{n \to \infty} \frac{\log (\mu/k^\frac{d}{2})^{-1} }{\log n} = \frac{\frac{d}{2} (\frac{2}{d}\beta + \gamma)\log n - \frac{d}{2}\gamma \log n }{\log n} = \beta, \quad\lim_{n \to \infty} \frac{\log k}{\log n} = \frac{2}{d}\beta + \gamma =:\alpha.
\end{split}
\end{equation*}

If $\beta > (\alpha -\frac{1}{2}) \frac{d}{2}$, there exists a sequence of polynomial-time tests $\{\phi_n\}$ such that $\liminf_{n \to \infty} \cE_{\ROHC_D}(\phi_n) < \frac{1}{2}$. Then by Lemmas \ref{lm: ROHS reduction guarantee} and 11 in \cite{luo2020tensorsupp} Section C, we have $\liminf_{n \to \infty} \cE_{\HPC_D}(\phi_n\circ \varphi ) < \frac{1}{2}$, i.e. $\phi_n \circ \varphi$ has asymptotic risk less than to $\frac{1}{2}$ in HPC detection. On the other hand, the size of the planted clique in HPC satisfies $\lim_{n \to \infty}\frac{\log \kappa}{\log n} = \gamma = \alpha - \frac{2}{d}\beta < \alpha - (\alpha - \frac{1}{2}) = \frac{1}{2}$. The combination of these two facts contradicts HPC detection conjecture \ref{conj: hardness of tensor clique detection}, so we conclude there are no polynomial-time tests $\{\phi_n \}$ that make $\liminf_{n \to \infty} \cE_{\ROHC_D}(\phi_n ) < \frac{1}{2}$ if $\beta > (\alpha -\frac{1}{2}) \frac{d}{2}$.

\item If $0<\alpha< \frac{1}{2}$, i.e., in the sparse cluster case, since $\CHC_D(\bk, \bn, \lambda)$ is a special case of $\ROHC_D(\bk,\bn, \mu)$ with $\mu = \lambda k^{d/2}$, the computational lower bound in $\CHC_D$ in Theorem \ref{thm: CHC_D comp lower bound} implies that if $\beta > 0$, then $\liminf_{n \to \infty} \cE_{\ROHC_D}(\phi_n) \geq \frac{1}{2}$ based on HPC conjecture \ref{conj: hardness of tensor clique detection}.
\end{itemize}
In summary, we conclude if $\beta > (\alpha -\frac{1}{2}) \frac{d}{2} \vee 0 := \beta^c_{\ROHC_D}$,  any sequence of polynomial-time tests has asymptotic risk at least $1/2$ for $\ROHC_D(\bn, \bk, \mu)$. This has finished the proof of computational lower bound for $\ROHC_D$.

Next we show the computational lower bound for $\ROHC_R$. Suppose there is a sequence of polynomial-time recovery algorithm $\{\phi_R\}_n$ such that $\liminf_{n \to \infty} \cE_{\ROHC_R}(\phi_R) < \frac{1}{2}$ when $\beta > (\alpha -\frac{1}{2}) \frac{d}{2} \lor 0$. In this regime, it is easy to verify $\mu \geq C k^{\frac{d}{4}}$ for some $C > 0$ in $\ROHC_D(\bn, \bk, \mu)$. By Lemma 10 in \cite{luo2020tensorsupp} Section B, 
we know there is a sequence polynomial-time detection algorithms $\{\phi_D\}_n$ such that $\liminf_{n \to \infty}\cE_{\ROHC_D}(\phi_D) < \frac{1}{2}$, which contradicts the computational lower bound established in the first part. This has finished the proof of the computational lower bound for $\ROHC_R$.

\section{Discussion and Future Work}\label{sec: discussion}

In this paper, we study the statistical and computational limits of tensor clustering with planted structures, including the constant high-order structure (CHC) and rank-one high-order structure (ROHC). We derive tight statistical lower bounds and tight computational lower bounds under the HPC/HPDS conjectures for both high-order cluster detection and recovery problems. For each problem, we also provide unconstrained-time algorithms and polynomial-time algorithms that respectively achieve these statistical and computational limits. The main results of this paper are summarized in the phase transition diagrams in Figure \ref{fig: detection and recovery phase transition} and Table \ref{tab: phase transition table}.

There are a few directions worth exploring in the future. First, this paper mainly focuses on the full high-order clustering in the sense that the signal tensor is sparse along all modes. In practice, the partial cluster also commonly appears (e.g., tensor biclustering \citep{feizi2017tensor}), where the signal is sparse only in part of the modes. It is interesting to investigate the statistical and computational limits for high-order partial clustering. Second, in addition to the exact recovery discussed in this paper, we think our results can be extended to other variants of recovery, such as partial recovery and weak recovery.  Third, in the ROHC model, the non-zero components of the signal are required to have the similar magnitudes as this assumption is essential for support recovery. Another interesting problem is to \emph{estimate $(\bv_1,\ldots,\bv_d)$ without the constraint on the component magnitudes of the signal}, which can be seen a rank-one case of the sparse tensor SVD/PCA problem \citep{zhang2019optimal-statsvd,sun2017provable,niles2020all}. For this problem, the signal-to-noise ratio lower bounds we established in Theorems \ref{thm: ROHC detection lower bound} and \ref{thm: ROHC_D and ROHC_R comp lower bound} still hold by virtue of the estimation-to-detection reduction. However, the ROHC$_R$ Search and Power-iteration algorithms studied in this paper may no longer be suitable for estimating $(\bv_1,\ldots,\bv_d)$. A natural unconstrained-time estimator is the maximum likelihood estimator, while to our best knowledge its guarantee is unexplored. \cite{zhang2019optimal-statsvd} developed efficient algorithms which can achieve the minimax optimal error rate in sparse tensor estimation. However, it is unclear if the required signal-to-noise in \cite{zhang2019optimal-statsvd} is tight. It is interesting to develop algorithms with optimal guarantees for sparse tensor SVD/PCA under the tight signal-to-noise ratio requirement. Finally, since our computational lower bounds of $\CHC$ and $\ROHC$ are based on HPC conjecture (Conjecture \ref{conj: hardness of tensor clique detection}) and HPDS conjecture (Conjecture \ref{conj: HPDS recovery conjecture}), it is interesting to provide more evidence for these conjectures.

\section*{Acknowledgement} We would like to thank Guy Bresler for the helpful discussions. We also thank the Editor, the Associated Editor, and two anonymous referees for their helpful suggestions, which helped improve the presentation and quality of this paper. This work was supported in part by NSF Grants CAREER-1944904, NSF DMS-1811868, NSF DMS-2023239, NIH Grant R01 GM131399, and Wisconsin Alumni Research Foundation (WARF).


\bibliographystyle{imsart-nameyear}
\bibliography{reference}

\newpage
\setcounter{page}{1}

\begin{appendix}
	\begin{center}
		{\LARGE Supplement to "Tensor Clustering with Planted Structures: Statistical}
		\medskip
		{\LARGE  Optimality and Computational Limits"	

		}		
		\medskip
		
		{\large Yuetian Luo \quad and\quad Anru R. Zhang}
		\medskip
	\end{center}
	
		 In this supplement, we provide a table of contents
		 and all technical proofs.
\tableofcontents

\bigskip

\section{Proofs of Statistical Bounds for $\CHC_R$ and $\ROHC_R$} \label{sec: proof stat recovery}
We begin by introducing a few notations that will be used throughout the proof sections. Given $\bk = (k_1,\cdots,k_d)$, let $\breve{\bk} = \prod_{i=1}^d k_i$, and $\breve{\bk}^{(-j)} = \prod_{i=1, i \neq j}^d k_i$. For a random variable $\X$, we use notation $\cL(\X)$ to denote the distribution of $\X$. Given a tensor $\bcX \in \bbR^{n_1 \times \cdots \times n_d}$, $\cL( \bcX + N(0,1)^{\otimes n_1 \times \cdots \times n_d})$ denotes the distribution of $\bcX + \bcZ$ where $\bcZ$ is drawn from $N(0,1)^{\otimes n_1 \times \cdots \times n_d}$.

\subsection{Proof of Theorem \ref{thm: CHC and ROHC recovery stat lower bound}}
Since $\CHC_R(\bk, \bn, \lambda)$ is a special case of $\ROHC_R(\bk, \bn, \lambda \sqrt{\prodk})$, the statistical lower bound for $\CHC_R$ implies the lower bound for $\ROHC_R$. So we only need to show the statistical lower bound for $\CHC_R$. Also to show the statistical lower bound of $\CHC_R$, it is enough to show for the case $\bcX = \lambda \mathbf{1}_{I_1} \circ \cdots \circ \mathbf{1}_{I_d}$.

Note that for any $\{\bcX_0, \ldots, \bcX_M \}\subseteq \sX_{\CHC}(\bk, \bn, \lambda)$,
\begin{equation*}
	 \inf_{\phi_R \in \AllAlg^R} \sup_{\bcX \in \sX_{\CHC}(\bk, \bn, \lambda) } \bbP_{\bcX} (\phi_R(\bcY) \neq S(\bcX)) \geq \inf_{\phi_R \in \AllAlg^R} \sup_{\bcX \in \{\bcX_0, \ldots, \bcX_M \} } \bbP_{\bcX} (\phi_R(\bcY) \neq S(\bcX)).
\end{equation*}

Next, we aim to select an appropriate set $\{\bcX_0, \ldots, \bcX_M \}$ and give a lower bound for
$$\inf_{\phi_R \in \AllAlg^R} \sup_{\bcX \in \{\bcX_0, \ldots, \bcX_M \} } \bbP_{\bcX} (\phi_R(\bcY) \neq S(\bcX)).$$ 
We set $\bcX_0 = \lambda 1_{I^{(0)}_1} \circ \cdots \circ 1_{I^{(0)}_d}$, where $I^{(0)}_i = [k_i]$ is the signal support of mode $i$. $\bcX_1, \ldots, \bcX_M$ $(M = n_1 - k_1)$ are constructed in the following way. Assume for mode 2 to mode $d$, $\bcX_1, \ldots, \bcX_{M}$ have the same signal support as $\bcX_0$, but for $\bcX_i$, its signal support on mode 1 is $I^{(i)}_1:=[k_1 - 1] \bigcup \{(i + k_1)\}$.

By such construction, $\bcX_0, \ldots, \bcX_M$ are close to each other. Now we calculate the KL divergence between $\bbP_{\bcX_i}$ and $\bbP_{\bcX_0}$. Observe that $\bbP_{\bcX_i}$ and $\bbP_{\bcX_0}$ are the same except that one of the index of signal support of $\bbP_{\bcX_0}$ on mode 1 is changed from $\{ k_1\}$ to $\{ (i+k_1) \}$,
\begin{equation}
\begin{split}
	\KL(\bbP_{\bcX_i}, \bbP_{\bcX_0}) & = \sum_{i_2 = 1}^{k_2} \cdots \sum_{i_d=1}^{k_d} \KL ( Z_\lambda, Z_0 ) + \sum_{i_2 = 1}^{k_2} \cdots \sum_{i_d=1}^{k_d} \KL ( Z_0, Z_\lambda )\\
	& = \prod_{j=2}^d k_j (\lambda^2/2) + \prod_{j=2}^d k_j (\lambda^2/2)\\
	& = \lambda^2 \prod_{j=2}^d k_j,
\end{split}
\end{equation}
where $Z_\lambda \sim N(\lambda, 1)$ and $Z_0\sim N(0,1)$. So if $\lambda \leq \sqrt{\frac{\eta \log (n_1 - k_1)}{\prod_{j=2}^d k_j} }$, then $\KL(\bbP_{\bcX_i}, \bbP_{\bcX_0}) \leq \eta \log (n_1 - k_1)$. 
So,
\begin{equation}
	\frac{1}{M} \sum_{j=1}^M \KL (\bbP_{\bcX_i}, \bbP_{\bcX_0}) \leq \eta \log M.
\end{equation}

By Theorem 2.5 of \cite{tsybakov2009introduction}, for $0 < \eta < 1/8$, we have
\begin{equation}
	\inf_{\phi_R \in \AllAlg^R} \sup_{\bcX \in \{\bcX_0, \ldots, \bcX_M \} } \bbP_{\bcX} (\phi_R(\bcY) \neq S(\bcX)) \geq \frac{\sqrt{M}}{1+ \sqrt{M}}\left(1- 2\eta - \frac{2\eta}{\log M}\right)\to 1 - 2\eta,
\end{equation} 
where the limit is taken under the asymptotic regime \eqref{assum: asymptotic assumption}.

Notice that in the above construction, we only change the signal support in mode 1. Similarly, we could construct a parameter set which only differs at the signal support at mode $j$, $2\leq j \leq d$. By repeating the argument above, we have if 
\begin{equation}
		\lambda \leq \max\left( \left\{ \sqrt{\frac{\eta \log(n_i - k_i)}{\prod_{j=1, j\neq i}^d k_j}} \right\}_{i=1}^d \right),
	\end{equation} 
	then the minimax estimation error converges to $1- 2\eta$.
	
	Finally, the second part of the conclusion follows by considering the asymptotic regime \eqref{eq: comp lower bound asymptotic regime}. For $\CHC_R$, we have
	\begin{equation*}
	\begin{split}
	\beta > (d-1)\alpha/2 =: \beta^s_{\CHC_R} \Longrightarrow n^{-\beta} \lesssim \sqrt{\frac{1}{n^{\alpha(d-1)}}}  \overset{ \eqref{eq: comp lower bound asymptotic regime} }\Longrightarrow	\lambda \leq \max\left( \left\{ \sqrt{\frac{\eta \log(n_i - k_i)}{\prod_{z=1, z\neq i}^d k_z}} \right\}_{i=1}^d \right). 
	\end{split}
\end{equation*}Similar derivation can be done for ROHC$_R$ as well and this has finished the proof.

\subsection{Proof of Theorem \ref{thm: CHC recovery stat upper bound}}
Denote $\bi = (i_1, \ldots, i_d) $. Given any signal support set $(I_1, \ldots, I_d)$, define $F(I_1, \ldots, I_d) = \sum_{i_1 \in I_1} \cdots \sum_{i_d \in I_d} \bcY_{[i_1, \ldots, i_d]}$. Since the problem only becomes easier as $\lambda$ increases, it is enough to show the algorithm succeeds in the hardest case $\bcX = \lambda \cdot 1_{I_1} \circ \cdots 1_{I_d}$.

Suppose the true signal support of $\bcY \sim \bbP_{\bcX}$ is $(I_1^*, \ldots, I_d^*)$. By construction, the output of Algorithm \ref{alg: CHC recovery combina search} satisfies
\begin{equation} \label{eq: recovery MLE upper bound}
\begin{split}
	\bbP( (\hat{I}_1, \ldots, \hat{I}_d) \neq (I_1^*, \ldots, I_d^*) ) & = \bbP\left( \underset{(\tilde{I}_1, \ldots \tilde{I}_d) \neq (I_1^*, \ldots, I_d^*)}\bigcup \left\{ F(\tilde{I}_1, \ldots, \tilde{I}_d) > F(I_1^*, \ldots, I_d^*) \right\} \right)\\
	& \leq \sum_{i_1=0}^{k_1} \cdots \sum_{i_d=0}^{k_d} P_{i_1, \ldots, i_d} - P_{k_1, \ldots, k_d},
\end{split}
\end{equation}
where 
\begin{equation*}
\begin{split}
	P_{i_1, \ldots, i_d} & = \bbP\left(  F(\tilde{I}_1, \ldots, \tilde{I}_d) > F(I_1^*, \ldots, I_d^*) \right)\\
	& = \left(\prod_{j=1}^d  \binom{k_j}{i_j} \binom{n_j-k_j}{k_j - i_j} \right) \bar{\Phi} \left(\lambda \sqrt{(\breve{\bk} - \breve{\bi} })/2 \right),
\end{split}
\end{equation*} here $\bar{\Phi}( \cdot )$ is the survival function of c.d.f. of a standard Gaussian distribution and $i_j = | \tilde{I}_j \bigcap I_j^*|\, (1\leq j \leq d)$. 

To bound the right hand side of \eqref{eq: recovery MLE upper bound}, we first decompose $\left(\sum_{i_1 =0}^{k_1} \cdots \sum_{i_d=0}^{k_d} P_{i_1, \ldots, i_d} \right)$ into $(d+1)$ different groups. Specifically:
\begin{itemize}
	\item In group 0, $i_j \leq k_j-1$ for all $1 \leq j \leq d$, we denote the summation of terms in this group as $T_0$.
	\item In group 1, exists one $j^* \in [d]$ such that $i_{j^*} = k_{j^*}$ and for $j \neq j^*$, $i_j \leq k_j -1$. We denote the summation of terms in this group as $T_1$.
	\item In group 2, there exists two distinct indices $j_1^*, j_2^* \in [d]$ such that $i_{j_1^*} = k_{j_1^*}$ and $i_{j_2^*} = k_{j_2^*}$. For the rest of indices $j \neq j_1^*, j_2^*$, $i_j \leq k_j -1$. We denote the summation of terms in this group as $T_2$.
	\item Similarly we can define for group $j\, (3 \leq j \leq d-1)$. Denote the summation of terms in group $j$ as $T_j$.
	\item In group d, there is only one term $P_{k_1, \ldots, k_d}$.
\end{itemize}

First notice that the term in group $d$ cancels with the $\left(-P_{k_1, \ldots, k_d}\right)$ in \eqref{eq: recovery MLE upper bound}. Next we are going to give an upper bound for $T_0, T_1, \ldots, T_{d-1}$. Since the strategy to bound each of them is similar, we only need to demonstrate how to bound $T_0$.
\begin{equation}\label{eq: bound of T0}
\begin{split}
	T_0 & \leq  \left(\prod_{j=1}^d (k_j - 1) \right) \underset{\substack{i_z = 0, \ldots, k_z-1\\ z=1, \ldots, d}}\max P_{[i_1, \ldots, i_d ]} \\
	& \leq \left(\prod_{j=1}^d (k_j - 1)\right) \underset{\substack{i_z = 0, \ldots, k_z-1\\ z=1, \ldots, d}}\max \left( \prod_{j=1}^d (n_j - k_j)^{2(k_j -i_j)} \right) \bar{\Phi} \left(\lambda \sqrt{(\breve{\bk} - \breve{\bi} })/2 \right)\\
	& \leq \left(\prod_{j=1}^d (k_j - 1)\right) \underset{\substack{i_z = 0, \ldots, k_z-1\\ z=1, \ldots, d}}\max \left( \prod_{j=1}^d (n_j - k_j)^{2(k_j -i_j)}\right) \exp \, \left( -\frac{\lambda^2}{4} \left( \breve{\bk}- \breve{\bi}  \right)  \right)\\
	& \leq \underset{\substack{i_z = 0, \ldots, k_z-1\\ z=1, \ldots, d}}\max \left(\prod_{j=1}^d (n_j - k_j)^{3(k_j -i_j)}\right) \exp \, \left( -\frac{\lambda^2}{4} \left( \breve{\bk}- \frac{1}{d} \sum_{z=1}^d i_z \breve{\bk}^{(-z)}  \right)  \right),
\end{split}
\end{equation}
where the second inequality is the result of plugging in $P_{i_1, \ldots, i_d}$ and the fact ${k_j \choose i_j} {n_j-k_j \choose k_j - i_j} \leq (n_j - k_j)^{2(k_j -i_j)}$ when $k_i \leq \frac{1}{2}n_i$, the third inequality is by concentration bound for Gaussian random variable $\bbP(Z_0 > t) \leq \frac{1}{t} \exp(-t^2/2)$ and the last inequality is due to $k_i \leq \frac{1}{2}n_i$ and $\frac{1}{d} \sum_{z=1}^d i_z \breve{\bk}^{(-z)} \geq \breve{\bi}$.

The maximum value of right hand side of \eqref{eq: bound of T0} is achieved when $i_1 = k_1 - 1, \ldots, i_d = k_d - 1$ and we have
\begin{equation}
	T_0 \leq \left(\prod_{j=1}^d (n_j - k_j)^3\right) \exp \left(- \frac{\lambda^2}{4d} \sum_{z=1}^d \breve{\bk}^{(-z)} \right).
\end{equation}

So when 
\begin{equation*}
	\lambda^2 \geq C \frac{\sum_{i=1}^d \log (n_i - k_i)}{\min_{1\leq i\leq d} \{\breve{\bk}^{(-i)} \} },
\end{equation*}
for large enough constant $C$ (only depend on $d$), we have 
\begin{equation*}
	T_0 \leq \sum_{i=1}^d (n_i -k_i)^{-c},
\end{equation*}
for some constant $c > 0$. 

Similar analysis holds for $T_j$ ($1 \leq j \leq d-1$). This has finished the proof. 

\subsection{Proof of Theorem \ref{thm: ROHC recovery stat upper bound}}
We first introduce a few notations and the rest of the proof is divided into three steps. Suppose $\bcY \sim \cL( \mu \cdot \bv_1 \circ \cdots \circ \bv_d + N(0,1)^{\otimes n_1 \times \cdots \times n_d})$ and $\bv_i \in \cV_{n_i, k_i}$. First it is easy to check that $\bcA, \bcB$ are independent and have the same distribution $\cL(\frac{\mu}{\sqrt{2}} \cdot \bv_1 \circ \cdots \circ \bv_d + N(0,1)^{\otimes n_1 \times \cdots \times n_d})$. Denote $k_i^* = |\bv_i|$, $\bu_i^* = \mathbf{1}_{S(\bv_i)}$. Since $(\bv_i)_{j} \leq \frac{C}{\sqrt{k_i}}$, the number of non-zero entries in $\bv_i$ is at least $ck_i$ for some small $c > 0$.

{\noindent \bf Step 1}. In this step, we show with high probability the marked pairs in Step 3 of Algorithm \ref{alg: ROHC recovery combinatorial search} are supported on $(S(\bv_1),\ldots, S(\bv_d))$. First for $(\bu_1, \ldots, \bu_d) \in S_{\bar{k}_1}^{n_1} \times \cdots \times S_{\bar{k}_d}^{n_d}$, if one of $\bu_i$ is not supported on corresponding $S(\bv_i)$, we show when $\mu \geq C\sqrt{k \log n}$, such $(\bu_1, \ldots, \bu_d)$ will not be marked in Algorithm \ref{alg: ROHC recovery combinatorial search} Step 3(b) with probability at least $1 - n^{-(d+1)}$. 

Without loss of generality, suppose $S(\bu_1) \nsubseteq S(\bv_1)$ and let $j \in S(\bu_1) \setminus S(\bv_1)$. Notice
\begin{equation*}
    \left(\bcB \times_2 \bu_2^\top \times \cdots \times_d \bu_d^\top\right)_{j} (\bu_1)_{j} \sim N(0, \prod_{i=2}^d \bar{k}_i).
\end{equation*} Then by the Gaussian tail bounds, we have
\begin{equation} \label{ineq: ROHC Recovery null hypothesis uppper bound}
\begin{split}
    &\bbP \left( \left( \bcB \times_2 \bu_2^\top \times \cdots \times_d \bu_d^\top\right)_{j} (\bu_1)_{j} \geq \frac{1}{2\sqrt{2}} \frac{\mu}{\sqrt{\prod_{i=1}^d k_i   }} \prod_{i=2}^d \bar{k}_i \right)\\
    \leq &  \exp \left( -\frac{\mu^2 (\prod_{i=2}^d \bar{k}_i)^2 }{16 \prod_{i=2}^d \bar{k}_i \prod_{i=1}^d k_i }  \right)\\
    \leq &  \exp( -c \frac{\mu^2}{k_1}  ) \leq n^{-(d+1)},
\end{split}
\end{equation} where the last inequality holds because $\mu \geq C\sqrt{k \log n}$ for sufficient large $C > 0$. Similar analysis holds for other modes.

If $(\bu_1, \ldots, \bu_d)$ is marked, let $(\bu_{\bar{k}_1}, \ldots, \bu_{\bar{k}_d}) = (\bu_1, \ldots, \bu_d)$, otherwise let $(\bu_{\bar{k}_1}, \ldots, \bu_{\bar{k}_d}) = (0,\ldots,0)$. Since $\bcB$ is independent of $\bcA$, applying a union bound, we have
\begin{equation}
    \begin{split}
        &\bbP\left(  S(\bu_{\bar{k}_1}) \subseteq S(\bv_1), \ldots, S(\bu_{\bar{k}_d}) \subseteq S(\bv_d) \text{ for all } \bar{k}_i \in [1, k_i], 1\leq i \leq d \right)\\
        \geq & 1- \sum_{ \bar{k}_i \in [1, k_i], 1 \leq i \leq d } \bbP \left( S(\bu_{\bar{k}_i}) \nsubseteq S(\bv_i) \right)\\
        \geq & 1 - (\prodk) n^{-(d+1)} \geq 1 - n^{-1}.
    \end{split}
\end{equation}

{\noindent \bf Step 2}. Let 
\begin{equation}
     (\hat{\bu}_1, \ldots, \hat{\bu}_d) = \arg\max_{(\bu_1, \ldots, \bu_d) \in S_{k_1^*}^{n_1} \times \cdots \times S_{k_d^*}^{n_d}} \bcA \times_1 \bu_1^\top \times_2 \ldots \times_d \bu_d^\top.
\end{equation} In this step, we show that $(S(\hat{\bu}_1), \ldots, S(\hat{\bu}_d)) = (S(\bv_1), \ldots, S(\bv_d)  )$ with high probability. Let $\A_i \in \bbR^{n_i \times n_i}\, (1 \leq i \leq d)$ be a diagonal matrix with its diagonal values
\begin{equation*}
    (\A_i)_{jj} = \left\{ \begin{array}{cc}
        0 & \text{ if } (\bv_i)_j = 0 \\
        1 & \text{ if } (\bv_i)_j > 0 \\
        -1 & \text{ if } (\bv_i)_j < 0. 
    \end{array}\right.
\end{equation*}
So $\bcA$ can be rewritten as $ \frac{\mu}{\sqrt{2}} \cdot \A_1^2\bv_1 \circ \cdots \circ \A_d^2\bv_d + \bcZ $ where $\bcZ \sim N(0,1)^{\otimes n_1 \times \cdots \times n_d}$. So
\begin{equation}\label{eq: reformulate optimization problem}
\begin{split}
    (\hat{\bu}_1, \ldots, \hat{\bu}_d)& = \arg\max_{(\bu_1, \ldots, \bu_d) \in S_{k_1^*}^{n_1} \times \cdots \times S_{k_d^*}^{n_d}} \bcA \times_1 \bu_1^\top \times_2 \ldots \times_d \bu_d^\top\\
    & =\arg \max_{(\bu_1, \ldots, \bu_d) \in S_{k_1^*}^{n_1} \times \cdots \times S_{k_d^*}^{n_d}} \left(\frac{\mu}{\sqrt{2}} \cdot \A_1\bv_1 \circ \cdots \circ \A_d \bv_d \right) \times_1 (\A_1\bu_1)^\top \times \cdots \times_d (\A_d \bu_d)^\top\\
    & \quad \quad \quad \quad \quad \quad \quad  + \bcZ \times_1 \bu_1^\top \times \cdots \times_d \bu_d^\top.
\end{split}
\end{equation}

Notice that in optimization problem \eqref{eq: reformulate optimization problem}, $\A_i \bv_i$ has positive values at its support and these positive values have magnitude at least $\frac{1}{\sqrt{k_i}}$. Also since the diagonal entries of $\A_i$ capture the exact support of $\bv_i$, problem in \eqref{eq: reformulate optimization problem} is a modified version of constant high-order clustering problem with $\lambda = \frac{\mu}{\sqrt{2\prodk}}$ and by a similar argument of Theorem \ref{thm: CHC recovery stat upper bound}, $(S(\hat{\bu}_1), \ldots, S(\hat{\bu}_d)) = (S(\bv_1), \ldots, S(\bv_d)  )$ and sign$(\bu_i) = $sign$(\bv_i)$ with probability at least $1-\sum_{i=1}^d(n_i-k_i)^{-1}$ when $	\lambda = \frac{\mu}{\sqrt{2\prodk}} \geq C_0  \sqrt{\frac{\sum_{j=1}^d \log (n_j - k_j) }{\underset{1\leq i \leq d}\min \{ \prod_{z=1, z\neq i}^d k_z \} }}$ i.e., $\mu \geq C\sqrt{k \log n}$ for some $C > 0$.

{\noindent \bf Step 3}. In this last step, we show $(\bu_1^*, \ldots, \bu_d^*)$ will be marked with high probability. Consider the analysis for mode-1 first. If $j \in S(\bv_1)$,
\begin{equation*}
    \left( \bcB \times_2 \bu_2^{*\top} \times \cdots \times_d \bu_d^{*\top} \right)_{j} (\bu_1^*)_{j} \sim N\left( \frac{\mu}{\sqrt{2}} (\bv_1)_{j} \sum_{i_2=1}^{n_2} \cdots \sum_{i_d=1}^{n_d} (\bv_2)_{[i_2]} \ldots (\bv_d)_{[i_d]} , \prod_{i=2}^d k_i^* \right),
\end{equation*} notice
\begin{equation*}
    \mu_1 := \frac{\mu}{\sqrt{2}} (\bv_1)_{j} \sum_{i_2=1}^{n_2} \cdots \sum_{i_d=1}^{n_d} (\bv_2)_{[i_2]} \ldots (\bv_d)_{[i_d]} \geq \frac{\mu}{\sqrt{2\prodk}}  \prod_{i=2}^d k^*_i.
\end{equation*} By Gaussian tail bound, we have
\begin{equation*}
    \begin{split}
         &\bbP \left( \left( \bcB \times_2 \bu_2^{*\top} \times \cdots \times_d \bu_d^{*\top}\right)_{j} (\bu_1^*)_{j} \leq \frac{1}{2\sqrt{2}} \frac{\mu}{\sqrt{\prodk}} \prod_{i=2}^d k_i^* \right)\\
         \leq & \bbP \left( \left( \bcB \times_2 \bu_2^{*\top} \times \cdots \times_d \bu_d^{*\top} \right)_{j} (\bu_1^*)_{j} - \mu_1 \leq -\frac{1}{2\sqrt{2}} \frac{\mu}{\sqrt{\prodk}} \prod_{i=2}^d k_i^*   \right) \leq n^{-(d+1)},
    \end{split}
\end{equation*} where the last inequality is obtained in the same way as in \eqref{ineq: ROHC Recovery null hypothesis uppper bound}. Similarly, if $j \notin S(\bv_1)$, then $$\left(\bcB \times_2 \bu_2^{*\top} \times \cdots \times_d \bu_d^{*\top}\right)_{j} (\bu_1^*)_{j} \sim N(0, \prod_{i=2}^d k_i^*),$$ and by the same argument of \eqref{ineq: ROHC Recovery null hypothesis uppper bound}, we have 
\begin{equation*}
    \bbP \left( \left( \bcB \times_2 \bu_2^{*\top} \times \cdots \times_d \bu_d^{*\top}\right)_{j} (\bu_1^*)_{j} \geq \frac{1}{2\sqrt{2}} \frac{\mu}{\sqrt{\prodk}} \prod_{i=2}^d k_i^* \right) \leq n^{-(d+1)}.
\end{equation*}
By a union bound, we have 
\begin{equation*}
\begin{split}
     &\bbP \left( S(\bv_1) =  \left\{j: \left(\bcB \times_1 \bu_2^{*\top} \times \cdots \times_d \bu_d^{*\top} \right)_{j} (\hat{\bu}_1)_{j} \geq \frac{1}{2\sqrt{2}} \frac{\mu}{\sqrt{\prodk}} \prod_{i=2}^{d} k^*_i    \right\}  \right)\\
     \geq & 1 - \sum_{j \notin S(\bv_i)} \bbP \left( \left( \bcB \times_2 \bu_2^{*\top} \times \cdots \times_d \bu_d^{*\top}\right)_{j} (\bu_1^*)_{j} \geq \frac{1}{2\sqrt{2}} \frac{\mu}{\sqrt{\prodk}} \prod_{i=2} k_i^* \right)\\
     & - \sum_{j \in S(\bv_i)} \bbP \left( \left( \bcB \times_2 \bu_2^{*\top} \times \cdots \times_d \bu_d^{*\top}\right)_{j} (\bu^*_1)_{j} \leq \frac{1}{2\sqrt{2}} \frac{\mu}{\sqrt{\prodk}} \prod_{i = 2 } k_i^* \right) \\
     \geq & 1  - (n_1 - k_1) n^{-(d+1)}- k_1 n^{-(d+1)}.
\end{split}
\end{equation*} Similar analysis holds for other modes, so a union bound yields that $(\bu_1^*, \ldots, \bu_d^*)$ is marked in Step (3b) with probability at least $1 - n^{-d}$.

Summarize the result so far, with probability at least $1 - \sum_{i=1}^d (n_i-k_i)^{-1} - n^{-d}$, we have
\begin{itemize}
    \item $S(\bu_{\bar{k}_i}) \subseteq S(\bv_i)$ for all $\bar{k}_i \in [1, k_i]$ and all $ 1 \leq i \leq d$.
    \item \begin{equation*}
		        (\bu_1^*, \ldots, \bu^*_d) = \arg \max_{(\bu_1, \ldots, \bu_d) \in S_{k_1^*}^{n_1} \times \cdots \times S_{k_d^*}^{n_d}} \bcA \times_1 \bu_1^\top \times \cdots \times_d \bu_d^\top.
		    \end{equation*}
    \item $(\bu_1^*, \ldots, \bu_d^*)$ is marked in Step (3b).
\end{itemize} With these three points, we conclude with probably at least $1 - \sum_{i=1}^d (n_i-k_i)^{-1} - n^{-d}$, the algorithm can output the true support.

\subsection{Proof of Theorem \ref{thm: CHC and ROHC recovery thresholding}}
The idea to prove this theorem is to use the Gaussian tail bound $\bbP(Z_0 \geq t) \leq \frac{1}{\sqrt{2\pi}}\frac{1}{t} \exp(-\frac{t^2}{2})$, where $Z_0 \sim N(0,1)$. Suppose $I_i = S(\bv_i)\, (1 \leq i \leq d)$. If $(i_1, \ldots, i_d) \notin I_1 \times \cdots \times I_d$, then $\bcY_{[i_1, \ldots, i_d]} \sim N(0,1)$, so 
\begin{equation*}
    \bbP\left(|\bcY_{[i_1, \ldots, i_d]}| \geq \sqrt{2(d+1)\log n} \right) = 2\bbP\left( Z_0 \geq \sqrt{2(d+1)\log n} \right) \leq \frac{2}{\sqrt{2\pi}} n^{-(d+1)} = O(n^{-(d+1)}).
\end{equation*}

For $(i_1, \ldots, i_d) \in I_1 \times \cdots \times I_d$, then $\bcY_{[i_1,\ldots, i_d]} \sim N(\lambda, 1)$ with $\lambda \geq 2\sqrt{2(d+1)\log n}$ for $\CHC_R$ and $\bcY_{[i_1,\ldots,i_d]} \sim N(\mu \cdot (\bv_1)_{i_1} \ldots (\bv_d)_{i_d}, 1 )$ for $\ROHC_R$. Since $(\bv_j)_{i_j} \geq \frac{1}{\sqrt{k_j}}$ by assumption $\bv_j \in \cV_{n_j, k_j}$, $\mu \cdot (\bv_1)_{i_1} \ldots (\bv_d)_{i_d} \geq 2\sqrt{2(d+1)\log n}$. This implies for both $\CHC_R$ and $\ROHC_R$,
\begin{equation*}
    \bbP\left(|\bcY_{[i_1, \ldots, i_d]}| \leq \sqrt{2(d+1)\log n}  \right) \leq \bbP\left( Z_0 \leq -\sqrt{2(d+1)\log n} \right) \leq O(n^{-(d+1)}).
\end{equation*}
So the probability that the set $(i_1, \ldots, i_d)$ with $|\bcY_{[i_1, \ldots, i_d]}| \geq \sqrt{2(d+1)\log n}$ is not exactly $I_1 \times \cdots \times I_d$ is, by union bound, at most
\begin{equation*}
\begin{split}
    &\sum_{(i_1, \ldots, i_d) \notin I_1 \times \cdots \times I_d} \bbP\left( |\bcY_{[i_1, \ldots, i_d]}| \geq \sqrt{2(d+1)\log n} \right) \\
    +& \sum_{(i_1, \ldots, i_d) \in I_1 \times \cdots \times I_d} \bbP\left( |\bcY_{[i_1, \ldots, i_d]}| \leq \sqrt{2(d+1)\log n} \right) = O(n^{-1}),
\end{split}
\end{equation*} which completes the proof of this theorem.

\subsection{Proof of Theorem \ref{thm: CHC and ROHC recovery HOOI upper bound} }
For any tensor $\bcW \in \bbR^{n_1 \times \cdots \times n_d}$,denote $\bcW_{k}^i \in \bbR^{n_1 \times \cdots \times n_{k-1} \times n_{k+1} \times \cdots \times n_d}$ as the subtensor of $\bcW$ by fixing the index of $k^{th}$ mode of $\bcW$ to be $i$ and range over all indices among other modes. 

First, in Algorithm \ref{alg: ROHC and CHC recover via HOOI}, we observe in both models, we have $\bcA \sim \bcX + \widetilde{\bcZ}$ and $\bcB \sim \bcX + \bcZ$, where $\widetilde{\bcZ},\bcZ$ are independent random variable with distribution $N(0,1)^{\otimes n_1 \times \cdots \times n_d}$. In CHC model, $\bcX = \lambda/\sqrt{2} \cdot 1_{I_1} \circ \cdots \circ 1_{I_d}$; in ROHC model, $\bcX = \frac{\mu}{\sqrt{2}} \cdot \bv_1 \circ \cdots \circ \bv_d$.

By the proof of Theorem 1 of \cite{zhang2018tensor}, when
\begin{equation}\label{eq: initialization SNR for HOOI}
    \lambda \sqrt{\breve{\bk}} \geq n^{\frac{d}{4}} \quad \text{ or } \quad \mu \geq Cn^{\frac{d}{4}},
\end{equation} where $\breve{\bk} = \prod_{i=1}^d k_i$, then w.p. at least $1- C\exp(-cn)$, we have
\begin{equation} \label{eq: HOOI recovery singular vector bound}
	\left\| P_{\hat{\bu}_i} - P_{\bu_i} \right\| \leq 2 \left\| \sin \Theta( \hat{\bu}_i, \bu_i ) \right\| \leq C \frac{\sqrt{n_i}}{\lambda \sqrt{ \breve{\bk} }} \text{ or } C\frac{\sqrt{n_i}}{\mu} \quad \text{ for } 1 \leq i \leq d,
\end{equation} here $P_{\bu} = \bu \bu^\top$ denotes the projection operator onto the subspace expanded by $\bu$ and $\|\cdot\|$ is the spectral norm of a matrix. Next, we consider $\CHC_R$ and $\ROHC_R$ separately.

{\noindent \bf For $\CHC_R(\bn, \bk, \lambda)$}. We consider the analysis for mode 1 signal support recovery. First
\begin{equation} \label{eq: error decomp of Y and X}
\begin{split}
&	\left\|  \bcB_1^i \times_2 P_{\hat{\bu}_2} \times \cdots \times_d P_{\hat{\bu}_d} - \bcX_1^i \right\|_{\HS}\\
= & \left\| (\bcB - \bcX)_1^i \times_2 P_{\hat{\bu}_2} \times \cdots \times_d P_{\hat{\bu}_d} + \bcX_1^i  \times_2 P_{\hat{\bu}_2} \times \cdots \times_d P_{\hat{\bu}_d} - \bcX_1^i \right\|_{\HS}\\
\leq&  \left\| \bcZ_1^i  \times_2 P_{\hat{\bu}_2} \times \cdots \times_d P_{\hat{\bu}_d} \right\|_{\HS} + \left\|\bcX_1^i  \times_2 P_{\hat{\bu}_2} \times \cdots \times_d P_{\hat{\bu}_d} - \bcX_1^i \right\|_{\HS}\\
\leq & \left\| \bcZ_1^i  \times_2 P_{\hat{\bu}_2} \times \cdots \times_d P_{\hat{\bu}_d} \right\|_{\HS} + \sum_{k=2}^d \left\| \bcX_1^i \right\|_{\HS} \left\| P_{\bu_k} - P_{\hat{\bu}_k} \right\|,
\end{split}
\end{equation} here $\|\bcA\|_{\HS}=\left(\sum_{i_1,\ldots, i_d} \bcA_{[i_1,\ldots, i_d]}^2\right)^{1/2}$ is the Hilbert-Schmidt norm for tensor $\bcA$. The last inequality is due to triangle inequality and the following decomposition of $\bcX_1^i$,
\begin{equation*}
\begin{split}
	\bcX_1^i =  &\bcX_1^i  \times_2 P_{\hat{\bu}_2} \times \cdots \times_d P_{\hat{\bu}_d} \\
	& + \bcX_1^i  \times_2 P_{\hat{\bu}_2} \times \cdots \times_{d-1} P_{\hat{\bu}_{d-1}} \times_d (\I_{n_d} - P_{\hat{\bu}_d})\\
	& + \bcX_1^i  \times_2 P_{\hat{\bu}_2} \times \cdots \times_{d-2} P_{\hat{\bu}_{d-2}} \times_{d-1} (\I_{n_{d-1}} - P_{\hat{\bu}_{d-1}})\\
	& + \ldots + \bcX_1^i \times_2 (\I_{n_2} - P_{\hat{\bu}_2})\\
	= & \bcX_1^i  \times_2 P_{\hat{\bu}_2} \times \cdots \times_d P_{\hat{\bu}_d} + \sum_{i=2}^d \bcX_1^i \times_{ 1< j < i} P_{\hat{\bu}_j} \times_i \left( P_{\bu_i} -  P_{\hat{\bu}_i} \right).
\end{split}
\end{equation*}
Since $\bcA$ and $\bcB$ are independent, $\bcZ_1^i$ and $P_{\hat{\bu}_2}, \ldots, P_{\hat{\bu}_d}$ are independent. So we have
\begin{equation*}
	\left\| \bcZ_1^i  \times_2 P_{\hat{\bu}_2} \times \cdots \times_d P_{\hat{\bu}_d} \right\|_{\HS} = \left|\bcZ_1^i \times_2 \hat{\bu}_2 \times \cdots \times_d \hat{\bu}_d \right| > \sqrt{2(c+1) \log(n_1)},
\end{equation*}
with probability at most $n_1^{-(c+1)}$. By union bound, we have 
\begin{equation} \label{ineq: upper bound for Z part}
	\bbP \left( \max_{i=1, \ldots, n_1}  \left\| \bcZ_1^i  \times_2 P_{\hat{\bu}_2} \times \cdots \times_d P_{\hat{\bu}_d} \right\|_{\HS} \geq  \sqrt{2(c+1) \log(n_1)} \right) \leq n_1^{-c}.
\end{equation} 

Combining \eqref{ineq: upper bound for Z part} with \eqref{eq: HOOI recovery singular vector bound}, we get an upper bound for \eqref{eq: error decomp of Y and X},
\begin{equation}\label{eq: difference upper bound}
\begin{split}
	\max_i \left\| \bcB_1^i \times_2 P_{\hat{\bu}_2} \times \cdots \times_d P_{\hat{\bu}_d} - \bcX_1^i \right\|_{\HS} & \leq \sqrt{2(c+1) \log(n_1)} + C(d-1) \frac{\sqrt{n_1}}{ \sqrt{ k_1 }},
\end{split}
\end{equation}
with probability at least $1 - n_1^{-c} - C \exp (-c n)$. Here we use the fact $\|\bcX_1^i\|_{\HS} \leq \lambda \sqrt{\prod_{i=2}^d k_i}$.

Then for $i \in I_1, i' \notin I_1$, condition on \eqref{eq: difference upper bound}, we have
\begin{equation}\label{ineq: in cluster and out cluster diff}
	\begin{split}
		& \left| \bcB_1^i \times_2 \hat{\bu}_2^\top  \times \cdots \times_d \hat{\bu}_d^\top - \bcB_1^{i'} \times_2 \hat{\bu}_2^\top  \times \cdots \times_d \hat{\bu}_d^\top \right|\\
	 = 	&\left\| \bcB_1^i \times_2 P_{\hat{\bu}_2} \times \cdots \times_d P_{\hat{\bu}_d} -  \bcB_1^{i'}\times_2 P_{\hat{\bu}_2} \times \cdots \times_d P_{\hat{\bu}_d}  \right\|_{\HS}\\
	 = 	&\left\| \bcB_1^i \times_2 P_{\hat{\bu}_2} \times \cdots \times_d P_{\hat{\bu}_d} - \bcX_1^i + \bcX_1^i - \bcX_1^{i'} + \bcX_1^{i'}  - \bcB_1^{i'}\times_2 P_{\hat{\bu}_2} \times \cdots \times_d P_{\hat{\bu}_d}  \right\|_{\HS}\\
	 \geq 	&- \left\| \bcB_1^i \times_2 P_{\hat{\bu}_2} \times \cdots \times_d P_{\hat{\bu}_d} - \bcX_1^i  \right\|_{\HS} + \left\| \bcX_1^i - \bcX_1^{i'} \right\|_{\HS} - \left\| \bcB_1^{i'} \times_2 P_{\hat{\bu}_2} \times \cdots \times_d P_{\hat{\bu}_d} - \bcX_1^{i'}  \right\|_{\HS}\\
	\geq 	& \frac{\lambda}{\sqrt{2}} \sqrt{\breve{\bk}^{(-1)}} - 2 \left( \sqrt{2(c+1) \log(n_1)} + C(d-1) \frac{\sqrt{n_1}}{ \sqrt{ k_1 }} \right),
	\end{split}
\end{equation} here $\breve{\bk}^{-1} = \prod_i k_i/k_1$.
Similarly, for $i \in I_1, i' \in I_1$ (or $i, i' \notin I_1$ ),
\begin{equation}\label{ineq: two in cluster difference}
	\begin{split}
		& \left| \bcB_1^i \times_2 \hat{\bu}_2^\top  \times \cdots \times_d \hat{\bu}_d^\top - \bcB_1^{i'} \times_2 \hat{\bu}_2^\top  \times \cdots \times_d \hat{\bu}_d^\top \right|\\
	\leq 	& \left\| \bcB_1^i \times_2 P_{\hat{\bu}_2} \times \cdots \times_d P_{\hat{\bu}_d} - \bcX_1^i  \right\|_{\HS} + \left\| \bcX_1^i - \bcX_1^{i'} \right\|_{\HS} + \left\| \bcB_1^{i'} \times_2 P_{\hat{\bu}_2} \times \cdots \times_d P_{\hat{\bu}_d} - \bcX_1^{i'}  \right\|_{\HS}\\
	\leq 	& 2 \left( \sqrt{2(c+1) \log(n_1)} + C(d-1) \frac{\sqrt{n_1}}{ \sqrt{ k_1 }} \right).
	\end{split}
\end{equation}
So when 
\begin{equation*}
	\frac{\lambda}{\sqrt{2}} \sqrt{\breve{\bk}^{(-1)}} \geq 6 \left( \sqrt{2(c+1) \log(n_1)} + C(d-1) \frac{\sqrt{n_1}}{\sqrt{k_1}} \right),
\end{equation*}
we have 
\begin{equation*}
\begin{split}
	& 2 \underset{i, i' \in I_1 \text{ or }i, i' \notin I_1}\max \left| \bcB_1^i \times_2 \hat{\bu}_2^\top  \times \cdots \times_d \hat{\bu}_d^\top - \bcB_1^{i'} \times_2 \hat{\bu}_2^\top  \times \cdots \times_d \hat{\bu}_d^\top \right| \\
	& \leq \underset{i \in I_1, i' \notin I_1 } \min  \left| \bcB_1^i \times_2 \hat{\bu}_2^\top  \times \cdots \times_d \hat{\bu}_d^\top - \bcB_1^{i'} \times_2 \hat{\bu}_2^\top  \times \cdots \times_d \hat{\bu}_d^\top \right|.
\end{split}
\end{equation*}
So a simple cutoff at the maximum gap at the ordered value of $$\left\{ \bcB_1^i \times_2 \hat{\bu}_2^\top  \times \cdots \times_d \hat{\bu}_d^\top \right\}_{i=1}^{n_1}$$ can identify the right support of mode 1 with probability at least $1 - n_1^{-c} - C \exp (-c n)$.

Similar analysis holds for mode $2, 3, \ldots, d$ and combined with SNR requirement in the initialization \eqref{eq: initialization SNR for HOOI}, when
\begin{equation*}
	\lambda \geq C\frac{n^{\frac{d}{4}}}{\sqrt{\breve{\bk}}} \lor 6 \left( \frac{\sqrt{2(c+1) \log(n)}}{\min_{i=1, \ldots, d} \sqrt{\breve{\bk}^{(-i)} }} + C(d-1) \frac{\sqrt{n}}{\sqrt{\breve{\bk}}} \right),
\end{equation*}
then with probability at least $ 1- \sum_{i=1}^d (n_i^{-c} +C \exp(-c n)) $, we have $\hat{I}_i = I_i^*, (1 \leq i \leq d)$. Notice the initialization SNR requirement will be the dominate one when $d \geq 3$, so the final SNR requirement for $\lambda$ is $\lambda \geq C \frac{n^{\frac{d}{4}}}{\sqrt{\breve{\bk}}}$.

\textbf{Next consider $\ROHC_R(\bk, \bn, \mu)$}. The proof of this part is similar to the $\CHC_R$ part.

First, when $\mu \geq C n^{d/4} $, then w.p. at least $1 - C\exp (-cn)$, we have \eqref{eq: HOOI recovery singular vector bound}. The rest of the proof is the same as the first part $\CHC_R$ proof after equation \eqref{eq: HOOI recovery singular vector bound} except that in \eqref{ineq: in cluster and out cluster diff} and \eqref{ineq: two in cluster difference}, we need to consider separately of the support index with positive values and negative values.

Suppose $I_1 = I_{1+} \bigcup I_{1-}$ where $I_{1+}$ denotes the indices in $S(\bv_1)$ that have positive values in $\bv_1$ and $I_{1-}$ denotes the indices in $S(\bv_1)$ that have negative values in $\bv_1$.

By the same argument as \eqref{eq: difference upper bound} and the fact $\bv_i \in \cV_{n_i,k_i}$, with probability at least $1 - n_1^{-c} - C \exp (-c n)$, we have
\begin{equation}\label{ineq: ROHC recovery difference upper bound}
    \max_i \left\| \bcB_1^i \times_2 P_{\hat{\bu}_2} \times \cdots \times_d P_{\hat{\bu}_d} - \bcX_1^i \right\|_{\HS} \leq \sqrt{2(c+1) \log(n_1)} + C(d-1) \frac{\sqrt{n_1}}{ \sqrt{ k_1 }}.
\end{equation}

Similarly, if $i \in I_1, i' \notin I_1$, condition on \eqref{ineq: ROHC recovery difference upper bound}, we have
\begin{equation}\label{ineq: ROHC in cluster and out cluster diff}
	\begin{split}
		& \left| \bcB_1^i \times_2 \hat{\bu}_2^\top  \times \cdots \times_d \hat{\bu}_d^\top - \bcB_1^{i'} \times_2 \hat{\bu}_2^\top  \times \cdots \times_d \hat{\bu}_d^\top \right| \geq  \frac{\mu}{\sqrt{2 k_1}} - 2 \left( \sqrt{2(c+1) \log(n_1)} + C(d-1) \frac{\sqrt{n_1}}{ \sqrt{ k_1 }} \right).
	\end{split}
\end{equation}
For $i \in I_{1+}$ ( $i \in I_{1-}$ or $i \notin I_1$) and $i' \in I_{1+}$ ( $i \in I_{1-}$  or $i \notin I_1$),
\begin{equation}\label{ineq: ROHC two in cluster difference}
	\begin{split}
		& \left| \bcB_1^i \times_2 \hat{\bu}_2^\top  \times \cdots \times_d \hat{\bu}_d^\top - \bcB_1^{i'} \times_2 \hat{\bu}_2^\top  \times \cdots \times_d \hat{\bu}_d^\top \right| \leq 2 \left( \sqrt{2(c+1) \log(n_1)} + C(d-1) \frac{\sqrt{n_1}}{ \sqrt{ k_1 }} \right).
	\end{split}
\end{equation}
For $i \in I_{1+}$ and $i \in I_{1-}$, 
\begin{equation}\label{ineq: ROHC positive and negative cluster difference}
	\begin{split}
		& \left| \bcB_1^i \times_2 \hat{\bu}_2^\top  \times \cdots \times_d \hat{\bu}_d^\top - \bcB_1^{i'} \times_2 \hat{\bu}_2^\top  \times \cdots \times_d \hat{\bu}_d^\top \right| \geq  2\frac{\mu}{\sqrt{2k_1}} - 2 \left( \sqrt{2(c+1) \log(n_1)} + C(d-1) \frac{\sqrt{n_1}}{ \sqrt{ k_1 }} \right).
	\end{split}
\end{equation}
So when \begin{equation*}
	\frac{\mu}{\sqrt{2 k_1}} \geq 6 \left( \sqrt{2(c+1) \log(n_1)} + C(d-1) \frac{\sqrt{n}}{\sqrt{k_1}} \right),
\end{equation*}
a simple two cuts at the top two maximum gaps at the ordered values of $$\left\{ \bcB_1^i \times_2 \hat{\bu}_2^\top  \times \cdots \times_d \hat{\bu}_d^\top \right\}_{i=1}^{n_1}$$ can identify the right support of $I_{1+},I_{1-}$ and $[n_1] \setminus I_1$ with probability at least $1 - n_1^{-c} - C \exp (-c n)$. Since $|I_{1+}|,|I_{1+}| \leq k_1 \ll n_1$, pick two small clusters could yield the right support. Similar analysis holds for other modes, combining with the SNR requirement in initialization \eqref{eq: initialization SNR for HOOI}, when 
\begin{equation*}
	\mu \geq Cn^{\frac{d}{4}}\lor 6 \left( \sqrt{2(c+1) \log(n) k} + C(d-1) \sqrt{n} \right),
\end{equation*}
with probability at least $ 1- \sum_{i=1}^d n_i^{-c} - C \exp(-c n) $, we have $\hat{I}_i = I_i^*, (1 \leq i \leq d)$. Notice the initialization SNR requirement will be the dominate one when $d \geq 3$, the final SNR requirement for $\mu$ is $\mu \geq Cn^{\frac{d}{4}}$. This has finished the proof.

\subsection{Proof of Theorem \ref{thm: CHC recovery aggregated svd} }

Assume $\bcY \sim \bbP_\bcX$ from model \eqref{eq: overall model} and the true support $S(\bcX) = (I_1^*, \ldots, I_d^*)$. We first take a look at mode 1 signal support recovery. Recall $1^*$ is the index such that $1^* = \arg\min_{j \neq 1} n_j$ and
\begin{equation} \label{eq: get Y from mode summation}
	\Y_{[k_1, k_2]}^{(1,1^*)} := \frac{\SUM(\bcY_{k_1, k_2}^{(1,1^*)})}{\sqrt{\prod_{j\neq 1, 1^*} n_j}} \quad  \text{ for } \quad 1 \leq k_1 \leq n_1, 1 \leq k_2 \leq n_{1^*}.
\end{equation}

It is not hard to check that if $\bcY = \lambda \cdot \mathbf{1}_{I_1^*} \circ \cdots \circ \mathbf{1}_{I_d^*} + \bcZ$, then 
\begin{equation*}
	\Y^{(1,1^*)} = \frac{\lambda  \prod_{j\neq 1, 1^*} k_j }{\sqrt{\prod_{j \neq 1, 1^*} n_j}} \mathbf{1}_{I_1^*} \cdot \mathbf{1}_{I_{1^*}^*}^\top + \Z,
\end{equation*}
where $\Z$ has i.i.d. $N(0,1)$ entries.

Now the problem reduces to the submatrix localization problem studied in literature given parameter $$(n_1, n_{1^*}, k_1, k_{1^*}, \frac{\lambda  \prod_{j\neq 1, 1^*} k_j }{\sqrt{\prod_{j \neq 1,1^*}^d n_j}}).$$

By Lemma 1 of \cite{cai2017computational}, if 
\begin{equation*}
	\frac{\lambda \prod_{j \neq 1, 1^*} k_j}{\sqrt{\prod_{j \neq 1, 1^*} n_j}} \sqrt{k_1^*} \geq C \left(  \sqrt{\frac{ n_{1^*}}{k_1}} + \sqrt{\log n_1 } \right),
\end{equation*}
for a large $C > 0$, then w.p. at least $1-n_1^{-c} - C\exp (-c n_1)$, we have $\hat{I}_1 = I_1^*$, and here $c, C > 0$ are some universal constants.

Analysis for other modes are similar and this has finished the proof.

\section{Proofs of Statistical bounds for $\CHC_D$ and $\ROHC_D$}\label{sec: proof-detection}

\subsection{Proof of Theorem \ref{thm: CHC testing statistical lower bound}}
The proof of Theorem \ref{thm: CHC testing statistical lower bound} is fairly long and the main idea is to reduce the minimax testing risk to a Bayesian testing risk with uniform prior over the set of parameters. The main technical difficulty is to bound the second moment of the truncated likelihood ratio, see Lemma \ref{lm: bound on the second moment of likelihood ratio}. For $d = 2$, the lower bound for constant matrix clustering detection has been proved in \cite{butucea2013detection}, however it is much more challenging to show it in order-$d$ case. In this section, we will first prove the main theorem and its subsections are devoted to prove Lemmas used in Theorem.

First, we define some convenient notations. Given a vector $\bx = (x_1, \ldots, x_d) \in \bbR^d$, let $G^{\bn}_{\bx} : = {n_1 \choose x_1} \cdots {n_d \choose x_d} $ and let $\bcY^{sum}_C = \frac{\sum_{(i_1, \ldots, i_d) \in C} \bcY_{[i_1, \ldots, i_d]} }{ \sqrt{k_1 \cdots k_d} }$ for any $C \in \cS_{\bk,\bn}$, here $\cS_{\bk, \bn}$ is the collection of all possible choices of signal locations in the big tensor defined in \eqref{eq: def of S-k-n}. Given $\bcY$, we use notation $\bbP_C:= \bbP_{\bcX}$ where $C = \cS(\bcX)$ to denote the distribution of a tensor with the high-order cluster supported on $C$. Let $\pi$ be a uniform prior on all elements in $\cS_{\bk, \bn}$, i.e.,
\begin{equation*}
	\pi = \left(G^{\bn}_{\bk}\right)^{-1} \sum_{C \in \cS_{\bk, \bn}} \delta_{C},
\end{equation*}
and $\bbP_{\pi}$ be the mixture of distributions $\bbP_\pi = \left(G^{\bn}_{\bk}\right)^{-1} \sum_{C \in \cS_{\bk, \bn}} \bbP_{C}$. Denote the likelihood ratio to be: 
\begin{equation*}
\begin{split}
	\LR_\pi(\bcY) := \frac{d\bbP_\pi}{d\bbP_0} (\bcY) = &\left(G^{\bn}_{\bk}\right)^{-1} \sum_{C \in \mathcal{\cS}_{\bk, \bn}} \exp(-\lambda^2 \breve{\bk}/ 2 + \lambda \sqrt{\breve{\bk}} \bcY^{sum}_C)\\
	= & \left(G^{\bn}_{\bk}\right)^{-1} \sum_{C \in \mathcal{\cS}_{\bk, \bn}} \exp(-b^2/ 2 + b \bcY^{sum}_C).
\end{split}
\end{equation*}
where $ \breve{\bk} = \prod_{i=1}^d k_i$ and $b^2 = \lambda^2 \breve{\bk}$. To show the lower bound, it suffices to show 
\begin{equation} \label{eq: likelihood of Y go to 1}
	\bbP_0\left(\left|\LR_\pi(\bcY) - 1 \right| \geq \epsilon \right) \to 0, \quad \forall \epsilon > 0. 
\end{equation}
Since 
\begin{equation}\label{ineq: risk lower bound}
\begin{split}
	\cE_{\CHC_D}^s  & :=  \inf_{\phi_D \in \AllAlg^D} \left( \bbP_0 (\phi_D(\bcY) = 1) + \sup_{\bcX \in \sX_{\CHC_D}(\bk, \bn, \lambda)} \bbP_\bcX (\phi_D(\bcY) = 0) \right)\\
	& \geq \inf_{\phi_D \in \AllAlg^D} \left( \bbP_0 (\phi_D(\bcY) = 1) + \left(G^{\bn}_{\bk}\right)^{-1} \sum_{C \in \cS_{\bk, \bn}} \bbP_C (\phi_D(\bcY) = 0) \right)\\
	& = \inf_{\phi_D \in \AllAlg^D} \left( \bbE_0 (\phi_D(\bcY)) + \left(G^{\bn}_{\bk}\right)^{-1} \sum_{C \in \cS_{\bk, \bn}} \bbE_0 [(1-\phi_D(\bcY)) \frac{d\bbP_C}{d\bbP_0} (\bcY) ] \right)\\
	& = \inf_{\phi_D \in \AllAlg^D} \left(\bbE_0 (\phi_D(\bcY)) + \bbE_0 [(1-\phi_D(\bcY)) \LR_\pi (\bcY) ] \right)\\
	& \geq \bbE_0 (\phi^*_{D}(\bcY)) + \bbE_0 [(1-\phi^*_{D}(\bcY)) \LR_\pi (\bcY) ],
\end{split}
\end{equation}
where $\phi^*_{D}(\bcY) = 1(\LR_\pi(\bcY) > 1)$ is the likelihood ratio test, take liminf at both side of \eqref{ineq: risk lower bound} and by Fatou's lemma, it is easy to get $\liminf_{n \to \infty} \cE_{\CHC_D}^s \to 1$ if $\LR_\pi(\bcY) \to 1$ in $\bbP_0$ probability.

One canonical way to show $\LR_{\pi}(\bcY) \to 1$ in $\bbP_0$ probability is to show $\bbE_0 (\LR_\pi^2 (\bcY)) \to 1$ and then use chebyshev's inequality. However, the direct calculation of $\bbE_0 (\LR_\pi^2 (\bcY))$ does not work here and we replace $\LR_\pi(\bcY)$ by its truncated version,
\begin{equation} \label{eq: truncated likelihood}
	\widetilde{\LR}_\pi(\bcY) = \left(G^{\bn}_{\bk}\right)^{-1} \sum_{C \in \mathcal{\cS}_{\bk, \bn}} \frac{d \bbP_{C}}{d\bbP_0} (\bcY) 1_{\Gamma_C}.
\end{equation}
Here $\Gamma_C$ is defined as follows: take small $\delta_1 > 0$ (will be specified later) and for $\bv = (v_1, v_2, \ldots, v_d)$, let $\cS_{\bv, C} = \{ V \in \cS_{\bv, \bn}: V \subset C \}$ be the sub-support set of $C$ which are in $\cS_{\bv, \bn}$ and define
\begin{equation*}
	\Gamma_C =\bigcap_{\substack{\delta_1  k_i\leq v_i \leq k_i\\ i=1, \ldots, d}} \underset{V \in \mathcal{\cS}_{\bv, C}}\bigcap \{ \bcY^{sum}_V \leq T_{\bv,\bn}  \},
\end{equation*}
where $T_{\bv, \bn} = \sqrt{2(\log G^\bn_{\bv} + \log \prod_{i=1}^d k_i )}$. It is easy to check, under asymptotic regime \eqref{assum: asymptotic assumption} and $\delta_1 k_i \leq v_i \leq k_i$, we have 
\begin{equation}\label{eq: asymp of T_vb}
	T_{\bv, \bn}^2 \sim 2\left(\sum_{i=1}^d v_i \log(\frac{n_i}{v_i})\right) \sim  2\left(\sum_{i=1}^d v_i \log(\frac{n_i}{k_i})\right),
\end{equation} here $a_n \sim b_n$ if $\lim_{n \to \infty} \frac{a_n}{b_n} = 1$.

Now we introduce the first Lemma.
\begin{Lemma}\label{lm: truncation set is prob 1}
Set $\Gamma_{\bk} = \bigcap_{C \in \mathcal{\cS}_{\bk, \bn}} \Gamma_C$, we have $\bbP_0(\Gamma_{\bk}) = 1$.
\end{Lemma}

This yields $\bbP_0(\LR_\pi(\bcY) = \widetilde{\LR}_{\pi} (\bcY)) \to 1$. So in place of checking \eqref{eq: likelihood of Y go to 1}, it is sufficient to check $\widetilde{\LR}_\pi(\bcY) \to 1$ in $\bbP_0$ probability. To show this, we change it to show the following two Lemmas.

\begin{Lemma}\label{lm: expect of likelihood ratio goes to 0}
	$\bbE_0(\widetilde{\LR}_\pi) \to 1$.
\end{Lemma}

\begin{Lemma}\label{lm: bound on the second moment of likelihood ratio}
	$\bbE_0(\widetilde{\LR}_\pi^2) \leq 1 + o(1)$.
\end{Lemma}

Lemma \ref{lm: expect of likelihood ratio goes to 0} and \ref{lm: bound on the second moment of likelihood ratio} imply that 
\begin{equation}
	\bbE_0 ( \widetilde{\LR}_\pi - 1)^2 = (\bbE_0(\widetilde{\LR}_{\pi}^2) - 1 ) - 2 ( \bbE_0(\widetilde{\LR}_{\pi}) - 1  ) \leq o(1).
\end{equation}
This has finished the proof of Theorem \ref{thm: CHC testing statistical lower bound}.

\subsubsection{Proof of Lemma \ref{lm: truncation set is prob 1}}
It suffices to check $\bbP_0 (\Gamma_\bk^c) \to 0$, where $\Gamma_\bk^c$ is the complement of event $\Gamma_\bk$. First,
\begin{equation}\label{eq: Gamma set complement}
\begin{split}
	\Gamma_\bk^c = & \bigcup_{C \in \cS_{\bk, \bn}} \bigcup_{\substack{\delta_1  k_i\leq v_i \leq k_i\\ i=1, \ldots, d}} \bigcup_{V \in \mathcal{\cS}_{\bv, C}} \{ \bcY^{sum}_V > T_{\bv,\bn}  \}\\
	= & \bigcup_{\substack{\delta_1  k_i\leq v_i \leq k_i\\ i=1, \ldots, d}} \bigcup_{V \in \cS_{\bv, \bn}} \{ \bcY^{sum}_V > T_{\bv,\bn}\}.
\end{split}
\end{equation}
Since $\bcY_V^{sum} \sim N(0, 1)$ under $\bbP_0$, by definition of $T_{\bv, \bn}$ and using the asymptotics $\Phi(-x) \sim e^{-x^2/2} / \sqrt{2\pi}x$ as $x \to \infty$, we have 
\begin{equation*}
	\begin{split}
		\bbP_0( \Gamma_\bk^c ) & \leq \sum_{\substack{\delta_1  k_i\leq v_i \leq k_i\\ i=1, \ldots, d}} \sum_{V \in \cS_{\bv, \bn}} \Phi(-T_{\bv, \bn}) = \sum_{\substack{\delta_1  k_i\leq v_i \leq k_i\\ i=1, \ldots, d}} G^\bn_{\bv} \Phi(-T_{\bv, \bn}) \\
		& \leq \sum_{\substack{\delta_1  k_i\leq v_i \leq k_i\\ i=1, \ldots, d}} \frac{1 + o(1)}{\breve{\bk} T_{\bv, \bn} \sqrt{2\pi} } \to 0.
	\end{split}
\end{equation*}
This has finished the proof.

\subsubsection{Proof of Lemma \ref{lm: expect of likelihood ratio goes to 0} }
	In view of symmetry for $C$ in \eqref{eq: truncated likelihood}, it suffices to check that, for any fixed $C \in \cS_{\bk, \bn}$,
\begin{equation*}
	\bbE_0 \left( \frac{d\bbP_C}{d\bbP_0} 1_{\Gamma_C}\right) = \bbP_C(\Gamma_C) \to 1,
\end{equation*} or equivalently, $\bbP_C \left( \Gamma_C^c \right) \to 0$. Since $\bcY_V^{sum} \sim N(z_{\bv}, 1)$ under the $\bbP_C$, where $z_\bv^2 = \lambda^2 \breve{\bv} := \lambda^2 \prod_{i=1}^d \mathbf{\bv}_{i}$, we have
\begin{equation*}
	\bbP_C(\Gamma_C^c) \leq \sum_{\substack{\delta_1  k_i\leq v_i \leq k_i\\ i=1, \ldots, d}} \sum_{V \in \cS_{\bv, C}} \Phi(z_\bv - T_{\bv, \bn}) = \sum_{\substack{\delta_1  k_i\leq v_i \leq k_i\\ i=1, \ldots, d}}G_{\bv}^{\bk} \Phi(z_\bv - T_{\bv, \bn}),
\end{equation*} where $G_{\bv}^{\bk}= {k_1 \choose v_1} \cdots {k_d \choose v_d}$. By the condition \eqref{cond: CHC detection stat lower bound condition}, there exists $\delta > 0$, 
\begin{equation} \label{ineq: upper bound for b^2}
	b^2 = \lambda^2 \breve{\bk} < (2 - \delta) \sum_{i=1}^d k_i \log \frac{n_i}{k_i}.
\end{equation}
Let $\delta_1$ small enough such that when $\delta_1 k_1 \leq v_1 \leq k_1, \ldots, \delta_1 k_d \leq v_d \leq k_d$, combining with \eqref{eq: asymp of T_vb} we have, 
 \begin{equation*}
	\begin{split}
		z_{\bv}^2 = \lambda^2 \breve{\bv} & < \frac{2-\delta}{ \breve{\bk} }(\sum_{i=1}^d k_i \log\frac{n_i}{k_i}) \breve{\bv} \leq (2-\delta) \sum_{i=1}^d v_i \log \frac{n_i}{k_i}\\
		& \sim (1- \frac{\delta}{2}) T_{\bv, \bn}^2.
	\end{split}
	\end{equation*}
	Thus there exists $\delta > 0$, 
	\begin{equation*}
		\Phi(z_\bv - T_{\bv, \bn}) \leq \exp(-\frac{\delta}{2} T^2_{\bv, \bn}).
	\end{equation*}
	By Stirling's formula, $\log \left(G_{\bv}^\bk\right) = \sum_{i=1}^d \log ({k_i \choose v_i}) \sim \sum_{i=1}^d v_i \log \frac{k_i}{v_i} = O(\sum_{i=1}^d k_i)$, where the last equality is because $\frac{k_i}{v_i} \leq \frac{1}{\delta_1}$. On the other hand, $T_{\bv, \bn}^2 \sim \sum_{i=1}^d v_i \log \frac{n_i}{k_i} \gg \sum_{i=1}^d k_i$ under the asymptotic regime \eqref{assum: asymptotic assumption}, so 
	\begin{equation*}
		\sum_{\substack{\delta_1  k_i\leq v_i \leq k_i\\ i=1, \ldots, d}}G_{\bv}^{\bk} \Phi(z_\bv - T_{\bv, \bn}) \leq \sum_{\substack{\delta_1  k_i\leq v_i \leq k_i\\ i=1, \ldots, d}} \exp( O(\sum_{i=1}^d k_i) - \frac{\delta}{2} T_{\bv, \bn}^2 ) \to 0.
	\end{equation*} 
	This has finished the proof.

\subsubsection{Proof of Lemma \ref{lm: bound on the second moment of likelihood ratio}}
	First,
	\begin{equation*}
		\bbE_0(\widetilde{\LR}_\pi^2) = \left(G_{\bk}^{\bn}\right)^{-2} \sum_{C_1, C_2 \in \mathcal{\cS}_{\bk, \bn}} \bbE_0 ( \exp( -b^2 + b(\bcY^{sum}_{C_1} + \bcY^{sum}_{C_2} ) ) 1_{\Gamma_{C_1} \cap \Gamma_{C_2} } ).
	\end{equation*}
	Denote two latent supports as $C_1 = A_1 \times A_2\times \cdots \times A_d$ and $C_2 = B_1 \times B_2 \times \cdots \times B_d$ where $A_i, B_i \subseteq [n_i]$ and $|A_i| = |B_i| = k_i$ for $1 \leq i \leq d$. Denote the intersection part of $C_1, C_2$ as $V$ and its dimension as $\bv = (v_1, \ldots, v_d)$, i.e., $V = (A_1 \bigcap B_1) \times \cdots \times (A_d \bigcap B_d)$ and $v_1 = |A_1 \bigcap B_1 |, \ldots, v_d = |A_d \bigcap B_d|$.
	
	Notice that the value of $\bbE_0 ( \exp( -b^2 + b(\bcY^{sum}_{C_1} + \bcY^{sum}_{C_2} ) ) 1_{\Gamma_{C_1} \bigcap \Gamma_{C_2} } )$ only depends on the size of $V$. So given $V = C_1 \bigcap C_2$, let 
	\begin{equation}\label{eq: g_v expression}
		g(\bv) := \bbE_0 ( \exp( -b^2 + b(\bcY^{sum}_{C_1} + \bcY^{sum}_{C_2} ) ) 1_{\Gamma_{C_1} \bigcap \Gamma_{C_2} }  ).
	\end{equation} Then \begin{equation*}
	\begin{split}
		\bbE_0(\widetilde{\LR}_\pi^2) & = \sum_{v_1=0}^{k_1} \cdots \sum_{v_d=0}^{k_d} \frac{\# \left( (C_1, C_2) \in \mathcal{\cS}_{\bk,\bn}^2: \text{ size}(V) = (v_1, \ldots, v_d) \right) }{G_{\bk, \bn}^2} g(\bv)\\
		& = \sum_{v_1=0}^{k_1} \cdots \sum_{v_d=0}^{k_d} \left(\prod_{i=1}^d \frac{{k_i \choose v_i} {n_i - k_i \choose k_i - v_i} }{{n_i \choose k_i}} \right) g(\bv)\\
		& = \bbE_{H\cG_1 \times \cdots \times H\cG_{d}} g\left(X_1, \ldots, X_d\right),
	\end{split}
	\end{equation*}
	where $X_i$ independently follows the hypergeometric distribution $H\cG (n_i, k_i, k_i)$.
	
	The goal is to show
	\begin{equation}\label{eq: two moment is order 1}
		\bbE_{H\cG_1 \times \cdots \times H\cG_{d}}\left[ g\left(X_1, \ldots, X_d\right) \right] = 1+ o(1).
	\end{equation}
	Under \eqref{cond: CHC detection stat lower bound condition}, there exists $\delta > 0$, $b^2 = \lambda^2 \breve{\bk} \leq (2 - \delta) \sum_{i=1}^d k_i \log \frac{n_i}{k_i} $, so
	\begin{equation}\label{eq: asymp of lambda square}
		\lambda^2 \asymp \frac{\sum_{i=1}^d k_i \log \frac{n_i}{k_i}}{\breve{\bk} }.
	\end{equation}
	
	To prove \eqref{eq: two moment is order 1}, we consider the value of $\bbE_{H\cG_1 \times \cdots \times H\cG_{d}} \left[g\left(X_1, \ldots, X_d\right) \right]$ on different events and the rest of the proof can be divided into three steps.
	
	{\noindent \bf Step 1}. First let $E_1$ be the value of $\bbE_{H\cG_1 \times \cdots \times H\cG_{d}} \left[g\left(X_1, \ldots, X_d\right) \right]$ on event $\{ \lambda^2 X_1 \cdots X_{d-1} \leq 1 \}$, in this step, we show $E_1 = 1 + o(1)$. Notice
	\begin{equation}\label{eq: hypergeo expect calculate}
	\begin{split}
		E_1 & = \bbE_{H\cG_1 \times \cdots \times H\cG_d} \left[ g(X_1, \ldots, X_d) 1(\lambda^2 X_1 \cdots X_{d-1} \leq 1 ) \right]\\
		& \overset{(a)}\leq \bbE_{H\cG_1 \times \cdots \times H\cG_d} \left[ \exp\left(\lambda^2 X_1 \cdots X_d \right) 1(\lambda^2 X_1 \cdots X_{d-1} \leq 1 ) \right]\\
		& = \bbE_{H\cG_1 \times \cdots \times H\cG_{d-1}} \left[ \bbE_{H\cG_d} \left(\exp\left( \lambda^2 X_1 \cdots X_d \right)\right) 1(\lambda^2 X_1 \cdots X_{d-1} \leq 1 )\right]\\
		& \overset{(b)}\leq \bbE_{H\cG_1 \times \cdots \times H\cG_{d-1}} \left[ \bbE_{Bin_d} \left(\exp( \lambda^2 X_1 \cdots X_d )\right) 1(\lambda^2 X_1 \cdots X_{d-1} \leq 1 )\right]\\
		& \overset{(c)}= \bbE_{H\cG_1 \times \cdots \times H\cG_{d-1}} \left[ \left( 1 + \tilde{q}_d \left( e^{\lambda^2 X_1 \cdots X_{d-1}} - 1 \right) \right)^{k_d} 1(\lambda^2 X_1 \cdots X_{d-1} \leq 1 ) \right]\\
		& \overset{(d)}\leq \bbE_{H\cG_1 \times \cdots \times H\cG_{d-1}} \left[ \exp\left( k_d \tilde{q}_d (e^{\lambda^2 X_1 \cdots X_{d-1}} - 1)  \right) 1(\lambda^2 X_1 \cdots X_{d-1} \leq 1 )\right]\\
		& \overset{(e)}\leq \bbE_{H\cG_1 \times \cdots \times H\cG_{d-1}} \left[ \exp\left( B k_d q_d \lambda^2 X_1 \cdots X_{d-1} \right) \right],
	\end{split}
	\end{equation}
	where $\tilde{q}_d = \frac{k_d}{n_d - k_d}, q_d = \frac{k_d}{n_d}$. (a) is due to Lemma \ref{lm: lm of CHC detection lower bound} \eqref{ineq: first ineq for g_v}, $(b)$ is due to the stochastic dominance of binomial distributed random variable to hypergeometric distributed random variable( see Lemma 5.2 of \cite{butucea2013detection}) and $\Bin_d$ denotes the distribution of a binomial distribution $\Bin(k_d, \frac{k_d}{n_d - k_d})$, (c) is due to the moment generating function of a Binomial distribution, (d) is due to the fact that $(1 + x) \leq e^x$ for any $x$ and (e) is because condition on event $1(\lambda^2 X_1 \cdots X_{d-1} \leq 1 )$, there exists $B >0$ such that $(e^{\lambda^2 X_1 \cdots X_{d-1}} - 1) \leq B\lambda^2 X_1 \cdots X_{d-1}$. 
	
	We can apply the argument in \eqref{eq: hypergeo expect calculate} sequentially and get an upper bound for $E_1$,
	\begin{equation*}
	\begin{split}
		E_1 & \leq \bbE_{H\cG_1 \times \cdots \times H\cG_{d-1}} \left[ \exp\left( B k_d q_d \lambda^2 X_1 \cdots X_{d-1} \right) \right]\\
		& \leq \bbE_{H\cG_1 \times \cdots \times H\cG_{d-2}} \left[ \exp\left( B k_d k_{d-1} q_d q_{d-1} \lambda^2 X_1 \cdots X_{d-2} \right) \right]\\
		& \ldots \\
		& \leq \bbE_{H\cG_1 \times H\cG_{2}} \left[ \exp\left( B k_3 \cdots k_{d-1} k_d q_3 \cdots q_{d-1} q_d \lambda^2 X_1 X_{2} \right) \right]\\
		& \leq \bbE_{H\cG_1} [\exp(Bk_2 \cdots k_d q_2 \cdots q_d \lambda^2 X_1)  ]\\
		& \leq \bbE_{\Bin_1} [ \exp(Bk_2 \cdots k_d q_2 \cdots q_d \lambda^2 X_1) ] \\
		& \leq \left(1 + \tilde{q}_1  ( e^{Bk_2 \cdots k_d q_2 \cdots q_d \lambda^2} - 1 ) \right)^{k_2}\\
		& \leq \exp( Bk_1 k_2 \cdots k_d q_1 q_2 \cdots q_d \lambda^2 )\\
		& = 1+ o(1),
	\end{split}
	\end{equation*}
	where $q_i = \frac{k_i}{n_i}$ and the constant $B$ may vary from line to line. The last equality is due to assumption in \eqref{cond: CHC detection stat lower bound condition} that $\frac{\left(k_1 \cdots k_d\right)^2}{ n_1 \ldots n_d} \lambda^2 \to 0$.
	
	{\noindent \bf Step 2.} To show \eqref{eq: two moment is order 1}, there is left to show \begin{equation}\label{eq: E2 is o(1)}
	    E_2 := \bbE_{H\cG_1 \times \cdots \times H\cG_{d}} \left[g(X_1, \ldots, X_d) 1(\lambda^2 X_1 \cdots X_{d-1} \geq 1) \right] = o(1).
	\end{equation} In this step, we show there exists $\delta_1 > 0$ such that $E_2 = o(1)$ on any of event $\{X_i \leq \delta_1 k_i \}$, $1 \leq i \leq d$.
	
	Denote $\bbP_i(x) : = \bbP_{H\cG_i}(X_i = x) = \frac{{k_i \choose x} {n_i - k_i \choose k_i - x}}{{n_i \choose k_i}}$. Since \eqref{eq: asymp of lambda square} and \eqref{assum: CHC detection lower bound assumption}, we can find small enough $\delta_1$ such that 
	\begin{equation}\label{ineq: delta 1 requirement}
		\delta_1 \lambda^2 \frac{ \breve{\bk}}{k_i} \leq \log(\frac{n_i}{k_i})/2 \quad \text{ for } i = 1, \ldots, d.
	\end{equation}
	So by the definition in \eqref{eq: g_v expression} and \eqref{ineq: first ineq for g_v} of Lemma \ref{lm: lm of CHC detection lower bound}, we have \begin{equation}
	\begin{split}
	 E_2 = & \bbE_{H\cG_1 \times \cdots \times H\cG_d} [ g(X_1, \ldots, X_d) 1(\lambda^2 X_1 \cdots X_{d-1} \geq 1) ] \\
		\leq & \underset{v_1 \cdots v_d \geq \frac{1}{\lambda^2}}{\sum_{v_1 =0}^{k_1} \cdots \sum_{v_d = 0}^{k_d}} \exp(\lambda^2 v_1 \cdots v_d) \prod_{j=1}^d \bbP_j(v_j).
	\end{split}
	\end{equation}
	To analysis the above term, we again consider its value on different events. First we consider the value of $E_2$ on event $W_1 = \{ X_1 \leq \delta_1 k_1 \}$. Notice on $W_1 \bigcap \{ \lambda^2 X_1 \cdots X_{d-1} \geq 1 \}$, for sufficient small $\delta_1$, 
	\begin{equation}
	\begin{split}
		X_2 & \geq \frac{1}{\lambda^2 \delta_1 k_1 X_3\cdots X_{d-1}}\\
		& \geq \frac{1}{\lambda^2 \delta_1 k_1 k_3 \cdots k_{d-1}} \geq \frac{k_2}{\log (\frac{n_2}{k_2})}, 
	\end{split}
	\end{equation} 
	where the least inequality is due to \eqref{eq: asymp of lambda square} and \eqref{assum: CHC detection lower bound assumption}. So by Lemma \ref{lm: lemma 5.3 of Butucea} (with $p = \frac{k_2}{n_2}, r(p) = \log \frac{n_2}{k_2}$), for $v_2 \geq \frac{k_2}{\log (\frac{n_2}{k_2})}$, we have
	\begin{equation}
		\bbP_2(v_2) \leq \exp\left(-v_2 \log\left(\frac{n_2}{k_2}\right)(1+o(1))\right).
	\end{equation}
	So on $W_1$, 
	\begin{equation}\label{ineq: E2 on event W1 goes to 0}
		\begin{split}
			&\bbE_{H\cG_1 \times \cdots \times H\cG_d} \left[ g(X_1, \ldots, X_d) 1(\lambda^2 X_1 \cdots X_{d-1} \geq 1) 1(X_1 \leq \delta_1 k_1) \right] \\
			\leq &\sum_{v_1 \leq \delta_1 k_1} \sum_{v_2 \geq \frac{k_2}{\log (\frac{n_2}{k_2})}}\sum_{v_3 = 1}^{k_3} \cdots \sum_{v_d=1}^{k_d}  \exp(\lambda^2 v_1 \cdots v_d) \prod_{j=1}^d \bbP_j(v_j) \\
			\leq & \sum_{v_1 \leq \delta_1 k_1} \sum_{v_2 \geq \frac{k_2}{\log (\frac{n_2}{k_2})}}\sum_{v_3 = 1}^{k_3} \cdots \sum_{v_d=1}^{k_d}  \exp(\lambda^2 v_1 \cdots v_d) \exp\left(-v_2 \log\left(\frac{n_2}{k_2}\right)\left(1+o(1)\right)\right) \\
			 \leq &\sum_{v_1 \leq \delta_1 k_1} \sum_{v_2 \geq \frac{k_2}{\log (\frac{n_2}{k_2})}}\sum_{v_3 = 1}^{k_3} \cdots \sum_{v_d=1}^{k_d} \exp\left(v_2 (\lambda^2 v_1 v_3\cdots v_d - \log \frac{n_2}{k_2} )\left(1+ o(1)\right) \right)\\
			\leq & \sum_{v_1 \leq \delta_1 k_1} \sum_{v_2 \geq \frac{k_2}{\log (\frac{n_2}{k_2})}}\sum_{v_3 = 1}^{k_3} \cdots \sum_{v_d=1}^{k_d} \exp\left(v_2 (\delta_1 \lambda^2 k_1 k_3\cdots k_d - \log \frac{n_2}{k_2} )\left(1+ o(1)\right) \right) \\
			\leq & \breve{\bk} \exp\left(-\frac{1}{2} v_2 \log \frac{n_2}{k_2} \right) \to 0,
		\end{split}
	\end{equation}
	where the last inequality is due to \eqref{ineq: delta 1 requirement} and the last term goes to 0 is because $v_2 \log \frac{n_2}{k_2} \gg \log( \breve{\bk} )$. 
	
	By the same argument of \eqref{ineq: E2 on event W1 goes to 0}, we can show $E_2$ goes to 0 on any of the event $\{ X_i \leq \delta_1 k_i \}$ for $2 \leq i \leq d$. 
	
	{\noindent \bf Step 3.} In this step, we show the value of $E_2$ is $o(1)$ on event 
	\begin{equation}
		\cH = \left\{(X_1, \ldots, X_d): X_1 \geq \delta_1 k_1, X_2 \geq \delta_1 k_2, \ldots, X_d \geq \delta_1 k_d\right\}.
	\end{equation} Combining this results with the results in Step 2, we have shown \eqref{eq: E2 is o(1)}.
	
	We decompose $\cH$ into two part $\cH_1$ and $\cH_2$ and consider the value of $E_2$ on $\cH_1$ and $\cH_2$ separately. Here $\cH_1$ and $\cH_2$ are defined as
	\begin{equation}
	\begin{split}
		\cH_1 = \left\{(v_1, \ldots, v_d) \in \cH: \sum_{i=1}^d v_1 \log ( \frac{n_i}{k_i} ) \geq 2 \rho_{\bv}\cdot \left(\sum_{i=1}^d k_i \log ( \frac{n_i}{k_i} ) \right)   \right\},\\
		\cH_2 = \left\{(v_1, \ldots, v_d) \in \cH: \sum_{i=1}^d v_1 \log ( \frac{n_i}{k_i} ) < 2 \rho_{\bv}\cdot \left(\sum_{i=1}^d k_i \log ( \frac{n_i}{k_i} ) \right)   \right\},
	\end{split}
	\end{equation}
	where $\rho_{\bv} := \frac{\breve{\bv}}{\breve{\bk}}$. 
	
Denote the value of $E_2$ on $\cH_1$ and $\cH_2$ as $E_{21}$ and $E_{22}$, respectively. We first bound $E_{21}$. 

Observe that $\rho_\bv \geq \delta_1^d$ for $\bv \in \cH$. Recall \eqref{ineq: upper bound for b^2}, observe that we could take $\delta > 0$ small enough such that $t = T_{\bk, \bn} - b(1+ \rho_{\bv}) < 0$. Applying Lemma \ref{lm: lm of CHC detection lower bound} \eqref{ineq: second ineq for g_v} we have
\begin{equation} \label{ineq: bound for E21}
	\begin{split}
		E_{21} & \leq \sum_{\bv \in \cH_1} \exp \left( - ( T_{\bk, \bn} - b)^2 +\frac{\rho_\bv T_{\bk,\bn}^2}{1 + \rho_\bv}  \right) \left( \prod_{i=1}^d \bbP_j(v_j) \right)\\
		& \leq \sum_{\bv \in \cH_1} \exp \left( - ( T_{\bk, \bn} - b)^2  +\frac{\rho_\bv T_{\bk,\bn}^2}{1 + \rho_\bv}- \sum_{i=1}^d v_i \log\frac{n_i}{k_i} + o(T^2_{\bk, \bn}) \right),
	\end{split}
\end{equation} where the second inequality is due to Lemma \ref{lm: lemma 5.3 of Butucea} and observe that $X_i \geq \frac{k_i}{\log \frac{n_i}{k_i}}\, (1 \leq i \leq d)$.

Note that when $\delta > 0$ is small enough, we can take $\delta_2 = \delta_2(\delta) > 0$ such that $(T_{\bk, \bn} - b)^2 \geq \delta_2 T_{\bk, \bn}^2$ for the first term at the right hand side of \eqref{ineq: bound for E21}. Recall $T_{\bk, \bn} \sim 2\left(\sum_{i=1}^d k_i \log \frac{n_i}{k_i}\right)$, so on $\cH_1$, 
\begin{equation*}
	\begin{split}
		& \frac{\rho_\bv T_{\bk,\bn}^2}{1 + \rho_\bv}- \sum_{i=1}^d v_i \log\frac{n_i}{k_i} + o(T^2_{\bk, \bn}) \\
		= &  \frac{2\rho_{\bv}}{1 + \rho_\bv} \left( \sum_{i=1}^d k_i \log ( \frac{n_i}{k_i} ) \right) - \sum_{i=1}^d v_i \log\frac{n_i}{k_i} + o(T^2_{\bk, \bn})\\
		 \leq & \left( \frac{1}{1 + \rho_\bv} -1 \right) \sum_{i=1}^d v_i \log ( \frac{n_i}{k_i} ) + o(T^2_{\bk, \bn}) \leq o(T^2_{\bk, \bn}),
	\end{split}
\end{equation*} here the first inequality is due to the construction of $\cH_1$.
Therefore
\begin{equation*}
	E_{21} \leq \breve{\bk} \exp \left( - (\delta_2+o(1)) T_{\bk, \bn}^2 \right) = o(1).
\end{equation*}

Now we consider $E_{22}$. Recall \eqref{eq: asymp of T_vb}, \eqref{ineq: upper bound for b^2} and $z_{\bv}^2 = \rho_{\bv} \lambda^2 \breve{\bk} = \rho_\bv (2-\delta)(\sum_{i=1}^d k_i \log \frac{n_i}{k_i})$, observe that on $\cH_2$, for small enough $\delta, \delta_1$, we have $T_{\bv, \bn} - 2 z_{\bv} < 0$. So by Lemma \ref{lm: lm of CHC detection lower bound} \eqref{ineq: third ineq for g_v}, we have 
\begin{equation*}
\begin{split}
	E_{22} &\leq \sum_{\bv \in \cH_2} \exp \left( T_{\bv, \bn}^2/2 - (T_{\bv, \bn} - z_{\bv})^2 \right) \prod_{i=1}^d \bbP_j(v_j)\\
	& \leq \sum_{\bv \in \cH_2} \exp \left( T_{\bv, \bn}^2/2 - (T_{\bv, \bn} - z_{\bv})^2  - \left(\sum_{i=1}^d v_i \log \frac{n_i}{k_i} \right)  \right), 
\end{split}
\end{equation*} where the last inequality is by Lemma \ref{lm: lemma 5.3 of Butucea}. 

Since $\sum_{i=1}^d v_i \log \frac{n_i}{k_i} \sim T^2_{\bv, \bn}/2$, the power in the exponent is of the form $$ - (T_{\bv, \bn} - z_{\bv})^2 + o(T_{\bv, \bn}^2). $$
Result $E_{22} = o(1)$ is due to the following result by observing \eqref{ineq: upper bound for b^2} and \eqref{eq: asymp of T_vb}
\begin{equation*}
	\begin{split}
		T^2_{\bv, \bn} - z^2_{\bv} = 2\sum_{i=1}^d v_i\log \frac{n_i}{k_i} - (2-\delta) \sum_{i=1}^d v_i \frac{\prod_{j \neq i} v_i}{\prod_{j \neq i} k_i} \log \frac{n_i}{k_i} \geq \delta' T^2_{\bv, \bn},  
	\end{split}
\end{equation*} for some $\delta' > 0$.

So we have shown $E_{21} = o(1)$, $E_{22} = o(1)$. This has finished the proof.

\begin{Lemma}\label{lm: lm of CHC detection lower bound}
Let $z^2_{\bv} = \lambda^2 \prod_{i=1}^d v_i$.

	(1) \begin{equation}\label{ineq: first ineq for g_v}
	    g(\bv) \leq \bbE_0 \left(\exp(- b^2 + b(\bcY^{sum}_{C_1}  + \bcY^{sum}_{C_2} ) ) \right)= \exp(z_{\bv}^2 ) =: g_1(\bv).
	\end{equation}
	
	(2) If $b \geq \frac{T_{\bk,\bn}}{1 + \rho_{\bv}}$, then 
	\begin{equation} \label{ineq: second ineq for g_v}
	\begin{split}
		g(\bv) & \leq \bbE_0 \left( \exp(-b^2 + b( \bcY^{sum}_{C_1}  + \bcY^{sum}_{C_2} ) ) 1_{\{ \bcY^{sum}_{C_1}\leq T_{\bk, \bn}, \bcY^{sum}_{C_2} \leq T_{\bk, \bn} \}}  \right)\\
		& \leq \exp( - (T_{\bk, \bn} - b)^2 + \frac{\rho_\bv T_{\bk,\bn}^2}{1 + \rho_\bv} ) =: g_2(\bv)
	\end{split}
	\end{equation}
	
	(3) Let $v_1 \geq \delta_1 k_1, \ldots, v_d \geq \delta_1 k_d$ and $T_{\bv, \bn}\leq 2 z_{\bv}$, then 
	\begin{equation}\label{ineq: third ineq for g_v}
	\begin{split}
		g(\bv) \leq & \bbE_0 \left( \exp(-b^2 + b( \bcY^{sum}_{C_1}  + \bcY^{sum}_{C_2} ) ) 1_{\{ \bcY_V^{sum} \leq T_{\bv, \bn} \}}  \right) \\
		= & \exp(T_{\bv,\bn}^2 /2 - (T_{\bv,\bn} - z_{\bv})^2 ) =: g_3(\bv)
	\end{split}
	\end{equation}
\end{Lemma}

\begin{Lemma}[ Lemma 5.3 of \cite{butucea2013detection} ]\label{lm: lemma 5.3 of Butucea}
Let $X_i \sim H\cG(n_i,k_i,k_i)$ and denote $\bbP_i(x) : = \bbP_{H\cG}(X_i = x) = \frac{{k_i \choose x} {n_i - k_i \choose k_i - x}}{{n_i \choose k_i}}$. Let $k_i \to \infty$, $p \to 0$, $v_i \geq k_i/r(p)$ where $r(p) \geq 1$ for $p > 0$ small enough, and $\log(r(p)) = o(log(p^{-1}))$. Then
\begin{equation*}
    \log(\bbP_i(v_i)) \leq v_i \log(p)(1 + o(1)).
\end{equation*}
\end{Lemma}

\subsubsection{Proof of Lemma \ref{lm: lm of CHC detection lower bound}}
The proof is based on the following fact: if $X \sim N(0, 1)$, then $\bbE(\exp(\tau X)) = \exp(\tau^2 /2)$.

We start with the proof of \eqref{ineq: first ineq for g_v}. Given two possible latent supports for the signal $C_1, C_2$, define $V_1 = C_1 \setminus C_2$, $V_2 = C_2 \setminus C_1$, $V = C_1 \bigcap C_2$. Notice $V_1, V_2, V$ are disjoint and $\# (V_1) = \# (V_2) = \breve{\bk} - \breve{\bv}$, $\# (V) = \breve{\bv}$, here $\breve{\bv}:= \prod_{i=1}^d \mathbf{\bv}_{i}$.

Recall $\rho_\bv := \frac{\breve{\bv}}{\breve{\bk}}, b^2 = \lambda^2 \breve{\bk}$, $z_{\bv} = \lambda^2 \breve{\bv} = b^2\rho_{\bv}$ and the definition of $\bcY_C^{sum}$ in the proof of Theorem \ref{thm: CHC testing statistical lower bound}, so
\begin{equation*}
\begin{split}
	\bcY_{C_1}^{sum} & = \sqrt{1- \rho_{\bv}} \bcY_{V_1}^{sum} + \sqrt{\rho_{\bv}} \bcY_{V}^{sum}\\
	\bcY_{C_2}^{sum} & = \sqrt{1- \rho_{\bv}} \bcY_{V_2}^{sum} + \sqrt{\rho_{\bv}} \bcY_{V}^{sum}.
\end{split}
\end{equation*}
So 
\begin{equation*}
\begin{split}
	\bbE_0 \left( \exp(-b^2 + b (\bcY^{sum}_{C_1}  + \bcY^{sum}_{C_2}) ) \right)& = \exp(-b^2)\bbE_0 ( b\sqrt{1- \rho_{\bv}} \bcY_{V_1}^{sum} + b\sqrt{1- \rho_{\bv}} \bcY_{V_2}^{sum} + 2b\sqrt{\rho_{\bv}} \bcY_{V}^{sum} )\\
	& = \exp( b^2 \rho_\bv ) = \exp(z_{\bv}^2).
\end{split}
\end{equation*}

The first inequality of \eqref{ineq: second ineq for g_v} is trivial and we focus on the second one. Let $h \geq 0$, we have 
\begin{equation*}
\begin{split}
	&\bbE_0 \left( \exp(-b^2 + b( \bcY^{sum}_{C_1}  + \bcY^{sum}_{C_2} ) ) 1_{\{ \bcY^{sum}_{C_1}\leq T_{\bk, \bn}, \bcY^{sum}_{C_2} \leq T_{\bk, \bn} \}}  \right)\\
= & \exp( -b^2 + 2 T_{\bk, \bn} h  )\\
& \times \bbE_0 ( \exp\left( (b-h)( \bcY^{sum}_{C_1}  + \bcY^{sum}_{C_2}) + h(\bcY^{sum}_{C_1}  + \bcY^{sum}_{C_2} - 2T_{\bk, \bn}  )  \right) 1_{\{ \bcY^{sum}_{C_1}\leq T_{\bk, \bn}, \bcY^{sum}_{C_2} \leq T_{\bk, \bn} \}}  )\\
\leq & \exp( -b^2 + 2 T_{\bk, \bn} h  ) \bbE_0 \left( \exp\left( (b-h)( \bcY^{sum}_{C_1}  + \bcY^{sum}_{C_2}) \right) \right)\\
= & \exp( -b^2 + 2 T_{\bk, \bn} h  ) \bbE_0 \left( \exp\left( (b-h)\sqrt{1-\rho_{\bv}} ( \bcY^{sum}_{V_1}  + \bcY^{sum}_{V_2}) + 2(b-h) \sqrt{\rho_\bv} \bcY_{V}^{sum} \right) \right) \\
= &  \exp\left( -b^2 + 2 T_{\bk, \bn} h + (b-h)^2(1-\rho_{\bv}) + 2(b-h)^2 \rho_\bv  \right)\\
= &  \exp\left( -b^2 + 2 T_{\bk, \bn} h + (b-h)^2(1+\rho_{\bv}) \right).
\end{split}
\end{equation*}
Plug in $h=  b- \frac{T_{\bk, \bn}}{1 + \rho_{\bv}}$, we get the result.

To prove \eqref{ineq: third ineq for g_v}, for $h \geq 0$, we have
\begin{equation*}
\begin{split}
	& \bbE_0 \left( \exp(-b^2 + b( \bcY^{sum}_{C_1}  + \bcY^{sum}_{C_2} ) ) 1_{\{ \bcY^{sum}_V \leq T_{\bv, \bn} \}}  \right)\\
	= & \exp(-b^2)\bbE_0 ( b\sqrt{1- \rho_{\bv}} \bcY_{V_1}^{sum} + b\sqrt{1- \rho_{\bv}} \bcY_{V_2}^{sum} + 2b\sqrt{\rho_{\bv}} \bcY_{V}^{sum} 1_{\{ \bcY^{sum}_V \leq T_{\bv, \bn} \}} )\\
	= &  \exp(-b^2 \rho_{\bv} + T_{\bv, \bn} h ) \bbE_0 \left( \exp\left( \left(2b\sqrt{\rho_{\bv}}-h\right) \bcY_{V}^{sum} + h(\bcY_{V}^{sum} - T_{\bv, \bn} )  \right) 1_{\{ \bcY^{sum}_V \leq T_{\bv, \bn} \}}  \right)\\
	\leq & \exp(-b^2 \rho_{\bv} + T_{\bv, \bn} h ) \bbE_0 \left( \exp\left( (2b\sqrt{\rho_{\bv}}-h) \bcY_{V}^{sum} \right)  \right)\\
	= & \exp( -b^2 \rho_{\bv} + T_{\bv, \bn} h + (2b\sqrt{\rho_{\bv}}-h)^2/2 ).
\end{split}	
\end{equation*}
Plug in $h = 2b\sqrt{\rho_\bv} - T_{\bv, \bn}$, we get the inequality.

\subsection{Proof of Theorem \ref{thm: CHC testing stat upped bound}}
We adopt the same notation as in the proof of Theorem \ref{thm: CHC testing statistical lower bound}. First note that $$\bbP_0(\psi_{\CHC_D}^s = 1) \leq \bbP_0(\psi_{sum} = 1) + \bbP_0(\psi_{scan} = 1)$$ and $$\bbP_{\bcX} (\psi_{\CHC_D}^s = 0) \leq \min(\bbP_\bcX (\psi_{sum} = 1), \bbP_\bcX (\psi_{scan} = 1) ).$$
So we could consider two tests $\psi_{sum}$ and $\psi_{scan}$ separately.

Let $T_{sum} = \frac{\sum_{i_1=1}^{n_1} \cdots \sum_{i_d=1}^{n_d} \bcY_{[i_1, \ldots, i_d]} }{\sqrt{n_1 \cdots n_d}}$ and denote $Z_\lambda$ to be a random variable following distribution $N(\lambda,1)$. Under $H_0$, $T_{sum} \sim N(0, 1)$, $\bbP_0(\psi_{sum} = 1) = \Phi (-W) \to 0$, since $W \to \infty$ by assumption.
Similarly we have \begin{equation*}
\begin{split}
	\bbP_0 (\psi_{scan} = 1) & = \bbP_0 (T_{scan} > T_{\bk,\bn}) \leq G_{\bk}^{\bn} \bbP_0 (Z_0 > \sqrt{ 2 \log (G^\bn_\bk) })\\
	& \leq G^\bn_{\bk} \cdot \exp(- \log(G^\bn_{\bk}))/ (2 \log (G^\bn_\bk) ) \to 0.
\end{split}
\end{equation*}
Here the second inequality is because $\bbP(Z_0 > t) \leq \frac{1}{t} \exp(-t^2/2)$. Therefore, we have $\bbP_0 (\psi_{\CHC_D}^s = 1) \to 0$.

Under $H_1$, $T_{sum} \sim N(\mu_{sum}, 1)$ with $\mu_{sum} := \frac{\lambda \breve{\bk}}{\sqrt{\breve{\bn}}}$. 

Then if \eqref{cond: sum test condition} holds,
\begin{equation}\label{ineq: linear test alternative}
\begin{split}
	\bbP_\bcX (T_{sum} \leq W) & = \bbP_\bcX(Z_{\mu_{sum} - W} \leq 0 )\\
	& \leq \bbP_\bcX ( Z_{(1-c)\mu_{sum}} \leq 0 ) \to 0,
\end{split}
\end{equation}
where the last inequality is by assumption $W= c \mu_{sum} (c < 1)$ when \eqref{cond: sum test condition} holds.

If \eqref{cond: scan test condition} holds, then
\begin{equation} \label{ineq: scan test guarantee}
\begin{split}
	\bbP_\bcX(T_{scan} \leq T_{\bk,\bn}) & \leq \bbP_\bcX (\bcY^{sum}_{C^*} \leq T_{\bk})\\
	& \leq \bbP_\bcX(Z_{ \lambda \sqrt{ \breve{\bk}} - T_{\bk,\bn} } \leq 0) = \Phi(T_{\bk,\bn} - \lambda \sqrt{\breve{\bk}}),
\end{split}
\end{equation}
where $C^* = \cS(\bcX)$ denote the true latent support of $\bcY$ given $\bcY \sim \bbP_{\bcX}$.

Under assumption \eqref{assum: asymptotic assumption}, by Stirling formula, we have $$T_{\bk,\bn} = \sqrt{2\log(G^\bn_{\bk})} \sim \sqrt{2 \sum_{i=1}^d k_i \log \frac{n_i}{k_i}}.$$
So $T_{\bk,\bn} - \lambda \sqrt{\breve{\bk}} \to - \infty$ and $\Phi(T_{\bk,\bn} - \lambda \sqrt{\breve{\bk}}) \to 0$ under condition \eqref{cond: scan test condition}.
 
Combining \eqref{ineq: linear test alternative} and \eqref{ineq: scan test guarantee}, we have $\bbP_\bcX (\psi_{\CHC_D}^s = 0) \to 0$ holds for every $\bcX$, so $$\sup_{\bcX \in \sX_{\CHC}(\bk, \bn, \lambda)} \bbP_\bcX (\psi_{\CHC_D}^s = 0) \to 0.$$ This has finished the proof.

\subsection{Proof of Theorem \ref{thm: CHC detection upper bound for poly test}}
	Similar to the proof of Theorem \ref{thm: CHC testing stat upped bound}, we only need to consider two tests $\psi_{max}$ and $\psi_{sum}$ separately.
	
	As we have shown in Theorem \ref{thm: CHC testing stat upped bound}, when $W \to \infty$, $\bbP_0(\psi_{sum} = 1) \to 0$. Also $\bbP_\bcX (\psi_{sum} = 0) \to 0, \, \forall \bcX \in \sX_{\CHC}(\bk, \bn, \lambda)$ when condition \eqref{cond: sum test condition} holds. So we only need to show $\bbP_0 (\psi_{max} = 1) \to 0$ and $\bbP_\bcX(\psi_{max} = 0) \to 0 \text{ for } \bcX \in \sX_{\CHC}(\bk, \bn, \lambda)$.
	
	First
	\begin{equation}
	\begin{split}
		\bbP_0 (\psi_{max} = 1) & = \bbP_0 \left(T_{max} \geq \sqrt{2 \sum_{i=1}^d \log n_i}  \right) \\
		& = \breve{\bn} \bbP(Z_0 \geq \sqrt{2 \log \breve{\bn} } )\\
		& \leq  \frac{1}{\sqrt{2 \log \breve{\bn} }} \to 0,
	\end{split}
	\end{equation}
	where the last inequality is due to $\bbP(Z_0 > t) \leq \frac{1}{t} \exp(-t^2/2)$. Also 
	\begin{equation}
	\begin{split}
		\bbP_\bcX (\psi_{max} = 0) & = \bbP_\bcX (T_{max} < \sqrt{2 \log \breve{\bn}}) \\
		& \leq \bbP_\bcX (Z_{\lambda} \leq \sqrt{2 \log \breve{\bn}})\\
		& = \bbP_\bcX (Z_0 \leq \sqrt{2 \log  \breve{\bn}} - \lambda) \to 0,
	\end{split}
	\end{equation}
    where the last term goes to $0$ because $\sqrt{2 \log \breve{\bn}} \to \infty$ and condition \eqref{cond: poly time max test condition}. This has finished the proof.

\subsection{Proof of Theorem \ref{thm: ROHC detection lower bound}}
To prove the lower bound, it is enough to prove for the following case
\begin{equation*}
    H_0: \bcX = 0 \quad \text{ vs } \quad H_1: \bcX = \mu \cdot \bv \circ \bv \circ \cdots \circ \bv,
\end{equation*} where $\bv$ is in the set $\cC$ of $k$-sparse unit vector with non-zero entries equal to $\pm \frac{1}{\sqrt{k}}$. In this case, $\bcY \sim N(0,1)^{\otimes (n^{\otimes d})}$ under $H_0$ and $\bcY \sim \bbP_\bv :=\cL(\mu \cdot \bv \circ \cdots \circ \bv + N(0,1)^{\otimes (n^{\otimes d})})$ under $H_1$.

Let $\bbP_1 = \frac{1}{|\cC|} \sum_{\bv \in \cC} \bbP_{\bv}$, so the likelihood ratio is
\begin{equation*}
    \begin{split}
        \frac{d\bbP_1}{d\bbP_0} &= \exp\left( - \frac{1}{2} \sum_{i_1=1}^n \cdots \sum_{i_d=1}^n \left( \bcY_{[i_1, \ldots, i_d]} - \mu \bv_{i_1} \ldots \bv_{i_d}  \right)^2 + \frac{1}{2} \sum_{i_1=1}^n \cdots \sum_{i_d=1}^n \bcY^2_{[i_1, \ldots, i_d]}    \right)\\
        &= \exp\left( \mu \langle \bcY, \bv\circ \cdots \circ \bv \rangle - \frac{\mu^2}{2} \| \bv\circ \cdots \circ \bv \|_{\HS}^2  \right)\\
      & = \exp\left( \mu \langle \bcY, \bv\circ \cdots \circ \bv \rangle - \frac{\mu^2}{2} \right),
    \end{split}
\end{equation*} where the last inequality is due to the fact that $\|\bv \circ \cdots \circ \bv\|_{\HS} = 1$.

As we illustrate in \eqref{ineq: risk lower bound}, the result follows by showing $\bbE_0 \left( \frac{d\bbP_1}{d\bbP_0}  \right)^2 \to 1$. Notice
\begin{equation*}
    \begin{split}
        \bbE_0 \left( \frac{d\bbP_1}{d\bbP_0}  \right)^2 &= \frac{1}{|\cC|^2} \sum_{\bu, \bv \in \cC} \bbE_0 \left( \frac{d\bbP_\bu}{d\bbP_0} \frac{d\bbP_\bv}{d\bbP_0}   \right)\\
        &= \frac{1}{|\cC|^2} \sum_{\bu, \bv \in \cC} \bbE_0 \exp \left(   \mu \langle \bcY, \bv\circ \cdots \circ \bv + \bu \circ \cdots \circ \bu \rangle - \mu^2 \right)\\
        &= \frac{1}{|\cC|^2} \sum_{\bu, \bv \in \cC} \exp \left( \frac{\mu^2}{2} \| \bu \circ \cdots \circ \bu + \bv \circ \cdots \circ \bv \|_{\HS}^2 - \mu^2  \right)\\
        &= \frac{1}{|\cC|^2} \sum_{\bu, \bv \in \cC} \exp \left( \mu^2 \langle \bu \circ \cdots \circ \bu, \bv \circ \cdots \circ \bv \rangle  \right)\\
        &=\frac{1}{|\cC|^2} \sum_{\bu, \bv \in \cC} \exp\left(  \mu^2 (\bu^\top \bv)^d  \right)\\
        &= \bbE_{\bu, \bv  \sim \text{Uinf}[\cC]} \left( \exp(\mu^2 (\bu^\top \bv)^d ) \right),
    \end{split}
\end{equation*} here the third equality follows from $\bbE(\exp(\langle t, \bcY )) = \exp(\frac{1}{2}\|t\|_{\HS}^2)$. Let $G_m$ denote a symmetric random walk on $\bbZ$ stopped at the $m^{th}$ step. When $\bu, \bv \sim \text{Unif}[\cC]$ are independent, $\bu^\top \bv$ follows distribution of $\frac{G_{H}}{k}$ where $H$ follows Hypergeometric distribution with parameter $(n, k, k)$. Thus $ \mu^2 (\bu^\top \bv)^d = \frac{\mu^2}{k^d} G_H^d = \frac{\mu^2}{k\log\frac{en}{k}} \frac{\log\frac{en}{k}}{k^{d-1}} G_H^d$. Then when $\frac{\mu}{\sqrt{k \log\frac{en}{k}}} \to 0$, by Lemma \ref{lm: MGF bound of power of a sum random walk}, we have 
\begin{equation*}
    \bbE_{\bu, \bv  \sim \text{Uinf}[\cC]} \left( \exp(\mu^2 (\bu^\top \bv)^d ) \right) = \bbE \left[\exp(\frac{\mu^2}{k\log\frac{en}{k}} \frac{\log\frac{en}{k}}{k^{d-1}} G_H^d)\right] \to 1.
\end{equation*}
So this has finished the proof of this theorem.

\begin{Lemma}\label{lm: MGF bound of power of a sum random walk}
Suppose $n \in \mathbb{N}$ and $k \in [n]$. Let $B_1, \ldots, B_k$ be independently Rademacher distributed. Denote the symmetric random walk on $\bbZ$ stopped at the $m^{th}$ step by
\begin{equation*}
    G_m = \sum_{i=1}^m B_i.
\end{equation*} Let $H \sim \text{Hypergeometric}(n,k,k)$ with $\bbP(H = i) = \frac{{k \choose i}{n-k \choose k-i}}{{n \choose k}}, i=1, \ldots, k$. Then there exists a function $g:(0, c) \to (1, \infty)$ with $c > 0$ is a fixed small constant and $g(0+)= 1$, such that for any $a < c$,
\begin{equation}
    \bbE \exp (t G_H^d) \leq g(a),
\end{equation} where $t = \frac{a}{k^{d-1}} \log \frac{en}{k}$.
\end{Lemma}

\subsubsection{Proof of Lemma \ref{lm: MGF bound of power of a sum random walk}}
When $d= 2$, this Lemma was proved in \cite{cai2015optimal}. Their proof heavily relies on the explicit formula of a moment generating function (MGF) of the second moment of Gaussian random variables, which can not be extended to the $d \geq 3$ case as the MGF of high-order moment of Gaussian random variables does not exist. Here, we introduce a new way of proving this lemma. Throughout the proof, we assume $a$ is sufficiently small and $k, n$ are sufficiently large.

Recall that for non-negative random variable $Y$, $\bbE (Y) = \int_{0}^{\infty} \bbP(Y \geq x) dx$. Condition on $H$, $0<\exp(t G_H^d) \leq \exp(tH^d)$, thus
\begin{equation}\label{ineq: upper bound of expectation condition on H}
    \begin{split}
        \bbE \left[ \exp(tG_H^d) | H  \right] & = \int_{0}^{\exp(tH^d)} \bbP \left( \exp(t G_H^d) \geq x  \right) dx \\
        & =\int_{0}^{\exp(tH^d)} \bbP \left( G_H \geq  \left(\frac{\log x}{t}\right)^{\frac{1}{d}}  \right) dx\\
        & \leq \int_{0}^{\exp(tH^d)} \bbP \left( |G_H| \geq  \left(\frac{\log x}{t}\right)^{\frac{1}{d}}  \right) dx\\
        & = \int_{1}^{\exp(tH^d)} \bbP \left( |G_H| \geq  \left(\frac{\log x}{t}\right)^{\frac{1}{d}}  \right) dx + 1\\
        & \overset{(a)}\leq \int_{1}^{\exp(tH^d)} 2 \exp \left( - \frac{  \left(\frac{\log x}{t}\right)^{\frac{2}{d}}   }{2H}   \right) dx + 1.
    \end{split}
\end{equation} Here (a) is due to Hoeffding's inequality.

By taking expectation over $H$, we have 
\begin{equation}\label{eq: integral over h}
    \bbE \left[ \exp(tG_H^d) \right] \leq \sum_{0\leq h \leq k} \bbE \left[ \exp(tG_H^d) | H = h  \right] \bbP(H = h) + 1.
\end{equation}

To prove the result, we only need to show the right hand side of \eqref{eq: integral over h} is upper bounded by $g(a)$, or equivalently we can show that there exists $f(a) > 0$ with $f(0+) = 0$ such that
\begin{equation}\label{eq: small summation part}
    \sum_{0\leq h \leq k }\int_{1}^{\exp(t h^d)} 2 \exp \left( - \frac{  \left(\frac{\log x}{t}\right)^{\frac{2}{d}}   }{2 h}   \right) dx \bbP(H = h) \leq f(a).
\end{equation}
 The idea to prove this is to divide the summation in \eqref{eq: small summation part} into three parts and bound each of them.
 
\textit{Small h}. Assume $h \leq C\left( \frac{k^{d-1}}{\log^2 \frac{en}{k}}  \right)^{\frac{1}{d}}$. In this regime,
\begin{equation}\label{ineq: small h upper bound}
    \begin{split}
         &\sum_{h\in \bbZ:0 \leq h \leq C\left( \frac{k^{d-1}}{\log^2 \frac{en}{k}}  \right)^{\frac{1}{d}}}  \int_{1}^{\exp(th^d)} 2 \exp \left( - \frac{  \left(\frac{\log x}{t}\right)^{\frac{2}{d}}   }{2h}   \right) dx \bbP(H = h)\\
         & \leq \max_{h\in \bbZ:0\leq h \leq C\left( \frac{k^{d-1}}{\log^2 \frac{en}{k}}  \right)^{\frac{1}{d}}}   \int_{1}^{\exp(th^d)} 2 \exp \left( - \frac{  \left(\frac{\log x}{t}\right)^{\frac{2}{d}}   }{2h}   \right) dx \\
         & \overset{(a)}\leq 2\left(\exp\left(Ct \frac{k^{d-1}}{\log^2 \frac{en}{k}}  \right) - 1\right)\\
         & \overset{(b)}\leq \exp(\frac{Ca}{e}) - 1,
    \end{split}
\end{equation} (a) is because the integrand is less or equal to 1 and the fourth inequality is due to the fact $\left(Ct \frac{k^{d-1}}{\log^2 \frac{en}{k}}  \right) = \frac{Ca \log \frac{en}{k} }{k^{d-1}} \frac{k^{d-1}}{\log^2 \frac{en}{k}} = \frac{Ca}{\log \frac{en}{k}} \leq \frac{Ca}{e} $.

\textit{Moderate large h}. Assume $C\left( \frac{k^{d-1}}{\log^2 \frac{en}{k}}  \right)^{\frac{1}{d}} \leq h \leq \frac{k}{a^{\frac{1}{d-1}} \left( \log \frac{en}{k}  \right)^{\frac{1}{d-1}} }$. In this regime,
\begin{equation}\label{ineq: moderate large h}
    \begin{split}
         &\sum_{h\in \bbZ:C\left( \frac{k^{d-1}}{\log^2 \frac{en}{k}}  \right)^{\frac{1}{d}}\leq h \leq \frac{k}{a^{\frac{1}{d-1}} \left( \log \frac{en}{k}  \right)^{\frac{1}{d-1}} } }  \int_{1}^{\exp(th^d)} 2 \exp \left( - \frac{  \left(\frac{\log x}{t}\right)^{\frac{2}{d}}   }{2h}   \right) dx \bbP(H = h) \\
         & \overset{(a)}\leq \sum_{h\in \bbZ:C\left( \frac{k^{d-1}}{\log^2 \frac{en}{k}}  \right)^{\frac{1}{d}}\leq h \leq \frac{k}{a^{\frac{1}{d-1}} \left( \log \frac{en}{k}  \right)^{\frac{1}{d-1}} } } \left(C \left(\exp(th^d - h)+ \exp(-K ) \right)  +  \frac{2}{K-1} (1 - e^{1-K}) \right) \bbP(H = h)\\
         & \overset{(b)}\leq \sum_{h\in \bbZ:C\left( \frac{k^{d-1}}{\log^2 \frac{en}{k}}  \right)^{\frac{1}{d}}\leq h \leq \frac{k}{a^{\frac{1}{d-1}} \left( \log \frac{en}{k}  \right)^{\frac{1}{d-1}} } }  \left(C \left((\exp( - ch) + \exp(-K)\right) + \frac{2}{K-1} (1 - e^{1-K}) \right) \bbP(H = h)\\
         & \overset{(c)}\leq C\exp(-K) + \frac{2}{K-1} (1 - e^{1-K}) + \sum_{h\in \bbZ:C\left( \frac{k^{d-1}}{\log^2 \frac{en}{k}}  \right)^{\frac{1}{d}}\leq h \leq \frac{k}{a^{\frac{1}{d-1}} \left( \log \frac{en}{k}  \right)^{\frac{1}{d-1}} } } C \exp( - ch)\\
         & \leq \exp(a) - 1,
    \end{split}
\end{equation} here the (a) is due to Lemma \ref{lm: upper bound of integral of exp of a log} and $K = \frac{ a^{\frac{2}{d}} k^{\frac{2d - 2}{d}}}{ 2(\log \frac{en}{k})^{\frac{2}{d}} h }$, (b) is due to the fact that when $h \leq \frac{k}{a^{\frac{1}{d-1}} \left( \log \frac{en}{k}  \right)^{\frac{1}{d-1}} }$, we have $t h^d \leq h$, (c) is due to the fact that $\bbP(H = h) \leq 1$ and $\sum_h \bbP(H = h) \leq 1$ and the last inequality holds for sufficiently large $k$.

\textit{Large h}. Assume $\frac{k}{a^{\frac{1}{d-1}} \left( \log \frac{en}{k}  \right)^{\frac{1}{d-1}} }\leq h\leq k$, 
\begin{equation}
    \begin{split}
        &\sum_{h\in \bbZ:\frac{k}{a^{\frac{1}{d-1}} \left( \log \frac{en}{k}  \right)^{\frac{1}{d-1}} }\leq h \leq k }  \int_{1}^{\exp(th^d)} 2 \exp \left( - \frac{  \left(\frac{\log x}{t}\right)^{\frac{2}{d}}   }{2h}   \right) dx \bbP(H = h) \\
        & \overset{(a)}\leq \sum_{h\in \bbZ:\frac{k}{a^{\frac{1}{d-1}} \left( \log \frac{en}{k}  \right)^{\frac{1}{d-1}} }\leq h \leq k }  \left(C \left(\exp(th^d - h)+ \exp(-K ) \right)  +  \frac{2}{K-1} (1 - e^{1-K}) \right) \bbP(H = h)\\
        & \overset{(b)}\leq \sum_{h\in \bbZ:\frac{k}{a^{\frac{1}{d-1}} \left( \log \frac{en}{k}  \right)^{\frac{1}{d-1}} }\leq h \leq k }  \left(C\left( \exp(c th^d ) + \exp(-K) \right) +  \frac{2}{K-1} (1 - e^{1-K}) \right)  \bbP(H = h)\\
        & \overset{(c)}\leq C \left(\exp(-\frac{k^{\frac{d-2}{d}}}{ (\log \frac{en}{k})^{\frac{2}{d}} }) + \frac{2(\log \frac{en}{k})^{\frac{2}{d}}  }{ k^{\frac{d-2}{d}} } \right) + \bbE_{H\sim H\cG} \left[C \exp(c t H^d ) \mathbf{1}(H \geq \frac{k}{a^{\frac{1}{d-1}} \left( \log \frac{en}{k}  \right)^{\frac{1}{d-1}} })\right] \\
        & \overset{(d)}\leq C \left(\exp(-\frac{k^{\frac{d-2}{d}}}{ (\log \frac{en}{k})^{\frac{2}{d}} }) + \frac{2(\log \frac{en}{k})^{\frac{2}{d}}  }{ k^{\frac{d-2}{d}} } \right) +  \bbE_{H \sim \Bin} \left[C \exp(c t H^d ) \mathbf{1}(H \geq \frac{k}{a^{\frac{1}{d-1}} \left( \log \frac{en}{k}  \right)^{\frac{1}{d-1}} })\right]\\
        & = C \left(\exp(-\frac{k^{\frac{d-2}{d}}}{ (\log \frac{en}{k})^{\frac{2}{d}} }) + \frac{2(\log \frac{en}{k})^{\frac{2}{d}}  }{ k^{\frac{d-2}{d}} } \right) \\
        &\quad + \sum_{h\in \bbZ:\frac{k}{a^{\frac{1}{d-1}} \left( \log \frac{en}{k}  \right)^{\frac{1}{d-1}} }\leq h \leq k }  C \exp(c th^d ) \left(\frac{k}{n-k}\right)^h \left(\frac{n-2k}{n-k}\right)^{k-h} {k \choose h}\\
        & \overset{(e)}\leq C \left(\exp(-\frac{k^{\frac{d-2}{d}}}{ (\log \frac{en}{k})^{\frac{2}{d}} }) + \frac{2(\log \frac{en}{k})^{\frac{2}{d}}  }{ k^{\frac{d-2}{d}} } \right) +  \sum_{h\in \bbZ:\frac{k}{a^{\frac{1}{d-1}} \left( \log \frac{en}{k}  \right)^{\frac{1}{d-1}} }\leq h \leq k } \exp(c t k^{d-1} h - h \log \frac{n}{2k} + h\log\frac{ek}{h} )\\
        & \overset{(f)}\leq  \sum_{h\in \bbZ:\frac{k}{a^{\frac{1}{d-1}} \left( \log \frac{en}{k}  \right)^{\frac{1}{d-1}} }\leq h \leq k } \exp\left(ca h \log \frac{en}{k}  - h \log \frac{n}{2k} + h \log \left(  e a^{\frac{1}{d-1}} \left( \log \frac{en}{k}  \right)^{\frac{1}{d-1}}  \right) \right)\\
        &\quad  +  C \left(\exp(-\frac{k^{\frac{d-2}{d}}}{ (\log \frac{en}{k})^{\frac{2}{d}} }) + \frac{2(\log \frac{en}{k})^{\frac{2}{d}}  }{ k^{\frac{d-2}{d}} } \right)\\
        & \overset{(g)}\leq \exp(a) - 1,
    \end{split}
\end{equation} here (a) is due to Lemma \ref{lm: upper bound of integral of exp of a log} with $K = \frac{ a^{\frac{2}{d}} k^{\frac{2d-2}{d}}}{2(\log \frac{en}{k} )^{\frac{2}{d}} h }$, (b) is due to the fact that when $h \geq a^{\frac{1}{d-1}} \left( \log \frac{en}{k}  \right)^{\frac{1}{d-1}}$, $t h^d \geq h$, (c) is due to the fact that $$C\exp(-K) + \frac{2}{K-1} (1 - e^{1-K}) \leq C \left(\exp \left(-\frac{k^{\frac{d-2}{d}}}{ (\log \frac{en}{k})^{\frac{2}{d}} }\right) + \frac{2(\log \frac{en}{k})^{\frac{2}{d}}  }{ k^{\frac{d-2}{d}} } \right),$$ (d) is due to stochastic dominance of Hypergeometric distribution Hypergeometric$(n,k,k)$ by Binomial distribution $\Bin(k, \frac{k}{n-k})$ (e.g., see Lemma 5.2 of \cite{butucea2013detection}), (e) is due to that ${k \choose h} \leq \left(\frac{ek}{h}\right)^h$ and $n \geq 3k$ and $\frac{n-2k}{n-k} < 1$, (f) is due to that in the summation range of $h$, $\frac{k}{h} \leq a^{\frac{1}{d-1}} \left( \log \frac{en}{k}  \right)^{\frac{1}{d-1}}$ and (g) holds for a sufficient small $a$ and sufficiently large $k$.

So we can take $f(a) = 2(\exp(a)-1) + \exp(\frac{Ca}{e}) - 1$, i.e., $g(a) = 2(\exp(a)-1) + \exp(\frac{Ca}{e})$ in the statement and this has finished the proof.

\begin{Lemma}\label{lm: upper bound of integral of exp of a log}
Suppose $M$, $K$ are sufficiently large fixed constants and set $\tau = \frac{2}{d}$, then for $d \geq 3$, we have
\begin{equation*}
    \int_{1}^{\exp(M)} 2 \exp(-K \log^\tau x) dx \leq \frac{2}{K-1} (1 - e^{1-K}) +  C\left(\exp(-K) + \exp(-M^\tau K + M)\right),
\end{equation*} for some $C > 0$.
\end{Lemma}

\subsubsection{Proof of Lemma \ref{lm: upper bound of integral of exp of a log}}
The idea to prove this lemma is to divide the integral range into different parts and then bound each part separately.

\begin{equation}\label{ineq: upper bound of integral of explog}
    \begin{split}
        \int_{1}^{\exp(M)} 2 \exp(-K \log^\tau x) dx &= \int_1^e 2 \exp\left(-K \log^\tau x\right) dx + \sum_{i=1}^{M^\tau - 1} \int_{\exp(i^{\frac{1}{\tau}})}^{\exp((i+1)^{\frac{1}{\tau}})} 2 \exp(-K \log^\tau x) dx \\
        & \overset{(a)}\leq \int_1^e 2 \exp\left(-K \log^\tau x\right) dx +  2\sum_{i=1}^{M^\tau - 1} \exp(-i K) \exp\left((i+1)^{\frac{1}{\tau}} \right),
    \end{split}
\end{equation} where (a) is due to the fact that the integrand decays as $x$ increases.

We first bound the first term in the right hand side of \eqref{ineq: upper bound of integral of explog}. Since $\log x \in [0,1]$ when $x \in [1, e]$ and $\tau < 1$, we have
\begin{equation*}
\begin{split}
    &\int_1^e 2 \exp\left(-K \log^\tau x\right) dx \leq \int_1^e 2 \exp(-K \log x) dx \\
    = & \int_1^e 2 x^{-K} dx = \frac{2}{K-1} (1 - e^{1-K}).
\end{split}
\end{equation*}
Now we move onto the second term in \eqref{ineq: upper bound of integral of explog}. We show that the value of the sequence of $\exp(-i K) \exp\left((i+1)^{\frac{1}{\tau}} \right)$ is a U-shape curve, specifically, it first geometrically decreases and then geometrically increases, so the summation of the sequence could be upper bounded by the sum of first term and the last term.

Compare two consecutive values in the summation $\sum_{i=1}^{M^\tau - 1} \exp(-i K) \exp\left((i+1)^{\frac{1}{\tau}} \right)$, when $i$ is large,
\begin{equation*}
\begin{split}
    \frac{\exp(-i K) \exp\left((i+1)^{\frac{1}{\tau}} \right)}{\exp(-(i-1) K) \exp\left(i^{\frac{1}{\tau}} \right)} & = \exp(-K) \exp\left(i^{\frac{1}{\tau}} \left[ (1+\frac{1}{i})^{\frac{1}{\tau}} - 1 \right] \right) \\
    & = \exp(-K) \exp\left( C \frac{1}{\tau} i^{\frac{1}{\tau} - 1}  \right),
\end{split}
\end{equation*} for some $C > 0$, where the last equation is due to Taylor expansion.

So the sequence $\exp(-i K) \exp\left((i+1)^{\frac{1}{\tau}} \right)$ geometrically decreases when $i \lesssim K^{\frac{\tau}{1-\tau}}$ and then geometrically increases when $i \gtrsim K^{\frac{\tau}{1-\tau}}$. At the same time, the beginning terms in the series $\exp(-i K) \exp\left((i+1)^{\frac{1}{\tau}} \right)$ are of order $\exp(-K)$. 

So by the summation property of geometric sequence, we have
\begin{equation*}
    \sum_{i=1}^{M^\tau - 1} \exp(-i K) \exp\left((i+1)^{\frac{1}{\tau}} \right) \leq C\left(\exp(-K) + \exp(-M^\tau K + M)\right) ,
\end{equation*} for some $C > 0$. Combining two parts, we finish the proof of this Lemma.

\subsection{Proof of Theorem \ref{thm: ROHC detection stat upper bound}}
Notice that by Theorem \ref{thm: ROHC recovery stat upper bound}, when $\mu \geq C \sqrt{k \log n}$, the $\ROHC_R$ search Algorithm \ref{alg: ROHC recovery combinatorial search} can identify the true support with probability goes to 1. Then by Lemma \ref{lm: ROHS recovery imply detection}, $\cE_{\ROHC_D}(\psi_{\ROHC_D}^s) \to 0$ follows by observing that $\psi_{\ROHC_D}^s$ is the algorithm used in Lemma \ref{lm: ROHS recovery imply detection}. Combining them together, we get the theorem.

\subsection{Proof of Theorem \ref{thm: ROHC detection poly test upper bound}}
First in regime $\lim_{n \to \infty} \frac{\mu}{\sqrt{2(\prodk)(\sum_{i=1}^d \log n_i)}  } > 1$, by the same argument as Theorem \ref{thm: CHC detection upper bound for poly test}, we can show $\bbP_0(\psi_{max} = 1) \to 1$ and $\bbP_{\bcX}(\psi_{max} = 0) \to 0$.

In regime $\mu \geq C n^{\frac{d}{4}}$, notice that by Theorem \ref{thm: CHC and ROHC recovery HOOI upper bound}, the Power-iteration algorithm can identify the true support with probability goes to 1. Also notice that $\mu \geq C k^{\frac{d}{4}}$ in this case, then by Lemma \ref{lm: ROHS recovery imply detection}, $\cE_{\ROHC_D}(\psi_{sing}) \to 0$ follows by observing that $\psi_{sing}$ is the algorithm used in Lemma \ref{lm: ROHS recovery imply detection}. Combining them together, we get the theorem.

\begin{Lemma}\label{lm: ROHS recovery imply detection}
Consider $\ROHC_D(\bk, \bn, \mu)$ and $\ROHC_R(\bk, \bn,\mu)$ under the asymptotic regime \eqref{assum: asymptotic assumption}. 

(1) If $\mu \geq C \sqrt{k \log n}$ for some $C > 0$ and there is a sequence of recovery algorithm $\{ \phi_R \}_n$ such that $\liminf_{n \to \infty} \cE_{\ROHC_R}(\phi_R) <  \eta$ for $\eta \in (0,1]$, then there exists a sequence of test $\{\phi_D \}_n$ such that $\liminf_{n \to \infty} \cE_{\ROHC_D}(\phi_D) < \eta$.

(2) If $\mu \geq Ck^{\frac{d}{4}}$ for some $C > 0$ and there is a sequence of polynomial-time recovery algorithm $\{ \phi_R \}_n$ such that $\liminf_{n \to \infty} \cE_{\ROHC_R}(\phi_R) <  \eta$ for $\eta \in (0,1]$, then there exists a sequence of polynomial-time algorithm $\{\phi_D \}_n$ such that $\liminf_{n \to \infty} \cE_{\ROHC_D}(\phi_D) < \eta$. 
\end{Lemma}

\subsubsection{Proof of Lemma \ref{lm: ROHS recovery imply detection}}
The idea to prove this lemma is to form two independent copies of the observation from the original observation, then use the first copy and the recovery algorithm $\phi_R$ to find the true support with high probability, finally use the estimated support and the second copy to do test.

First by the property of Gaussian, it is easy to check $\bcA := \frac{\bcY + \bcZ_1}{\sqrt{2}}$ and $\bcB:=\frac{\bcY - \bcZ_1}{\sqrt{2}}$ are two independent copies with distribution $\cL \left( \frac{\mu}{\sqrt{2}} \cdot\bv_1 \circ \cdots \circ \bv_d + N(0,1)^{\otimes n_1 \times \cdots \times n_d} \right)$ if $\bcZ_1 \sim N(0,1)^{\otimes n_1 \times \cdots \times n_d}$.

{\noindent \bf Proof of Statement (1)}. Based on $\bcA$ and assumption about algorithm $\phi_R$, we have event $E = \{ \phi_R(\bcY) = S(\bv_1) \times \cdots \times S(\bv_d) \}$ happens with probability more than $ 1- \eta$. Condition on $E$, denote $\widetilde{\bcB}$ to the part of $\bcB$ that restricted to the support $S(\bv_1) \times \cdots \times S(\bv_d)$.

Apply $\ROHC_R$ Search Algorithm \ref{alg: ROHC recovery combinatorial search} to $\widetilde{\bcB}$ and in the output of Step 4, we get $(\hat{\bu}_1, \ldots, \hat{\bu}_d)$ and let $\hat{\bv}_1 = \frac{\hat{\bu}_1}{\sqrt{k_1}}, \ldots, \hat{\bv}_d = \frac{\hat{\bu}_d}{\sqrt{k_d}}$. The test procedure is $\phi_D = \mathbf{1}\left(\widetilde{\bcB} \times_1 \hat{\bv}_1^\top \times \cdots \times_d^\top \hat{\bv}_d \geq C \sqrt{k} \right)$. 

Under $H_0$, $\widetilde{\bcB}$ has i.i.d. $N(0,1)$ entries, by Lemma 5 of \cite{zhang2018tensor}, we have
\begin{equation*}
    \bbP( |\widetilde{\bcB} \times_1 \hat{\bv}_1^\top \times \cdots \times_d \hat{\bv}_d^\top| \geq C\sqrt{k} ) \leq \exp(-ck),
\end{equation*} for some $C, c > 0$.

Under $H_1$ and regime $\mu \geq C \sqrt{k \log n}$, by Theorem \ref{thm: ROHC recovery stat upper bound}, we have the elementwise sign of $\hat{\bu}_i$ matches $\bv_i$ and $\hat{\bv}_i^\top \bv_i \geq c$ for some $c > 0$ with probability at least $1-O(n^{-1})$. So conditioned on $E$, with probability at least $1 - n^{-1}$,
\begin{equation}
\begin{split}
    \widetilde{\bcB} \times_1 \hat{\bv}_1^\top \times \cdots \times_d^\top \hat{\bv}_d &= \frac{\mu}{\sqrt{2}} \times \prod_{i=1}^d (\hat{\bv}_i^\top \bv) \geq \frac{\mu}{\sqrt{2}} c^d \geq C\sqrt{k},
\end{split}
\end{equation} So overall we have 
\begin{equation}
    \cE_{\ROHC_D}(\phi_D) \leq  C\exp(-ck) + Cn^{-1} \to 0,
\end{equation}
condition on $E$. On $E^c$ we have $\cE_{\ROHC_D}(\phi_D) < 1$. So the results follow since $E$ happens with probability more than $1 - \eta$.

{\noindent \bf Proof of Statement (2)}. Similarly based on $\bcA$ and assumption about polynomial-time algorithm $\phi_R$, we have event $E = \{ \phi_R(\bcY) = S(\bv_1) \times \cdots \times S(\bv_d) \}$ happens with probability more than $ 1- \eta$. Define $\widetilde{\bcB}$ in the same way as before.

Apply Power-iteration (Algorithm \ref{alg: ROHC and CHC recover via HOOI}) on $\widetilde{\bcB}$ and get the rank-one approximation $\hat{\bv}_1, \ldots, \hat{\bv}_d$ at the Step 5. The test procedure is $\phi_D = \mathbf{1}\left(\widetilde{\bcB} \times_1 \hat{\bv}_1^\top \times \cdots \times_d^\top \hat{\bv}_d \geq C \sqrt{k} \right)$. 

Under $H_0$, similarly we have
\begin{equation*}
    \bbP( |\widetilde{\bcB} \times_1 \hat{\bv}_1^\top \times \cdots \times_d \hat{\bv}_d^\top| \geq C\sqrt{k} ) \leq \exp(-ck),
\end{equation*} for some $C, c > 0$.

Under $H_1$ and regime $\mu \geq C k^{\frac{d}{4}}$, by Theorem 1 of \cite{zhang2018tensor}, we have $\| \sin \Theta( \hat{\bv}_i, \bv_i ) \|_2 \leq \frac{\sqrt{k}}{\mu} \leq \frac{1}{2}\, (1 \leq i \leq d)$ w.p. at least $1- C\exp(-ck)$. Here $\|\sin\Theta(\hat{\bv}_i, \bv_i)\| = \sqrt{1 - (\hat{\bv}_i^\top \bv_i)^2}$ is used measure the angle between $\hat{\bv}_i$ and $\bv_i$. So $\hat{\bv}_i^\top \bv_i \geq \sqrt{1 - \frac{1}{4}} = \frac{\sqrt{3}}{2}$. Conditioned on $E$, 
\begin{equation}
\begin{split}
    \widetilde{\bcB} \times_1 \hat{\bv}_1^\top \times \cdots \times_d^\top \hat{\bv}_d &= \frac{\mu}{\sqrt{2}} \times \prod_{i=1}^d (\hat{\bv}_i^\top \bv) \geq \frac{\mu}{\sqrt{2}} (\frac{\sqrt{3}}{2})^d \geq Ck^{\frac{d}{4}} \geq C\sqrt{k},
\end{split}
\end{equation} with probability at least $1- C\exp(-ck)$ for some $c, C > 0$. So overall we have \begin{equation}
    \cE_{\ROHC_D}(\phi_D) \leq  C\exp(-ck) \to 0,
\end{equation}
condition on $E$. On $E^c$, $\cE_{\ROHC_D}(\phi_D) < 1$ for sufficient large $n$. So the results follow since event $E$ happens with probability more than $1 - \eta$ and it is polynomial-time computable since both $\phi_R$ and Power-iteration are polynomial-time computable.

\section{Average-Case Reduction} \label{sec: average case reduction}

We first list some existing tools in the literature and introduce some new results for establishing computational lower bounds.

The following lemma gives the rigorous results why mapping the conjecturally hard problem in distribution to the target problem is the key step in the idea of average-case reduction we have built in Section \ref{sec: ROHC_D and ROHC_R comp lower bound}.
\begin{Lemma}[Lemma 4 of \cite{brennan2018reducibility}] \label{lm: average case reduction simultanious mapping}
Let $\cP$ and $\cP'$ be detection problems with hypotheses $H_0, H_1, H_0', H_1'$ and let $X$ and $Y$ be instances of $\cP$ and $\cP'$, respectively. Suppose there is a polynomial-time computable map $\varphi$ satisfying
\begin{equation*}
    \TV\left( \cL_{H_0} (\varphi(X)), \cL_{H_0'}(Y) \right) + \sup_{\bbP \in H_1} \inf_{\pi \in \triangle(H_1')} \TV \left( \cL_\bbP( \varphi(X) ), \int_{H_1'} \cL_{\bbP'} (Y) d\pi(\bbP')   \right) \leq \delta,
\end{equation*} 
where $\triangle(H_1')$ denotes the set of priors on $H_1'$. If there is a polynomial-time algorithm solving $\cP'$ with Type I + II error at most $\epsilon$, then there is a polynomial-time algorithm solving $\cP$ with Type I+II error at most $\epsilon + \delta$.
\end{Lemma}

\begin{Lemma}[Data Processing \citep{csiszar1967information}] \label{lm: data processing}
Let $P$ and $Q$ be distributions on a measurable space $(\Omega, \mathcal{B})$ and let $f: \Omega \to \Omega'$ be a Markov transition kernel (see definition in \cite[Section 8.3]{klenke2013probability}). If $A \sim P$ and $B \sim Q$, then 
\begin{equation*}
    \TV\left( \cL(f(A)), \cL(f(B)) \right) \leq \TV(P, Q).
\end{equation*}
Here, $\cL(\cdot)$ is the distribution of any random variable $``\cdot"$.
\end{Lemma}

\begin{Lemma}[Tensorization (\cite{ma2015computational}, Lemma 7)] \label{lm:tensorization}
Let $P_1, \ldots, P_n$ and $Q_1, \ldots, Q_n$ be distributions on a measurable space $(\Omega, \mathcal{B})$. Then 
\begin{equation*}
    \TV\left(\prod_{i=1}^n P_i, \prod_{i=1}^n Q_i\right) \leq \sum_{i=1}^n \TV(P_i, Q_i).
\end{equation*}
\end{Lemma}

The following Lemma \ref{lm: diagonal permutation} will be used in proving the computational lower bound of constant high-order cluster recovery. 
Given an order-$d$ tensor $\bcW$, let $\bcW^{\sigma_1, \ldots, \sigma_d}$ be the tensor formed by permuting mode $i$ indices by permutation $\sigma_i$, i.e. 
$$\bcW_{[\sigma_1(i_1), \ldots, \sigma_d(i_d)] }^{\sigma_1, \ldots, \sigma_d} = \bcW_{[i_1, \ldots, i_d]}.$$

\begin{Lemma}\label{lm: diagonal permutation}
Let $P$ and $Q$ be two distributions such that $Q$ dominates $P$, (i.e. for any event $A$, $Q(A) = 0$ implies $P(A) = 0$) and $\chi^2(P, Q) \leq 1$. Suppose $\bcW \in \bbR^{n^{\otimes d}}$ is an order-$d$ dimension-$n$ tensor with all its non-diagonal entries i.i.d. sampled from $Q$ and all of its diagonal entries i.i.d. sampled from $P$, where the set of diagonal entries of $\bcW$ is $\{\bcW_{i, i, \ldots, i}\}$. Suppose that $\sigma_1, \ldots\, \sigma_{d-1}$ are independent permutations on $[n]$ chosen uniformly at random.  Then 
\begin{equation*}
    \TV \left(\cL (\bcW^{\id, \sigma_1, \ldots, \sigma_{d-1}}), Q^{\otimes (n^{\otimes d})}   \right) \leq \sqrt{ \chi^2 (P, Q) },
\end{equation*}
where 
$$\chi^2(P, Q) = \int \frac{(P(x) - Q(x))^2}{Q(x)}dx$$
is the $\chi^2$ divergence between distributions $P$ and $Q$, and ``$\id$" be the identity permutation.
\end{Lemma}

The proof of this lemma is given in Section \ref{sec: average case lemmas}. Next, we introduce the rejection kernel algorithm.
\begin{algorithm}[h]\caption{Rejection Kernel}\label{alg: Rejection Kernel}
\begin{algorithmic}[1]
\State {\bf Input:} $x \in \{0,1\}$, a pair of PMFs or PDFs $f_X$ and $g_X$ that can be efficiently computed and sampled, Bernoulli probabilities $p, q \in [0,1]$, number of iterations $T$.
	\State Initialize $Y = 0$.
	\State For $i = 1, \ldots, T$, do:
	\begin{enumerate}[label=(\alph*)]
	    \item If $x = 0$, sample $Z \sim g_X$ and if
	    \begin{equation*}
	        p \cdot g_X(Z) \geq q \cdot f_X(Z),
	    \end{equation*}
	    then with probability $1 - \frac{q \cdot f_X (Z)}{p\cdot g_X(Z)}$, update $Y = Z$ and break.
	    \item If $x = 1$, sample $Z \sim f_X$ and if 
	    \begin{equation*}
	        (1 - q) \cdot f_X(Z) \geq (1-p) \cdot g_X(Z)
	    \end{equation*} 
	    then with probability $1 - \frac{(1-p) \cdot g_X (Z)}{(1-q)\cdot f_X(Z)}$, update $Y = Z$ and break.
	\end{enumerate}
	\State {\bf Output}: $Y$.
\end{algorithmic}	
\end{algorithm}

 Recall we denote the above Rejection Kernel map as $\RK(p \to f_X, q \to g_X, T)$. The following lemma discusses the mapping from Bernoulli random variable to Gaussian random variable by rejection kernel. We omit the proof of Lemma \ref{lm: Rejection kernel for CHC recovery} here since the proof is essentially the same as the proof of Lemma 14 in \cite{brennan2018reducibility} except for some constant modifications.
\begin{Lemma} \label{lm: Rejection kernel for CHC recovery}
Let $n$ be a parameter and suppose that $p = p(n)$ and $q = q(n)$ satisfy $p > q$, $p, q \in [0, 1]$, $\max(q, 1-q) = \Omega(1)$ and $p - q \geq n^{-O(1)}$. Let $\delta = \min \{ \log \frac{p}{q}, \log \frac{1-q}{1-p}  \}$. Suppose $\xi = \xi(n) \in (0,1)$ satisfies 
$$\xi \leq \frac{\delta}{2 \sqrt{2(d+1)\log n + 2 \log\left((p - q)^{-1}\right)}}.$$ 
Then the map $\RK_G := \RK(p \to N(\xi, 1), q \to N(0,1), T)$ with $T = \lceil 2(d+1) \delta^{-1} \log n \rceil$ 
can be computed in poly$(n)$ time and satisfies
$$\TV\left( \RK_G(\Bern(p)), N(\xi,1) \right) = O(n^{-(d+1)})\quad \text{and}\quad  \TV(\RK_G(  \Bern(q), N(0,1))) = O(n^{-(d+1)}).$$
\end{Lemma}

The next lemma gives the property of tensor reflection cloning in Algorithm \ref{alg: Tensor Reflecting Clone}.
\begin{Lemma}[Tensor Reflecting Cloning]\label{lm: tensor reflecting cloning}
Suppose $n$ is even and $\ell = O(\log n)$. There is a randomized polynomial-time computable map $\varphi: \bbR^{n^{\otimes d}} \to \bbR^{n^{\otimes d}}$ given by Algorithm \ref{alg: Tensor Reflecting Clone} that satisfies
\begin{enumerate}[leftmargin=*]
    \item $\varphi\left(N(0,1)^{\otimes (n^{\otimes d})} \right) \sim N(0,1)^{\otimes (n^{\otimes d})}$.
    \item Consider any $\lambda > 0$ and any set of vectors $\bu_1, \ldots, \bu_d \in \bbZ^n$. There exists a distribution $\pi$ over vectors $\bu_1^{(\ell)}, \ldots, \bu_d^{(\ell)} \in \bbZ^{n}$ with $\| \bu_i^{(\ell)} \|_2^2 = 2^{\ell} \| \bu_i \|_2^2$, $2^\ell \| \bu_i \|_0 \geq \|\bu_i^{(\ell)}\|_0$ for $1 \leq i \leq d$ such that 
    \begin{equation*}
        \varphi\left(\lambda \cdot \bu_1 \circ \cdots \circ \bu_d + N(0,1)^{\otimes (n^{\otimes d})} \right) \sim  \int \cL \left( \frac{\lambda}{\sqrt{2}^{d\ell}}\cdot \bu_1^{(\ell)} \circ \cdots \circ \bu_d^{(\ell)} + \bcZ \right) d\pi (\bu_1^{(\ell)}, \ldots, \bu_d^{(\ell)}),
    \end{equation*} 
    where $\bcZ \sim N(0,1)^{\otimes (n^{\otimes d})}$. Furthermore, it holds with probability at least $1- 4(\sum_{i=1}^d \|\bu_i\|^{-1}_0 )$ over $\pi$ such that for $i \in [d]$, 
    \begin{equation*}
        2^\ell \|\bu_i\|_0 \geq \|\bu^{(\ell)}_i\|_0 \geq 2^\ell \|\bu_i\|_0 \left( 1- \max( \frac{2C \ell \log(2^{\ell} \|\bu_i\|_0 )}{\|\bu_i\|_0}, \frac{2^{\ell} \|\bu_i\|_0}{n} ) \right),
    \end{equation*} for some $C > 0$ if $\|\bu_i\|_0$ are sufficiently large and at most $2^{-\ell-1} n$ for $i \in [d]$. Finally, if $\bu_i = \bu_j$ are equal for some $i, j \in [d], i \neq j$ then $\bu_i^{(\ell)} = \bu_j^{(\ell)}$ holds almost surely.
\end{enumerate}
\end{Lemma}

\subsection{Proof of Lemma \ref{lm: diagonal permutation}} \label{sec: average case lemmas}
Let $\sigma_1', \ldots, \sigma_{d-1}'$ be independent permutations of $[n]$ chosen uniformly at random and also independent of $\sigma_1, \ldots, \sigma_{d-1}$. For convenience, we denote $\bsigma(i) = (\sigma_1(i), \ldots, \sigma_{d-1}(i)), \bsigma'(i) = (\sigma'_1(i), \ldots, \sigma'_{d-1}(i))$ and $\bj = (j_1, \ldots, j_{d-1})$. By property $\chi^2(P,Q) + 1 = \int \frac{P^2(x)}{Q(x)} dx$, we have
\begin{equation} \label{eq: chi-square divergence def1}
\begin{split}
    &\chi^2 \left( \cL(\bcW^{\id, \sigma_1, \ldots, \sigma_{d-1}}), Q^{\otimes (n^{\otimes d})}  \right) + 1 \\
    & = \int \frac{ \left( \bbE_{\sigma_1, \ldots, \sigma_{d-1}} \left[ \bbP_{\bcW^{\id, \sigma_1, \ldots, \sigma_{d-1}}} (\bcX | \sigma_1, \ldots, \sigma_{d-1}) \right] \right)^2 }{\bbP_{Q^{\otimes (n^{\otimes d})}}(\bcX)} d\bcX\\
    & = \bbE_{\sigma_1, \ldots ,\sigma_{d-1}, \sigma'_{1}, \ldots, \sigma'_{d-1} } \int \frac{ \bbP_{\bcW^{\id, \sigma_1, \ldots, \sigma_{d-1}}} (\bcX | \sigma_1, \ldots, \sigma_{d-1}) \bbP_{\bcW^{\id, \sigma'_1, \ldots, \sigma'_{d-1}}} (\bcX | \sigma'_1, \ldots, \sigma'_{d-1}) }{ \bbP_{Q^{\otimes (n^{\otimes d})}} (\bcX) } d\bcX.
\end{split}
\end{equation}
Notice 
\begin{equation*}
    \bbP_{\bcW^{\id, \sigma_1, \ldots, \sigma_{d-1}}} (\bcX | \sigma_1, \ldots, \sigma_{d-1}) = \prod_{i=1}^n \left\{  P(\bcX_{[i, \sigma_1(i), \ldots,\sigma_{d-1}(i) ]}) \prod_{\bj \neq \bsigma(i) } Q(\bcX_{ [i, j_1, \ldots, j_{d-1}] })  \right\},
\end{equation*}
so 
\begin{equation}\label{eq: chi-square distance decomposition}
\begin{split}
      & \frac{ \bbP_{\bcW^{\id, \sigma_1, \ldots, \sigma_{d-1}}} (\bcX | \sigma_1, \ldots, \sigma_{d-1}) \bbP_{\bcW^{\id, \sigma'_1, \ldots, \sigma'_{d-1}}} (\bcX | \sigma'_1, \ldots, \sigma'_{d-1}) }{ \bbP_{Q^{\otimes (n^{\otimes d})}} (\bcX) }\\
    =& \prod_{i: \bsigma(i) = \bsigma'(i) } \left\{  \frac{ P^2(\bcX_{[i, \sigma_1(i), \ldots,\sigma_{d-1}(i) ]}) }{ Q(\bcX_{ [i, \sigma_1(i), \ldots,\sigma_{d-1}(i) ] }) }  \prod_{\bj \neq \bsigma(i) } Q(\bcX_{ [i, j_1, \ldots, j_{d-1}] }) \right\}\\
    & \times \prod_{i: \bsigma(i) \neq \bsigma'(i) } \frac{P(\bcX_{[i, \sigma_1(i), \ldots,\sigma_{d-1}(i) ]}) P(\bcX_{[i, \sigma'_1(i), \ldots,\sigma'_{d-1}(i) ]}) \prod_{ \substack{\bj \neq \bsigma(i) \\ \bj \neq \bsigma'(i) }} Q^2( \bcX_{ [i, j_1, \ldots, j_{d-1}] } ) }{ \prod_{ \substack{\bj \neq \bsigma(i) \\ \bj \neq \bsigma'(i) }} Q(\bcX_{ [i, j_1, \ldots, j_{d-1}] })  }\\
    =&  \prod_{i: \bsigma(i) = \bsigma'(i) } \frac{ P^2(\bcX_{[i, \sigma_1(i), \ldots,\sigma_{d-1}(i) ]}) }{ Q(\bcX_{ [i, \sigma_1(i), \ldots,\sigma_{d-1}(i) ] }) } \prod_{i: \bsigma(i) \neq \bsigma'(i) } P(\bcX_{[i, \sigma_1(i), \ldots,\sigma_{d-1}(i) ]}) P(\bcX_{[i, \sigma'_1(i), \ldots,\sigma'_{d-1}(i) ]})  \\
    & \times \left( \prod_{ \substack{\bj \neq \bsigma(i) \\ \bj \neq \bsigma'(i) }} Q( \bcX_{ [i, j_1, \ldots, j_{d-1}] } )   \right).
\end{split}
\end{equation}
After integration, last two terms at the right hand side of \eqref{eq: chi-square distance decomposition} are integrated to $1$, only the first term left. So
\begin{equation} \label{eq: chi-square divergence simplify}
\begin{split}
    & \int \frac{ \bbP_{\bcW^{\id, \sigma_1, \ldots, \sigma_{d-1}}} (\bcX | \sigma_1, \ldots, \sigma_{d-1}) \bbP_{\bcW^{\id, \sigma'_1, \ldots, \sigma'_{d-1}}} (\bcX | \sigma'_1, \ldots, \sigma'_{d-1}) }{ \bbP_{Q^{\otimes (n^{\otimes d})}} (\bcX) } d\bcX \\
    =  &\prod_{i: \bsigma(i) = \bsigma'(i) } \left(  \int \frac{ P^2(\bcX_{[i, \sigma_1(i), \ldots,\sigma_{d-1}(i) ]}) }{ Q(\bcX_{ [i, \sigma_1(i), \ldots,\sigma_{d-1}(i) ] }) } d \bcX_{ [i, \sigma_1(i), \ldots,\sigma_{d-1}(i) ] } \right)\\
    =  &\left( 1 + \chi^2(P, Q)  \right)^{|\{i: \bsigma(i) = \bsigma'(i)\}|}.
\end{split}
\end{equation}
Let $Y = |\{i: \bsigma (i) = \bsigma'(i)\}|$ be the number of fixed coordinates in all $d-1$ permutations and let $\bar{Y} = |\{i: \sigma_1(i) = \sigma'_1(i) \}|$ be the number of fixed coordinate in the first permutation, clearly $Y \leq \bar{Y}$. As shown in \cite{pitman1997some}, the $i$th moment of $\bar{Y}$ is at most the $i$th Bell
number and possion distribution with rate $1$  has its $i$th moment given by $i$th Bell number for all $i$. So the moment generating function(m.g.f.) $\bbE(e^{t\bar{Y}})$ is at most the m.g.f. of possion distribution with rate $1$ which is $\exp(e^t -1)$. Set $t = \log (1 + \chi^2(P, Q)  )$,
\begin{equation*}
\begin{split}
    \chi^2 \left( \cL(\bcW^{\id, \sigma_1, \ldots, \sigma_{d-1}}), Q^{\otimes (n^{\otimes d})}  \right) &\overset{(a)}= \bbE[ \left( 1 + \chi^2(P, Q)  \right)^{Y} ] - 1\\
    & \leq \bbE[ \left( 1 + \chi^2(P, Q)  \right)^{\bar{Y}} ] -1\\
    & = \bbE [\exp(\bar{Y}\log(1 + \chi^2(P,Q) ))] -1 \\
    &  \leq \exp\left( \chi^2(P, Q) \right) - 1 \leq 2\cdot \chi^2(P, Q).
\end{split}
\end{equation*} Here (a) is due to \eqref{eq: chi-square divergence def1}\eqref{eq: chi-square divergence simplify} and the last inequality is because $e^x \leq 1 + 2x$ for $x \in [0,1]$. Finally, since $\TV(P, Q) \leq \sqrt{\frac{\chi^2(P, Q)}{2}}$ by \cite{tsybakov2009introduction} Lemma 2.7, we have
\begin{equation*}
    \TV \left(  \cL (\bcW^{\id, \sigma_1, \ldots, \sigma_{d-1}}), Q^{\otimes (n^{\otimes d})}   \right) \leq \sqrt{ \frac{ \chi^2 \left( \cL(\bcW^{\id, \sigma_1, \ldots, \sigma_{d-1}}), Q^{\otimes (n^{\otimes d})}  \right) }{2} } \leq \sqrt{ \chi^2 (P, Q) }.
\end{equation*} This has finished the proof.

\subsection{Proof of Lemma \ref{lm: tensor reflecting cloning}}
First notice $\frac{\A + \B}{\sqrt{2}}$ is an orthogonal matrix. If $\bcW_0 \sim N(0,1)^{\otimes (n^{\otimes d})}$ then by the orthogonal invariant property of Gaussian, it is easy to check $\bcW_0^{\sigma^{\otimes d}} \times_i \frac{\A + \B}{\sqrt{2}} \sim N(0,1)^{\otimes (n^{\otimes d})}$ for all $1 \leq i \leq d$. So the property 1 of the Lemma is established.

Now when $\bcW = \bcW_0 \sim \lambda \cdot \bu_1 \circ \cdots \circ \bu_d + \bcZ$ where $\bcZ \sim N(0,1)^{\otimes (n^{\otimes d})}$, we first consider its update after one step. 
\begin{equation}\label{eq: decomp of W}
\begin{split}
    \bcW' = & \bcW^{\sigma^{\otimes d}} \times_1 \frac{\A + \B}{\sqrt{2}} \times \cdots \times_d \frac{\A + \B}{\sqrt{2}}\\
    = & \lambda \times_1  \frac{\A + \B}{\sqrt{2}} \bu_1^\sigma \times \cdots \times_d  \frac{\A + \B}{\sqrt{2}} \bu_d^\sigma + \bcZ^{\sigma^{\otimes d}} \times_1 \frac{\A + \B}{\sqrt{2}} \times \cdots \times_d \frac{\A + \B}{\sqrt{2}}\\
    = & \frac{\lambda}{\sqrt{2}^d} \times_1  (\A + \B) \bu_1^\sigma \times \cdots \times_d  (\A + \B)  \bu_d^\sigma + \bcZ^{\sigma^{\otimes d}} \times_1 \frac{\A + \B}{\sqrt{2}} \times \cdots \times_d \frac{\A + \B}{\sqrt{2}}.
\end{split}
\end{equation} By the result of the first part, we know $\bcZ^{\sigma^{\otimes d}} \times_1 \frac{\A + \B}{\sqrt{2}} \times \cdots \times_d \frac{\A + \B}{\sqrt{2}} \sim N(0,1)^{\otimes (n^{\otimes d})}.$ Now we consider the first term at the right hand side of \eqref{eq: decomp of W}. Denote $\bu_i^{(0)} = \bu_i$ and $\bu_i^{(\ell)} = (\A + \B) (\bu_i^{(\ell-1)})^{\sigma}, \, (1 \leq i \leq d) $. Since $\frac{\A + \B}{\sqrt{2}}$ is orthogonal,
\begin{equation*}
    \|\bu_i^{(\ell)} \|_2^2 = 2\| (\bu_i^{(\ell-1)})^{\sigma} \|_2^2 = 2\|\bu_i^{(\ell-1)}\|_2^2.
\end{equation*} After $\ell$ iterations, we have $\|\bu_i^{(\ell)}\|_2^2 = 2^\ell \| \bu_i \|_2^2 $ for $1 \leq i \leq d$. Also $\bu_i^{(\ell)} \in \bbZ^n$ since each entry of $(\A + \B)$ belongs to $\{ -1, 0, 1 \}$. 
By the definition of $\bu_i^{(\ell)}$, after $\ell$ steps, we have
\begin{equation*}
    \varphi\left(  \lambda \cdot \bu_1 \circ \cdots \circ \bu_d + N(0,1)^{\otimes (n^{\otimes d})} \right) \sim  \cL \left( \frac{\lambda}{\sqrt{2}^{d\ell}}\cdot \bu_1^{(\ell)} \circ \cdots \circ \bu_d^{(\ell)} + N(0,1)^{\otimes (n^{\otimes d})} \right).
\end{equation*} The first part of statement (2) follows.

At the same time if $\bu_i = \bu_j$ for $i, j \in [d], i \neq j$, then $\bu_i^{(\ell)}=\bu_j^{(\ell)}$ holds almost surely by the definition of $\bu_i^{(\ell)}$. 

To finish the proof of the statement (2), we only need to show the desired bound for $\|\bu_i^{(\ell)}\|_0, (1 \leq i \leq d)$. Since $\frac{\A + \B}{\sqrt{2}}$ only has two non-zero values for each row and each column, $\|\bu_i^{(\ell)}\|_0 \leq 2 \|\bu_{i}^{(\ell-1)}\|_0$. Iterate this step, we have $\|\bu_i^{(\ell)}\|_0 \leq 2^\ell \|\bu_{i}\|_0$. The proof for the lower bound of $\|\bu_i^{(\ell)}\|_0$ follows the same as the one of Lemma 26 in \cite{brennan2018reducibility} and we omit it here for simplicity. This has finished the proof.

\section{Proofs of Computational Lower Bounds of CHC$_D$, CHC$_R$}

\subsection{Proof of Theorem \ref{thm: CHC_D comp lower bound}}
The high level idea to show the computational lower bound of $\CHC_D$ is given in Section \ref{sec: ROHC_D and ROHC_R comp lower bound}. The randomized polynomial-time reduction procedure we use shares the same idea as \cite{ma2015computational}, but are modified to handle high-order case. Also we note applying tensor version Gaussian distributional lifting and multivariate rejection kernel techniques in \cite{brennan2018reducibility,brennan2019universality} could probably also yield the same tight computational lower bound.  

First, we introduce a few necessary notation. Let $N = d n \ell$ with $\ell \in \mathbb{N}$ to be chosen depending on $n, k, \lambda$ of CHC model and $\bcA \in (\{ 0, 1 \}^{N})^{\otimes d}$ be the adjacency tensor of hypergraph $G$. 

\begin{algorithm}[h]\caption{Randomized Polynomial-time Reduction for CHC Detection}\label{alg: CHC detection reduction}
\begin{algorithmic}[1]
\State {\bf Input:} $\bcA \in (\{ 0, 1 \}^{N})^{\otimes d}$.
	\State Let $\bcA_0 = \bcA_{[1: n\ell, n\ell +1: 2n \ell, \ldots, (d-1)n \ell +1: dn \ell]}$.  Let $\RK_G = \RK(1 \to N(\xi,1), \frac{1}{2} \to N(0,1), T)$ with $T = \lceil 2(d+1)\log_2 (n\ell)\rceil$ and $\xi = \frac{\log 2}{ 2 \sqrt{ 2(d+1)\log (n\ell) + 2 \log 2 }}$. Compute the tensor $\bcB \in \bbR^{{(n\ell)}^{\otimes d}}$ with $\bcB_{[i_1, \ldots, i_d]} = \RK_G(\bcA_{0[i_1, \ldots, i_d]})$.
	\State Construct $\bcY \in \mathbb{R}^{{n}^{\otimes d}}$ as follows:
	\begin{equation}\label{eq: CHC detection reduction summation}
	\bcY_{[i_1, \ldots, i_d]} = \frac{1}{\ell^{ \frac{d}{2} }} \sum_{j_1 \in [\ell]} \cdots \sum_{j_d \in [\ell]} \bcB_{[(j_1-1)n +i_1, \ldots, (j_d-1)n + i_d]}, \quad 1\leq i_1, \ldots, i_d \leq n.
	\end{equation}
	\State {\bf Output:} $\bcY$.
\end{algorithmic}
\end{algorithm}

Algorithm \ref{alg: CHC detection reduction} summaries the randomized polynomial-time reduction procedure: 
\begin{equation}
\begin{split}
\varphi: (\{ 0, 1 \}^{N})^{\otimes d} \to \mathbb{R}^{{n}^{\otimes d}} :
\bcA \to \bcY.
\end{split}
\end{equation}

Lemma \ref{lem: CHC detection HPC reduction} shows that the randomized polynomial-time algorithm we construct in Algorithm \ref{alg: CHC detection reduction} maps $\HPC_D$ to $\CHC_D$ asymptotically. 
\begin{Lemma} \label{lem: CHC detection HPC reduction}
Given hypergraph $G$ and its adjacency tensor $\bcA$. When $N = dn\ell$ for some $\ell > 0$ integer, if $G \sim H_0^G$, 
\begin{eqnarray}
\TV(N(0,1)^{\otimes(n^{\otimes d})}, \mathcal{L}(\varphi(\bcA))) \leq \frac{1}{n}.
\end{eqnarray}

If $G \sim H_1^G$, and in addition assume $\xi = \frac{\log 2}{ 2 \sqrt{ 2(d+1)\log (n\ell) + 2 \log 2 }}$, $\kappa \geq 4d, \kappa = 4dk, n \geq 2\kappa$, then there exists a prior $\pi$ on $\sX_{\bk,\bn, \frac{\xi}{(\frac{N}{dn})^{d/2}}}$ such that
\begin{equation}
\TV(\cL(\varphi(\bcA)), \bbP_\pi) \leq \frac{1}{n} + 2d\, \exp(-\frac{c}{2d^2} \kappa) + dk\, \exp(-\frac{c}{4d}\kappa \log \frac{n}{k}),
\end{equation} where $\bbP_{\pi} (\cdot) = \int_{ \sX_{\bk,\bn, \frac{\xi}{(\frac{N}{dn})^{d/2}}} } \bbP_{\bcX}(\cdot) d\pi(\bcX)$, and $c > 0$ is some fixed constant.
\end{Lemma}
Lemma \ref{lem: CHC detection HPC reduction} specifically implies that if $N = dn\ell$, $\kappa \geq 4d, \kappa = 4dk, n \geq 2\kappa$ and $\mu = \frac{\xi}{(\frac{N}{dn})^{d/2}}$ with $\xi = \frac{\log 2}{ 2 \sqrt{ 2(d+1)\log (n\ell) + 2 \log 2 }}$, the reduction map $\varphi(G)$ we construct from Algorithm \ref{alg: CHC detection reduction} satisfies $$TV(\varphi(\HPC_D(N,1/2,\kappa)), \CHC_D(\bn,\bk,\mu)) \to 0$$ under both $H_0$ and $H_1$.

Now we prove the computational lower bound of CHC detection by contradiction. If when $\beta > (\alpha - \frac{1}{2})d \lor 0=: \beta^c_{\CHC_D}$, there exists a sequence of algorithm $\{\phi\}_{n}$ such that $\liminf \cE_{\CHC_D}(\phi_n) < \frac{1}{2}$. Then under this regime, we can find $\lambda$ such that 
\begin{equation} \label{ineq: lambda condition}
    \lambda \leq \frac{\log 2}{ 2 \sqrt{ 2(d+1)\log (n\ell) + 2 \log 2 }}\quad \text{ and }\quad \lambda \leq \frac{C n^{\frac{d}{2}}}{k^{d + \delta}},
\end{equation} for some $\delta > 0$. Here we let the sparsity and signal strength of the $\CHC_D$ satisfy $k = \tilde{\Theta}(n^\alpha)$, $\lambda = \tilde{\Theta}(n^{-\beta})$.

Under the first condition in \eqref{ineq: lambda condition}, there exists $\ell \geq 1$ such that $ \lambda \leq \frac{\xi}{\ell^{\frac{d}{2}}} = \frac{\xi}{(N/(dn))^{\frac{d}{2}}}$ by the definition of $\xi$. Thus $\sX_{\bk, \bn, \frac{\xi}{ (\frac{N}{dn})^{d/2} }}$ is supported on $\sX_{\bk, \bn, \lambda}$ and let $\ell$ to be the largest integer satisfies $\frac{\xi}{\ell^{\frac{d}{2}}} \geq \lambda$. 

By combining Lemma \ref{lem: CHC detection HPC reduction} and Lemma \ref{lm: average case reduction simultanious mapping}, we have 
\begin{equation*}
    \cE_{\HPC_D}(\phi_n \circ \varphi ) \leq \cE_{\CHC_D}(\phi_n) + \frac{2}{n} + 2d \exp(-\frac{c}{2d^2} \kappa) +d k \, \exp(-k \log \frac{n}{k}  ),
\end{equation*} and $\liminf_{n \to \infty} \cE_{\HPC_D}(\phi_n \circ \varphi) < \frac{1}{2}$, i.e. $\phi_n \circ \varphi$ has asymptotic risk less than $\frac{1}{2}$ in HPC detection.

On the other hand, note $d$ is fixed and $\ell$ is a largest integer satisfying $ \lambda \leq \frac{\xi}{\ell^{\frac{d}{2}}}$, so combining with the second condition in \eqref{ineq: lambda condition}, we have 
\begin{equation*}
     \frac{\ell^{\frac{d}{2}}}{\xi} \geq \frac{k^{d+ \delta}}{C n^{\frac{d}{2}}} \Longrightarrow k^{d + \delta} \leq C N^{\frac{d}{2}} \sqrt{ 2(d+1) \log N } \Longrightarrow \left(\frac{\kappa}{4d}\right)^{d + \delta} \leq  C N^{\frac{d}{2}} \sqrt{ 2(d+1) \log N },
\end{equation*} 
which implies $\liminf_{n \to \infty} \frac{\log \kappa}{\log N} \leq \frac{1}{2 + \frac{\delta}{2}} < \frac{1}{2}$. The above two facts together contradicts with Conjecture \ref{conj: hardness of tensor clique detection}. So by contradiction argument, we have finished the proof.

\subsection{Proof of Lemma \ref{lem: CHC detection HPC reduction}}
First, by Lemma \ref{lm: Rejection kernel for CHC recovery}, when $\xi = \frac{\log 2}{ 2 \sqrt{2(d+1)\log (n\ell) + 2 \log 2  } }$, we have 
\begin{equation} \label{eq: rejection kernel in CHC_D}
    \TV\left( \RK_G(\Bern(1)), N(\xi,1) \right) = O({(n\ell)}^{-(d+1)})\quad  \text{and }\TV(\RK_G(  \Bern(\frac{1}{2})), N(0,1)  ) = O({(n\ell)}^{-(d+1)}).
\end{equation}

	Under $H_0^G$, let $\widetilde{\bcB} \in \bbR^{{(n\ell)}^{\otimes d}}$ be a random tensor with i.i.d. $N(0,1)$ entries independent of $\bcA$ and $\widetilde{\bcY}$ be the quantity by applying \eqref{eq: CHC detection reduction summation} to $\widetilde{\bcB}$. It is straightforward to verify $\widetilde{\bcY}$ has i.i.d. $N(0,1)$ entries, i.e., $\cL(\widetilde{\bcY}) = N(0,1)^{\otimes (n^{\otimes d})}$. So we have
	\begin{equation}
		\TV(N(0,1)^{\otimes (n^{\otimes d})}, \mathcal{L}(\varphi(\bcA))) \leq \TV(\widetilde{\bcB}, \bcB) \leq (n\ell)^d (n \ell)^{-(d+1)} = \frac{1}{n},
	\end{equation}
	where the first inequality is due to Lemma \ref{lm: data processing} and the second inequality is due to the tensorization and \eqref{eq: rejection kernel in CHC_D}.
	
	Under $H_1^G$, we denote the index set of the planted clique to be $V$. Let $V_1 = V \bigcap [n\ell], V_2 =\left( V \bigcap [n\ell +1: 2 n \ell] \right) - n\ell, \ldots, V_d = \left(V \bigcap [(d-1) n\ell +1, dn \ell]\right) - (d-1) n\ell$, here the notation ``a set minus a value" means shift the indices of the set by that value. By assumption, we have $\sum_{i=1}^d |V_i| = \kappa$. Recall $\bcA_0 = \bcA_{[1: n\ell, n\ell +1: 2n \ell, \ldots, (d-1)n \ell +1: dn \ell]}$, so $\bcA_{0_{[i_1, \ldots, i_d]}} = 1$ if $[i_1,\ldots, i_d] \in V_1 \times \cdots \times V_d$ and $\bcA_{0_{[i_1, \ldots, i_d]}} \sim \Bern(1/2)$ otherwise.
	
	Denote the map $h(x): = 1 + (x-1) \text{ mod } n$. Let 
	\begin{equation}
		U_i = h(V_i), \, i = 1, 2, \ldots, d,
	\end{equation} where $U_i$ could be viewed as the latent signal support on $\bcY$.
	
	We define sets 
	\begin{equation}
		\begin{split}
			& N_{i_1, \ldots, i_d} : = [h^{-1}(i_1), \ldots, h^{-1}(i_d)] \setminus (V_1 \times \cdots \times V_d)\\
			& T_{i_1, \ldots, i_d} := [h^{-1}(i_1), \ldots, h^{-1} (i_d)]\bigcap (V_1 \times \cdots \times V_d)
		\end{split}
	\end{equation} for $1 \leq i_j \leq n \,( 1 \leq j \leq d)$. By construction, we have $\bcY_{[i_1, \ldots, i_d]}$ is the normalized summation of values of $\bcB$ on sets $N_{i_1, \ldots, i_d}$ and $T_{i_1, \ldots, i_d}$.

	We divide the rest of the proof into two steps:
	
	\textbf{Step 1.} In this step we show the event 
	\begin{equation*}
		E = \{|U_1| \geq k, |U_2| \geq k, \ldots, |U_d| \geq k  \},
	\end{equation*}happens with high probability.
	
	To show this, we first show $|V_1|, \ldots, |V_d|$ are concentrated around $\frac{\kappa}{2d}$. By symmetry, we only need to consider $|V_1|$. Notice $|V_1|$ follows the hypergeometric distribution $H\cG (N, \kappa, \frac{N}{d})$. By concentration result of hypergeometric distribution in \cite{hush2005concentration}, we have
	\begin{equation} \label{eq: concentration of |V|}
		\bbP(\left||V_1| - \frac{\kappa}{d}\right| \geq \frac{\kappa}{2d}) \leq \exp(- 2 \alpha_h (\frac{\kappa^2}{4d^2} - 1)) \leq \exp\left(-\frac{c}{2d^2} \kappa\right),
	\end{equation}
	for some constant $c > 0$ and $\alpha_h := \max( \frac{1}{N/d + 1} + \frac{1}{\frac{d-1}{d} N +1}, \frac{1}{\kappa+1} + \frac{1}{N - \kappa + 1} )$. The last inequality is due to the fact that when $N \geq dn$ and $n \geq 2 \kappa$, $\alpha_h$ is of order $\frac{1}{\kappa + 1}$.
	
	Since $\frac{\kappa}{d}= 4k$, $\frac{\kappa}{2d} \asymp \frac{3\kappa}{2d} \asymp k$, then for $\kappa_1 \in [\frac{\kappa}{2d}, \frac{3\kappa}{2d} ]$, 
	\begin{equation}
	\begin{split}
		\bbP( |U_1| \leq k | |V_1| = \kappa_1 ) & \leq \sum_{j = \lceil \kappa_1 / \ell \rceil }^{k} \frac{{n \choose j}{j\ell \choose \kappa_1}}{{n\ell \choose \kappa_1}}\\
		& \leq k \frac{{n \choose k} {k\ell \choose \kappa_1}}{{n\ell \choose \kappa_1}} \leq  k (\frac{en}{k})^k  (\frac{e k \ell}{\kappa_1})^{\kappa_1} (\frac{\kappa_1}{n\ell - \kappa_1})^{\kappa_1}\\
		& \leq k \exp\left(k \log \frac{en}{k} - \kappa_1 \log (\frac{n\ell - \kappa_1}{ek \ell} )\right)\\
		& \leq k \exp(k \log \frac{en}{k} - \frac{\kappa}{2d} \log(\frac{n\ell}{k \ell} )  )\\
		& \leq k\exp(-\frac{c}{4d}\kappa \log \frac{n}{k}),
	\end{split}
	\end{equation}
	for some $c > 0$ and here the first inequality is due to the fact that $j \longmapsto \frac{{n \choose j}{j\ell \choose \kappa_1}}{{n\ell \choose \kappa_1}}$ is an increasing function of $j$ when $n \geq 2 \kappa$ and the last inequality is because $k = \frac{\kappa}{4d}$.
	So \begin{equation}
		\begin{split}
			\bbP(|U_1| \leq k) & \leq \sum_{\kappa_1 = 0}^{\kappa} \bbP(|U_1| \leq k | |V_1| = \kappa_1 ) \bbP(|V_1| = \kappa_1)\\
			& \leq \bbP(|V_1| \leq \frac{\kappa}{2d} ) + \bbP(|V_1| \geq \frac{3}{2d} \kappa ) + \max_{\kappa_1 \in [ \frac{\kappa}{2d}, \frac{3\kappa}{2d} ]} \bbP (|U_1| \leq k | |V_1| = \kappa  )\\
			& \leq 2 \exp(-\frac{c}{2d^2} \kappa) + k\exp(-\frac{c}{4d}\kappa \log \frac{n}{k}).
		\end{split}
	\end{equation}
	Thus
	\begin{equation}
		\bbP(E^c) = \bbP\left(\left(\bigcap_i \{ U_i \geq k \} \right)^c\right) \leq 2d\, \exp(-\frac{c}{2d^2} \kappa) + dk\, \exp(-\frac{c}{4d}\kappa \log \frac{n}{k}),
	\end{equation} for some $c > 0$.
	
	\textbf{Step 2.} Now condition on $V$, generate $\widecheck{\bcB} \in \bbR^{{(n\ell)}^{\otimes d}}$ random tensor such that $\widecheck{\bcB}_{[V_1, \ldots, V_d]} \overset{i.i.d.}\sim N(\xi, 1)$ and the rest of the entries of $\widecheck{\bcB}$ are i.i.d. $N(0, 1)$. Denote $\widecheck{\bcY}$ as the quantity of applying $\widecheck{\bcB}$ to \eqref{eq: CHC detection reduction summation} of Algorithm \ref{alg: CHC detection reduction}. $\widecheck{\bcB}$ and $\widecheck{\bcY}$ are the ideal values we want $\bcB$ and $\bcY$ to be under $H_1^G$.
	
	By the construction of $\widecheck{\bcY}$, for any $(i_1, \ldots, i_d)$, we have $\bbE[ \widecheck{\bcY}_{[i_1, \ldots, i_d]} ] = \frac{\xi |T_{i_1, \ldots, i_d}|}{\ell^{\frac{d}{2}}}$. Since for any $[i_1, \ldots, i_d] \in V_1 \times V_2 \times \cdots \times V_d$, $|T_{i_1, \ldots, i_d}| \geq 1$, in this case $\bbE[ \widecheck{\bcY}_{[i_1, \ldots, i_d]} ] \geq \frac{\xi}{\ell^{\frac{d}{2}}}$.
	
	Also by the construction of $\widecheck{\bcY}$, we have the entries of $\widecheck{\bcY}$ are independent. Since $1_E$ is deterministic given $V$, for any $V$ such that $1_E = 1$, there exists some $\bcX = \bcX(V) \in \sX(\bk,\bn, \frac{\xi}{\ell^{\frac{d}{2}}})$ such that $\cL(\widecheck{\bcY}| V) = \bbP_{\bcX}$.
	
	Define the probability distribution $\pi = \cL(\bcX(V) | E)$ which is supported on the set $\sX(\bk,\bn, \frac{\xi}{\ell^{\frac{d}{2}}})$. Then $\cL(\widecheck{\bcY} | E) = \bbP_{\pi}$ is a mixture of distributions of $\{ \bbP_{\bcX}: \bcX \in \sX(\bk,\bn, \frac{\xi}{\ell^{\frac{d}{2}}})  \}$.
	
	Now we are ready to show that $\TV(\cL(\bcY), \bbP_\pi)$ is small.
	\begin{equation}
	\begin{split}
			\TV( \cL(\bcY), \bbP_\pi ) & \overset{(a)}\leq \TV ( \cL(\bcY), \cL(\widecheck{\bcY}) ) + \TV ( \cL(\widecheck{\bcY}), \bbP_\pi )\\
		& \leq \bbE_V [\TV ( \cL(\bcY|V), \cL(\widecheck{\bcY}|V) ) ] + \TV( \cL(\widecheck{\bcY}), \cL(\widecheck{\bcY} | E) )\\
		& \overset{(b)}\leq \sum_{i_1, \ldots, i_d = 1}^{n \ell} \TV( \cL(\bcB_{[i_1, \ldots, i_d]}|V), \cL(\widecheck{\bcB}_{[i_1, \ldots, i_d]}|V) ) + \bbP(E^c)\\
		& \overset{(c)}\leq (n\ell)^{d} (n\ell)^{-(d+1)}  + \bbP(E^c)\\
		& \leq \frac{1}{n} + \bbP(E^c),
	\end{split}
	\end{equation}	
where (a) is due to triangle inequality, (b) is due to Lemma \ref{lm: data processing} and Lemma 7 of \cite{brennan2018reducibility} and (c) is due to \eqref{eq: rejection kernel in CHC_D}. This has finished the proof of this lemma.

\subsection{Proof of Theorem \ref{thm: CHC_R computation lower bound}}
In this section, we prove the computational lower bound for $\CHC_R(\bk, \bn, \lambda)$ and it is enough to prove the lower bound for the case: $$\bcX \in \left\{\lambda \mathbf{1}_{I_1} \circ \cdots \circ \mathbf{1}_{I_d}: I_i \subseteq \{1,\ldots, n_i\}, |I_i| = k_i \right\}.$$

We first introduce a randomized polynomial-time Algorithm \ref{alg: CHC recovery reduction} to do reduction from $\HPDS(n,\kappa, \frac{1}{2} + \rho, \frac{1}{2})$ to $\CHC(\bn, \bk, \lambda)$.

\begin{algorithm}[h]\caption{Randomized Polynomial-time Reduction for CHC Recovery} \label{alg: CHC recovery reduction}
\begin{algorithmic}[1]
\State {\bf Input:} Hypergraph $G$ and the density bias $\rho$.
	\State Let $\RK_G = \RK( \frac{1}{2} + \rho \to N(\xi, 1), \frac{1}{2} \to N(0,1), T )$ where $\xi = \frac{\log (1 + 2 \rho) }{2\sqrt{2(d+1) \log n + 2 \log 2 }}$ and $T = \lceil 2(d+1)  \log_{1+ 2\rho} n \rceil$ and compute the symmetric tensor $\bcW \in \bbR^{n^{\otimes d}}$ with $\bcW_{[i_1, \ldots, i_d]} = \RK_G(\mathbf{1}((i_1, \ldots, i_d) \in E(G) ))$. Let the diagonal entries of $\bcW_{[i, \ldots, i]}$ to be i.i.d. $N(0,1)$.
	\State Generate $(d!-1)$ i.i.d. symmetric random tensor $\bcB^{(1)}, \ldots, \bcB^{(d!-1)}$ in the following way: their diagonal values are $0$ and non-diagonal values are i.i.d. $N(0,1)$. Given any non-diagonal index $\bi = (i_1, \ldots, i_d)$ ($i_1 \leq i_2 \leq \ldots \leq i_d$), suppose it has $D\, (D \leq d!)$ unique permutations and denote them as $\bi_{(0)}: = \bi, \bi_{(1)}, \ldots, \bi_{(D-1)}$, then we transform $\bcW$ in the following way
	\begin{equation*}
	    \left( \begin{array}{c}
	         \bcW_{\bi_{(0)}}  \\
	         \bcW_{\bi_{(1)}}\\
	         \vdots \\
	         \bcW_{\bi_{(D-1)}}
	    \end{array}  \right) = \left( \frac{\mathbf{1}}{\sqrt{d!}}, \left(\frac{\mathbf{1}}{\sqrt{d!}} \right)_{\perp}  \right)_{[1:D,:]} \times \left(  
	    \begin{array}{c}
	         \bcW_{[i_1, \ldots, i_d]}  \\
	         \bcB^{(1)}_{[i_1, \ldots, i_d]} \\
	         \vdots\\
	         \bcB^{(d!-1)}_{[i_1, \ldots, i_d]} 
	    \end{array} \right).
	\end{equation*} 
	Here $\frac{\mathbf{1}}{\sqrt{d!}}$ is a $\bbR^{d!}$ vector with all entries to be $\frac{1}{\sqrt{d!}}$ and $\left(\frac{\mathbf{1}}{\sqrt{d!}} \right)_{\perp} \in \bbR^{d! \times (d! - 1)}$ is any orthogonal complement of $\frac{\mathbf{1}}{\sqrt{d!}}$. 
	\State {\bf Output:} $\bcW^{\id, \sigma_1, \ldots, \sigma_{d-1}}$, where $\sigma_1, \ldots, \sigma_{d-1}$ are independent permutations of $[n]$ chosen uniformly at random.
\end{algorithmic}
\end{algorithm}

The Lemma \ref{lm: CHC recovery reduction} shows randomized polynomial-time mapping in Algorithm \ref{alg: CHC recovery reduction} maps HPDS to CHC asymptotically. 
\begin{Lemma}\label{lm: CHC recovery reduction}
Suppose that $n, \xi$ and $\rho \geq \frac{1}{n^{\frac{d-1}{2}}}$ are such that $$\xi = \frac{\log (1 + 2 \rho)}{2 \sqrt{2(d+1) \log n + 2 \log 2 }},$$ then the randomized polynomial-time computable map $\varphi : \cG_d(n) \to \bbR^{n^{\otimes d}}$ represented by Algorithm \ref{alg: CHC recovery reduction} holds the following: for any subset $S \subseteq [n]$ with $|S| = \kappa$,
\begin{equation*}
\begin{split}
    & \TV\left(  \varphi \left( \cG_d(n,\kappa, \frac{1}{2} + \rho, \frac{1}{2}, S ) \right) , \int \cL \left(  \frac{\xi}{\sqrt{d!}} \cdot \mathbf{1}_{S} \circ \mathbf{1}_{T_1} \circ \cdots \circ \mathbf{1}_{T_{d-1}} + \bcZ \right) d \pi(T_1) \cdots d\pi(T_{d-1}) \right) \\
    &= O\left(  \frac{1}{\sqrt{\log n}}   \right),
\end{split}
\end{equation*} where $\bcZ \sim N(0,1)^{\otimes (n^{\otimes d})}$ and $\pi$ is the uniform distribution on subsets of $[n]$ of size $\kappa$ and $\cG_d(n, \kappa,\frac{1}{2}+ \rho, \frac{1}{2}, S)$ represents the distribution of HPDS $\cG_d(n,\kappa ,\frac{1}{2}+ \rho, \frac{1}{2})$ with planted dense subgraph supported on set $S$.
\end{Lemma}

Lemma \ref{lm: CHC recovery reduction} specifically implies that if $k = \kappa$, $\lambda = \frac{\log (1 + 2 \rho)}{2\sqrt{d!} \sqrt{2(d+1) \log n + 2 \log 2}  }$, the reduction map $\varphi(G)$ in Algorithm \ref{alg: CHC recovery reduction} satisfies
\begin{equation*}
    \TV(\varphi(\HPDS_R(n,\kappa,1/2 +\rho, 1/2)), \cL(\CHC_R(\bn,\bk,\lambda))) \to 0.
\end{equation*}

Next, we prove Theorem \ref{thm: CHC_R computation lower bound} by contradiction argument and we divide the proof into dense and sparse regimes.

\begin{itemize}[leftmargin=*]
    \item \textbf{$\alpha \geq \frac{1}{2}$.} Let $k = \kappa \lceil n^{\alpha} \rceil$ and $\rho = n^{-\beta}$. Then it is easy to check in CHC model, the signal strength and sparsity levels are $\lim_{n \to \infty} \frac{\log (\lambda^{-1})}{\log n} = \beta$ and $\lim_{n \to \infty} \frac{\log k}{\log n} = \alpha$.
    
    Suppose when $(d-1)\alpha - \frac{d-1}{2} < \beta$, there is a sequence of polynomial-time algorithm $\{\phi\}_n$ such that w.p. more than $\frac{1}{2}$, it can identify the true planted latent cluster in $\bcY \sim \ROHC_R(\bn,\bk,\lambda)$. Denote the support in $\bcY$ as $S$. Let $\phi^1$ be the restriction of $\phi$ that only output the estimated support of $\bcY$ at mode 1.

Denote $\cL_{S}$ as the distribution of $$\int \cL \left(  \lambda \cdot \mathbf{1}_{S} \circ \mathbf{1}_{T_1} \circ \cdots \circ \mathbf{1}_{T_{d-1}} + N(0,1)^{\otimes (n^{\otimes d})} \right) d \pi(T_1) \cdots d\pi(T_{d-1}),$$ where $\pi$ is the uniform distribution on the subset of $[n]$ of size $\kappa$. By Lemma \ref{lm: CHC recovery reduction} and noticing $\lambda = \xi/\sqrt{d!}$, we have
\begin{equation*}
    \left| \bbP_{\bcW \sim \cL(\varphi(G))} \left( \phi_n^1(\bcW) = S  \right)  - \bbP_{\bcW \sim \cL_S} ( \phi_n^1(\bcW) = S ) \right| \leq \TV(\cL(\varphi(G)), \cL_S) = O(\frac{1}{\sqrt{\log n}}),
\end{equation*} where the inequality is due to the definition of total variation distance.

Note by assumption 
\begin{equation*}
    \bbP_{\bcW \sim \cL_S} ( \phi_n^1(\bcW) = S ) = \bbE_{T_1, \ldots, T_{d-1} \sim \pi} \bbP_{\bcW \sim \cL_{S,T_1, \ldots, T_{d-1}}} ( \phi_n^1(\bcW) = S ) > \frac{1}{2},
\end{equation*} so
\begin{equation*}
    \bbP \left(  \phi^1 \circ \varphi (G) = S  \right) \geq  \bbP_{\bcW \sim \cL_S} ( \phi_n^1(\bcW) = S ) - O(\frac{1}{\sqrt{\log n}}),
\end{equation*} and $ \liminf \cE_{\HPDS_R}(\phi^1\circ \varphi) < \frac{1}{2}$, i.e. $\phi^1 \circ \varphi$ can recovery the support of HPDS with asymptotic risk less than $\frac{1}{2}$.

On the other hand, the condition \eqref{eq: HPDS conjecture} at here becomes
\begin{equation*}
    \lim_{n \to \infty} \log_n \frac{\kappa^{d-1} \rho }{\sqrt{\frac{1}{4} - \rho^2}} = (d-1)\alpha - \beta < \frac{d-1}{2},
\end{equation*} 
where the inequality is because $(d-1)\alpha - \frac{d-1}{2} < \beta$. Combining $1/2 + \rho < 1-\Omega(1)$, these two facts together contradict the HPDS recovery conjecture \ref{conj: HPDS recovery conjecture}. This has finished the proof for the $\alpha \geq \frac{1}{2}$ region.

\item $0<\alpha < \frac{1}{2}$. The main idea to prove the computational lower bound for $\CHC_R$ in this regime is to use the established computational lower bound for $\CHC_D$. Suppose $\beta > 0$ (for example, we can let $\beta = \epsilon$ for some small enough $\epsilon > 0$), there is a sequence of polynomial-time algorithm $\{\phi_R \}_n$ such that $ \liminf \cE_{\CHC_R}(\phi_R) < \frac{1}{2}$. In this regime, we have $\lambda\sqrt{\prodk} \geq C  k^{\frac{d}{4}}$, then by the following Lemma \ref{lm: CHC recovery imply detection}, there exists a sequence of polynomial-time test $\{ \phi_D \}$ with $\lim_{n \to \infty}\cE_{\CHC_D}(\phi_D) < \frac{1}{2}$. This contradicts with the computational lower bound of $\CHC_D$ under the HPC conjecture \ref{conj: hardness of tensor clique detection}. So this has established the computational lower bound for $\CHC_R$ in regime $0<\alpha <1/2$.

In summary, combining the results in dense and sparse regime, we have established the computational lower bound for $\CHC_R$ in the regime $\beta > \beta^c_{\CHC_R} := (d-1)(\alpha -1/2) \vee 0 $.

\begin{Lemma}\label{lm: CHC recovery imply detection}
Consider $\CHC_D(\bk, \bn, \lambda)$ and $\CHC_R(\bk, \bn,\lambda)$ under the asymptotic regime \eqref{assum: asymptotic assumption}. If $\lambda \sqrt{\prodk} \geq Ck^{\frac{d}{4}}$ for some $C > 0$ and there exists a sequence of polynomial-time recovery algorithm $\{ \phi_R \}_n$ such that $\liminf_{n \to \infty} \cE_{\CHC_R}(\phi_R) < \eta$ for some $\eta \in (0,1]$, then there exists a sequence of polynomial-time algorithm $\{\phi_D \}_n$ such that $\liminf_{n \to \infty} \cE_{\CHC_D}(\phi_D) < \eta$. 
\end{Lemma}

\end{itemize}

\subsection{Proof of Lemma \ref{lm: CHC recovery reduction}}
Let $\varphi$ be the map of Algorithm \ref{alg: CHC recovery reduction} and $\bcW_2$, $\bcW_3$, $\bcW_4$ be the value of $\bcW$ after step 2, step 3 and step 4 of $\varphi$.

First by Lemma \ref{lm: Rejection kernel for CHC recovery}, we know that with value $\xi = \frac{\log (1 + 2 \rho)}{ 2 \sqrt{2(d+1)\log n + 2 \log 2  } }$ and $\rho \geq \frac{1}{n^{\frac{d-1}{2}}}$, $\TV\left( \RK_G(\Bern(\frac{1}{2} + \rho)), N(\xi,1) \right) = O(n^{-(d+1)})$  and $ \TV(\RK_G(  \Bern(\frac{1}{2}), N(0,1)  )) = O(n^{-(d+1)}).$ 

Denote $\cM_n : = \cM_n (S, P, Q)$ as the distribution of a $\bbR^{n^{\otimes d}}$ symmetric tensor $\bcA$, where its diagonal entries are independent $N(0,1)$ random variables, and for non-diagonal entries $\bcA_{[i_1,\ldots, i_d]}$, if set $(i_1, \ldots, i_d) \subseteq S$, $\bcA_{[i_1,\ldots, i_d]}$ is draw independently from distribution $P$ and if $(i_1, \ldots, i_d) \nsubseteq S$, $\bcA_{[i_1,\ldots, i_d]}$ is draw independently from $Q$. Here $S$ is a subset of $[n]$ of size $\kappa$.

Let $\bcM_2 \sim \cM_n(S, N(\xi, 1), N(0,1))$ and $\bcM_3$, $\bcM_4$ be the value of $\bcM_2$ after applying step 3 and step 3 and 4 of $\varphi$ to $\bcM_2$. Here $\bcM_2, \bcM_3, \bcM_4$ could be viewed as the ideal values we want $\bcW_2, \bcW_3, \bcW_4$ to be. 

After step 2 of $\varphi$, 
\begin{equation*}
    \begin{split}
        \TV\left( \cL(\bcW_2), \cL(\bcM_2) \right)  \leq & \kappa^d \TV\left(  \RK_G(\Bern(\frac{1}{2} + \rho)), N(\xi,1) \right) \\
        & + \left( n^d - \kappa^d \right) \TV\left(\RK_G(  \Bern(\frac{1}{2}), N(0,1)  )\right) \\
        = & O(\frac{1}{n}),
    \end{split}
\end{equation*} where $G \sim \cG_d(n,\kappa,  \frac{1}{2}+\rho,\frac{1}{2},S)$.

Also by data processing inequality, we have 
\begin{equation} \label{ineq: TV distance after step 3}
    \begin{split}
        \TV\left(  \cL(\bcM_3) , \cL(\bcW_3) \right)  &\leq \TV \left(  \cL(\bcM_2), \cL (\bcW_2)  \right) = O(n^{-1}).
    \end{split}
\end{equation}
Also notice that after applying step 3 of $\varphi$, the diagonal values of $\bcM_3$ are i.i.d. $N(0,1)$ and the non-diagonal values of $\bcM_3$ have the same distribution as the non-diagonal entries of $\frac{\xi}{\sqrt{d!}} \mathbf{1}_{S} \circ \cdots \circ \mathbf{1}_{S}  + N(0,1)^{\otimes (n^{\otimes d})}$. 

Now consider the distribution of $\bcM_4$ condition on permutations $\sigma_1, \ldots, \sigma_{d-1}$. It's entries have the ditribution as entries of $$\frac{\xi}{\sqrt{d!}} \mathbf{1}_S \circ \mathbf{1}_{T_1} \circ \cdots \circ \mathbf{1}_{T_{d-1}} + N(0,1)^{\otimes (n^{\otimes d})},$$ other than entries at indices $(i, \sigma_1(i), \ldots, \sigma_{d-1} (i))$ for $i \in S $. Here $T_1 = \sigma_1(S), \ldots, T_{d-1} = \sigma_{d-1} (S)$. Specifically, conditioning on $\sigma_1(S) = T_1, \ldots, \sigma_{d-1}(S) = T_{d-1}$, we have 
\begin{equation} \label{ineq: conditional TV distance}
    \begin{split}
    &\TV \left( \left( \cL(\bcM_4) | \sigma_1(S) = T_1, \ldots, \sigma_{d-1}(S) = T_{d-1} \right), \cL \left(  \frac{\xi}{\sqrt{d!}} \mathbf{1}_S \circ \mathbf{1}_{T_1} \circ \cdots \circ \mathbf{1}_{T_{d-1}} + \bcZ \right)  \right) \\
    & = \TV  \left(   \left(\cL(\bcM_4)[S \times T_1 \times \cdots \times T_{d-1}] | \sigma_1(S) = T_1, \ldots, \sigma_{d-1}(S) = T_{d-1} \right),\cL \left( N(\frac{\xi}{\sqrt{d!}},1)^{\otimes k^{\otimes d}} \right) \right)\\
    & \overset{(a)}\leq \sqrt{ \chi^2 \left( N(0, 1), N(\frac{\xi}{\sqrt{d!}})  \right)  } \\
      & = O( \frac{1}{\sqrt{\log n}} ).
    \end{split}
\end{equation}
Here $(a)$ is due to Lemma \ref{lm: diagonal permutation}.

Finally, we have 
\begin{equation*}
    \begin{split}
      & \TV\left( \varphi \left( \cG_d(n,\kappa, \frac{1}{2} + \rho, \frac{1}{2}, S ) \right), \int \cL \left(  \frac{\xi}{\sqrt{d!}} \cdot \mathbf{1}_{S} \circ \mathbf{1}_{T_1} \circ \cdots \circ \mathbf{1}_{T_{d-1}} + \bcZ \right) d \pi(T_1) \cdots d\pi(T_{d-1}) \right)\\
        \overset{(a)}\leq &\TV \left(  \varphi \left( \cG_d(n,\kappa, \frac{1}{2} + \rho, \frac{1}{2}, S ) \right), \cL(\bcM_4) \right) \\
        & + \TV \left( \cL(\bcM_4),  \int \cL \left(  \frac{\xi}{\sqrt{d!}} \cdot \mathbf{1}_{S} \circ \mathbf{1}_{T_1} \circ \cdots \circ \mathbf{1}_{T_{d-1}} + \bcZ \right) d \pi(T_1) \cdots d\pi(T_{d-1})  \right)\\
      \overset{(b)} \leq  & \TV \left( \cL(\bcW_3), \cL(\bcM_3) \right)\\
        & + \TV \left( \cL(\bcM_4),  \int \cL \left(  \frac{\xi}{\sqrt{d!}} \cdot \mathbf{1}_{S} \circ \mathbf{1}_{T_1} \circ \cdots \circ \mathbf{1}_{T_{d-1}} + \bcZ \right) d \pi(T_1) \cdots d\pi(T_{d-1})  \right)\\
      \overset{(c)}\leq  & O(\frac{1}{n}) + O( \frac{1}{\sqrt{\log n}} ) = O( \frac{1}{\sqrt{\log n}} ),
    \end{split}
\end{equation*} where $(a)$ is due to triangle inequality, (b) is due to data processing inequality and (c) is due to \eqref{ineq: TV distance after step 3}\eqref{ineq: conditional TV distance}. So we finish the proof of this Lemma.

\subsection{Proof of Lemma \ref{lm: CHC recovery imply detection}}
The proof idea is similar to the proof of statement (2) of Lemma \ref{lm: ROHS recovery imply detection}. Given $\bcY$ generated from $$\cL \left( \lambda \cdot \mathbf{1}_{I_1} \circ \cdots \circ \mathbf{1}_{I_d} + N(0,1)^{\otimes n_1 \times \cdots \times n_d} \right),$$ by the property of Gaussian, it is easy to check that $\bcA := \frac{\bcY + \bcZ_1}{\sqrt{2}}$ and $\bcB:=\frac{\bcY - \bcZ_1}{\sqrt{2}}$ are two independent copies with distribution $$\cL \left( \frac{\lambda}{\sqrt{2}} \cdot \mathbf{1}_{I_1} \circ \cdots \circ \mathbf{1}_{I_d} + N(0,1)^{\otimes n_1 \times \cdots \times n_d} \right),$$ if $\bcZ_1 \sim N(0,1)^{\otimes n_1 \times \cdots \times n_d}$ and independent of $\bcY$.

The rest of the proof is the same as the proof of statement (2) of Lemma \ref{lm: ROHS recovery imply detection} by replacing $\mu$ with $\lambda \sqrt{\prodk}$.

\section{Proofs of Computational Lower Bounds of $\ROHC_D$ and $\ROHC_R$}\label{sec:proof_lower bound}

\subsection{Proof of Lemma \ref{lm: ROHS reduction guarantee}}
Let $\varphi: \cG_d(n) \to \bbR^{n^{\otimes d}}$ be the map in Algorithm \ref{alg: ROHC reduction} and let $\bcW_4$, $\bcW_5$ be the values of $\bcW$ after step 4 and step 5 of Algorithm \ref{alg: ROHC reduction}. Also by Lemma \ref{lm: Rejection kernel for CHC recovery}, we know that with value $\xi = \frac{\log 2}{ 2 \sqrt{2(d+1)\log n + 2 \log 2  } }$, we have 
\begin{equation}\label{eq: ROHC rejection kernel}
    \TV\left( \RK_G(\Bern(1)), N(\xi,1) \right) = O(n^{-(d+1)}), \, \text{and }\TV(\RK_G(  \Bern(\frac{1}{2}), N(0,1)  )) = O(n^{-(d+1)}).
\end{equation}

If $G \sim \cG_d(n, \frac{1}{2})$, then
\begin{equation*}
    \TV \left( \varphi(\cG_d(n, \frac{1}{2})), N(0,1)^{\otimes (n^{\otimes d})} \right) \leq \TV\left(\varphi_1(\cG_d(n, \frac{1}{2})), N(0,1)^{\otimes (n^{\otimes d})} \right)\leq  O(\frac{1}{n}),
\end{equation*} where $\varphi_1$ denotes the step 1 of $\varphi$, the first inequality is due to data processing inequality and the second one is due to tensorization and \eqref{eq: ROHC rejection kernel}.

Now we consider the case $G \sim \cG_d(n, \frac{1}{2}, \kappa)$. First following the same proof of Lemma \ref{lm: CHC recovery reduction}, we have
\begin{equation} \label{ineq: ROHC TV distance 1}
\begin{split}
     &\TV\left(  \cL(\bcW_4) , \int \cL \left(  \frac{\xi}{\sqrt{d!}} \cdot \mathbf{1}_{T_1} \circ \cdots \circ \mathbf{1}_{T_{d}} + N(0,1)^{\otimes (n^{\otimes d})} \right) d \pi'(T_1, \ldots, T_d) \right)\\
     =& O\left(  \frac{1}{\sqrt{\log n}}   \right),
\end{split}
\end{equation} where $\pi'$ is the uniform distributions over pairs $(T_1, \ldots, T_d)$ of $\kappa$-subsets $T_1, \ldots, T_d \subseteq [n]$. 

Let $\bcM_4$ be a tensor distributed as 
\begin{equation*}
    \frac{\xi}{\sqrt{d!}} \cdot \mathbf{1}_{T_1} \circ \cdots \circ \mathbf{1}_{T_{d}} + N(0,1)^{\otimes (n^{\otimes d})},
\end{equation*} where $T_1, \ldots, T_d$ are $\kappa$-subsets of $[n]$ chosen uniformly at random. Also let $\bcM_5$ be the value of $\bcM_4$ after applying step 5 of $\varphi$ to $\bcM_4$. Again $\bcM_4$ and $\bcM_5$ are the ideal values we want $\bcW_4$ and $\bcW_5$ to be. By statement 2 of Lemma \ref{lm: tensor reflecting cloning}, the distribution of $\bcM_5$ condition on sets $T_1, \ldots, T_d$ is given by
\begin{equation*}
    \cL(\bcM_5|T_1, \ldots, T_d) \sim \int \cL \left( \frac{\xi }{\sqrt{d!} (\sqrt{2})^{d\ell} } \bv_1 \circ \cdots \circ \bv_d + N(0,1)^{\otimes (n^{\otimes d})} \right) d \bar{\pi}(\bv_1, \ldots, \bv_d),
\end{equation*} where $\bar{\pi} := \bar{\pi}_{T_1, \ldots, T_d}$ is a prior defined in Lemma \ref{lm: tensor reflecting cloning}. 

As shown in Lemma \ref{lm: tensor reflecting cloning}, $\bv_i$s $(1 \leq i \leq d)$ are supported on $\bar{\pi}$ and satisfy $\|\bv_i\|_2^2 = 2^\ell \|\mathbf{1}_{T_i} \|_2^2 = 2^\ell \kappa$,$\|\bv_i\|_0 \leq 2^\ell \kappa$. Moreover, with probability at least $1 - 4d/\kappa$ on $\bar{\pi}$, we have $\|\bv_i\|_0 \geq 2^{\ell} \kappa (1 - \max( \frac{2C \ell \log(2^{\ell} \kappa)}{\kappa}, \frac{2^{\ell}\kappa }{n} ) ) $. We denote this conditional distribution of $\bar{\pi}$ on this event as $\bar{\pi}'$. Since $2^{\ell}\kappa \leq n/C$ for $C > 1$,  if $\bu_i = \frac{1}{\sqrt{2^\ell \kappa}} \bv_i$, then $\bar{\pi}'$ induces a prior on pair $(\bu_1, \ldots, \bu_d)$ in $\cV_{n, 2^\ell \kappa}$. This is because $\bv_i \in \bbZ^n$ as shown in Lemma \ref{lm: tensor reflecting cloning}, and the nonzero entries of $\bu_i$ have magnitudes at least $\frac{1}{\sqrt{2^\ell \kappa}}$.

Let $\pi' = \bbE_{T_1, \ldots, T_d} \left( \bar{\pi}_{T_1, \ldots, T_d} \right) $ be a prior formed by marginalizing $T_1, \ldots, T_d$, and $\pi'$ is also supported on $\cV_{n, 2^\ell \kappa}$. So
\begin{equation} \label{ineq: ROHC TV distance 2}
\begin{split}
    \TV(\cL(\bcM_5),\int \cL \left( \frac{\xi}{\sqrt{d!}} \kappa^{\frac{d}{2}} \bu_1 \circ \cdots \circ \bu_d + N(0,1)^{\otimes (n^{\otimes d})} \right) d \pi'(\bu_1, \ldots, \bu_d)) \leq 4d/\kappa.
\end{split}
\end{equation}
Finally by triangle inequality and Lemma \ref{lm: data processing}, we have
\begin{equation*}
    \begin{split}
        & \TV\left(  \cL(\varphi(G)), \int \cL \left( \frac{\xi \kappa^{\frac{d}{2}}}{\sqrt{d!}} \bu_1 \circ \cdots \circ \bu_d + N(0,1)^{\otimes p^{(\otimes d)}} \right) d \pi'(\bu_1, \ldots, \bu_d) \right)\\
        \leq & \TV \left(  \cL( \bcW_5 ), \cL ( \bcM_5 ) \right) \\
        & + \TV\left( \cL ( \bcM_5 ), \int \cL \left( \frac{\xi}{\sqrt{d!}} \kappa^{\frac{d}{2}} \bu_1 \circ \cdots \circ \bu_d \right) d \pi' (\bu_1, \ldots, \bu_d) \right)\\
        \overset{(a)}\leq & O(\frac{1}{\sqrt{\log n}}),
    \end{split}
\end{equation*} 
where (a) is due to \eqref{ineq: ROHC TV distance 1}, \eqref{ineq: ROHC TV distance 2}. This has finished the proof.

\section{Proofs for the Evidence of HPC Conjecture \ref{conj: hardness of tensor clique detection} and HPDS Conjecture \ref{conj: HPDS recovery conjecture} }\label{sec: proof-hardness}

\subsection{Proof of Proposition \ref{prop: spectral method for HPDS} }
Without loss of generality, we can assume $N$ is a multiplier of $d$, otherwise we can replace $N = d \lfloor \frac{N}{d} \rfloor$. Recall the planted dense subgraph index set is $K$, and let $K_i = K \bigcap [(i-1)\frac{N}{d}+1: \frac{i N}{d}]$. By symmetry, we only need to consider recovering $K_1$ and $K_2$.

First by the same argument as \eqref{eq: concentration of |V|}, we can show $|K_1| \asymp |K_2| \asymp \frac{\kappa}{d}$. Denote $\{X_i^{(1)}\}$ as i.i.d. $\Bern(q_2)$ random variables and $\{X_i^{(2)} \}$ as i.i.d. $\Bern(q_1)$ random variables. Following the same notation in Algorithm \ref{alg: CHC recovery aggregated SVD}, for $(k_1, k_2) \notin K_1 \times K_2$, 
\begin{equation*}
    \Y^{(1,2)}_{[k_1, k_2]} = \frac{ \sum_{i=1}^{ (\frac{N}{d})^{d-2} }  \left(X_i^{(1)} - q_2\right)}{\sqrt{ \left(\frac{N}{d} \right)^{d-2} \left( q_2(1-q_2) \right) }},
\end{equation*} and for $(k_1, k_2) \in K_1 \times K_2$,
\begin{equation*}
    \Y^{(1,2)}_{[k_1, k_2]} =  \frac{\sum_{i=1}^{ (\frac{N}{d})^{d-2} - C(\frac{\kappa}{d})^{d-2} }\left(X_i^{(1)} - q_2\right)}{\sqrt{ \left(\frac{N}{d} \right)^{d-2} \left( q_2(1-q_2) \right) }} +  \frac{ \sum_{i=1}^{  C(\frac{\kappa}{d})^{d-2} } \left(X_i^{(2)} - q_2\right) }{\sqrt{ \left(\frac{N}{d} \right)^{d-2} \left( q_2(1-q_2) \right) }},
\end{equation*} for some constant $C > 0$.

By Chernoff bound, we have if $(k_1, k_2) \notin K_1 \times K_2$
\begin{equation*}
    \bbP\left(| \Y^{(1,2)}_{[k_1, k_2]} | > t \right) \leq \exp(-t^2),
\end{equation*} and if $(k_1, k_2) \in K_1 \times K_2$, 
\begin{equation*}
    \bbP\left(\left| \Y^{(1,2)}_{[k_1, k_2]} - C\frac{ (\frac{\kappa}{d})^{d-2} (q_1 - q_2) }{ \sqrt{ \left(\frac{N}{d} \right)^{d-2} \left( q_2(1-q_2) \right) } } \right| > t \right) \leq \exp(-t^2)
\end{equation*}

So we can write $\Y = \lambda\cdot 1_{K_1} 1_{K_2}^\top + \Z$ where $\lambda = C\frac{ (\frac{\kappa}{d})^{d-2} (q_1 - q_2) }{ \sqrt{ \left(\frac{N}{d} \right)^{d-2} \left( q_2(1-q_2) \right) } }$ and entries $\Z_{ij}$ are independent subgaussian random variable with variance $1$.

To recover $K_1, K_2$, it is the same as biclustering recovery problem studied in literature with parameters $(\frac{N}{d}, |K_1|, |K_2|, \lambda)$. By Lemma 1 of \cite{cai2017computational}, when 
\begin{equation*}
    \lambda \geq C'\frac{\sqrt{\frac{N}{d}}}{|K_1| \land |K_2| }, \text{ i.e., } \limsup_{N \to \infty} \log_N \left(\frac{\kappa^{d-1} (q_1-q_2)}{\sqrt{
			q_2(1-q_2)}}\right) \geq \frac{d}{2} -\frac{1}{2}, 
\end{equation*}
then with probability at least $1 - (\frac{N}{d})^{-c} - \exp(-C \frac{N}{d})$, the output of the Algorithm can exactly recover $K_1$ and $K_2$. Similar analysis holds for other modes.

However when $\limsup_{N \to \infty} \log_N \left(\frac{\kappa^{d-1} (q_1-q_2)}{\sqrt{
 			q_2(1-q_2)}}\right) < \frac{d}{2} -\frac{1}{2}$, then under a variant PC recovery conjecture \cite{cai2017computational}, the biclustering recovery procedure of Algorithm \ref{alg: CHC recovery aggregated SVD} fails with non-trivial probability, so as Algorithm \ref{alg: HPDS recovery aggregated SVD}.

\subsection{Proof of Lemma \ref{lm: HPC recovery hard imply detection hard}}
The proof of this lemma is by a contradiction argument. Suppose for sufficient large $N$, there is no polynomial-time recovery algorithm can output the right clique of $\cG_d(N,1/2,\kappa)$ with success probability at least $1-1/N$, but we can distinguish in polynomial time whether a hypergraph $G$ is drawn from $\cG_d(N, \frac{1}{2})$ or $\cG_d(N, \frac{1}{2},\kappa/3)$ with probability at least $\frac{1}{4N^d}$. Denote $\varphi$ as this distinguisher. 

In Algorithm \ref{alg: HPC recovery based on HPC detection}, we provide an polynomial-time algorithm which can find a clique of size $\kappa$ in $G$ using the distinguisher $\varphi$ for $G$ from $\cG_d(N, \frac{1}{2},\kappa)$. Denote the clique set in HPC as $K$. We first give a high level idea why the algorithm works. When $\bv \subsetneq K$, by construction, after remove $\bv$ and $X$, at most $2\kappa/3$ nodes in the clique is removed with high probability by concentration, thus the graph $G_x$ is a random graph with a planted clique of size at least $\kappa/3$, i.e., chosen from $\cG_d(N_x, 1/2, \kappa')$ for some $\kappa' > \kappa/3$. When $\bv \subset K$, then after remove $\bv$ and $X$, we remove the entire clique $K$, and the remaining hypergraph $G_i$ is a random graph from $\cG_d(N_x, \frac{1}{2})$. In this case we include all vertices in $\bv$ to $Q$ set.

Based on the idea above, we formalize the proof next. The proof can be divided into two steps.

{\noindent \bf Step 1.} In this step, we consider $\bv \subsetneq K$ in step 3 of Algorithm \ref{alg: HPC recovery based on HPC detection}. When $\bv \subsetneq K$, then each entry of $\bcA_{[v_1,\ldots,v_{d-1}, x]}$ are independent $\Bern(1/2)$ random variables. Thus each vertex in clique will be removed with probability $1/2$. By concentration result, for each $G_x$, we have $\bbP(\kappa' > \frac{\kappa}{3} ) \geq 1- \exp(-C \kappa)$ where $\kappa'$ is the clique size in graph $G_x$. So $\varphi$ will output $\cG_d(N_x, 1/2, \kappa/3)$ with probability at least $1 - \exp(-C \kappa) - 1/(4N^d)$.

{\noindent \bf Step 2.} In this step, we consider $\bv \subset K$. Given $\bv \subset K$, by the construction of Algorithm \ref{alg: HPC recovery based on HPC detection}, the whole clique set $K$ is removed in $G_x$. So $G_x \sim \cG_d(N_i, \frac{1}{2})$, and we have
\begin{equation*}
\begin{split}
    \bbP(\kappa'  > \frac{\kappa}{3}) \leq {N_x \choose \kappa/3} (\frac{1}{2})^{ {\kappa/3 \choose d} } \leq \exp(-C \kappa),
\end{split}
\end{equation*}
here$\kappa'$ is the clique size in graph $G_x$. So $\varphi$ outputs $G_d(N_x, \frac{1}{2})$ with probability at least $1 - \exp(-C\kappa) - \frac{1}{4N^d}$.

By the union bound over all possible choices of $\bv$, there exists $C  > 0$ that for $\kappa \geq C \log N$,
\begin{equation*}
\begin{split}
        \bbP( Q  = K ) &\geq 1 - {N \choose d-1} \exp(-C \kappa)  -  {N \choose d-1} \frac{1}{2N^d} \geq 1 - \frac{1}{N}.
\end{split}
\end{equation*}
This contracts the assumption. So we have finished the proof of this lemma.

\begin{algorithm}[h] \caption{Algorithm for HPC Recovery based on HPC Detection} \label{alg: HPC recovery based on HPC detection}
	\begin{algorithmic}[1]
	\State \textbf{Input:} hypergraph $G$ with adjacency tensor $\bcA$, HPC detection algorithm $\varphi$.
		\State Let $H = \{\bv: \bv \subset [N], |\bv| = d-1  \}$,  $Q = \emptyset$ (representing the current clique).
		\For{$\bv \in H$}
		    \State Denote all vertices in $\bv$ as $\bv_1, \ldots, \bv_{d-1}$. Let $X = \{ x \in [N] \setminus \bv \, |  \, \bcA_{[v_1,v_2,\ldots, v_{d-1}, x]} = 1 \}$, $G_x = G \setminus \{\bv \bigcup X\}$, and $N_x = |G_x|$.
		    \State Input $G_x$ to the HPC detection algorithm $\varphi$. If $\varphi$ outputs $\cG_d(N_x, \frac{1}{2})$, then set $Q = Q \bigcup \bv$.
		\EndFor
		\State \textbf{Output:} $Q$.
	\end{algorithmic}
\end{algorithm}

\subsection{Proof of Theorem \ref{thm: Metropolis process is hard to find large clique}}
We begin the theorem by introducing the concept "m-gateway" proposed in \cite{jerrum1992large}. A state $K$ is callled m-gateway if there exists a sequence of states $K_0, K_1, \ldots, K_s$ such that $K_0 = K, K_i \oplus K_{i+1} = 1, 0 \leq i \leq s-1$,  $|K_s| = m$ and $|K_i| > |K_0|$ for $1 \leq i \leq s$.

Next, we introduce Lemma \ref{lm: m-gateway prop in hypergraphical plantedclique} which is useful in the proof of this theorem.
\begin{Lemma} \label{lm: m-gateway prop in hypergraphical plantedclique}
Suppose $0< \epsilon < 1, 0 < \beta < \frac{1}{2}$ and $\frac{1}{3} \epsilon < \frac{1}{2} - \beta$. For $G \sim \mathcal{G}_d (N, \frac{1}{2}, N^\beta)$ and let $m = 2k - \lceil \left((1+ \frac{2}{3} \epsilon) (d-1)! \log_2 N\right)^{\frac{1}{d-1}} \rceil $ where $k = \lceil \left((1+\frac{2}{3}\epsilon) \frac{d!}{2} \log_2 N\right)^{\frac{1}{d-1}} \rceil $. Let $\rho(G)$ the proportion of size k cliques in $G$ that are m-gateway. Then $\rho(G) \leq N^{-\Omega \left(( \log_2 N)^{\frac{1}{d-1}} \right)}$ for almost every $G$.  
\end{Lemma}

Recall $\Gamma$ denotes the collection of all cliques in $G$ and $\Gamma_k \subseteq \Gamma$ is the set of k-cliques in $G$. Let $Q\subseteq \Gamma_k$ be the set of all k-cliques that are m-gateways, where $k, m$ are quantities defined in Lemma \ref{lm: m-gateway prop in hypergraphical plantedclique}. Divide cliques in $G$ into two sets $S$ and $\bar{S}:= \Gamma \setminus S$ where $S$ can be reached without passing through $Q$. It is easy to see that $\Gamma_k  \subseteq S$ and all m-cliques are in $\bar{S}$. The intuition underlying the proof is that $Q$ is rather small which makes it hard to transit from $S$ to m-cliques.

First we calculate the probability, in the stationary distribution, of transiting from $S$ to $\bar{S}$ conditional on being in $S$, 
\begin{equation}\label{eq: transition probability}
    \Phi_S := \bbP \left( \text{ transit from S to }\bar{S} | \text{ being in } S  \right) = \sum_{K \in S, K' \in \bar{S}} \pi(K) \bbP(K, K')/(\sum_{K \in S} \pi(K) ).
\end{equation}
By the definition of $\bar{S}$, to transit to $\bar{S}$, we need to pass nodes in $Q$. So the numerator in \eqref{eq: transition probability} is bounded by $\pi(Q)$. Since $\Gamma_k \subseteq S$, $\pi(S)$ in the denominator is greater or equal to $\pi(\Gamma_k)$. Thus the conditional transition probability in \eqref{eq: transition probability} is no more than $\frac{\pi(Q)}{\pi(\Gamma_k)}$, which is the proportion of $k-$clique that are m-gateways in Lemma \ref{lm: m-gateway prop in hypergraphical plantedclique}. Based on the results in Lemma \ref{lm: m-gateway prop in hypergraphical plantedclique}, we have
\begin{equation*}
    \Phi_S \leq N^{-\Omega( (\log_2 N)^{\frac{1}{d-1}}  )  }.
\end{equation*}

Next, we make rigorous argument that the restriction on $Q$ makes it hard to transit from $S$ to $\bar{S}$. We modify the Metropolis process to make states in $\bar{S}$ absorbing; this is done by setting $\bbP(K', K) = \delta_{K' K}$ for all $K' \in \bar{S}$ where $\delta_{K' K}$ is the Kronecker delta. Also define the initial distribution $\pi_0$ by 
\begin{equation*}
    \pi_0(K) = \left\{ \begin{array}{lc}
         \pi(K)/\pi(S), & \text{ if } K \in S;  \\
         0,& \text{otherwise}.
    \end{array}\right. 
\end{equation*}
Note that given the initial distribution $\pi_0$, the probability that the Metropolis process transit to $S$ in the first step is $\Phi_S$. Also since the states in $\bar{S}$ are absorbing, for any fixed $K \in S$, the probability $\pi_t(K)$ of being in state $K$ at time $t$ is a monotonically decreasing function of $t$. So the probability of transition from $S$ to $\bar{S}$ in each subsequent step is bounded above by $\Phi_S$. Hence the expected time of first entry into $\bar{S}$, given initial distribution $\pi_0$, is bounded below by $\frac{1}{2\Phi_S}$. Clearly, there must be some choice of initial state from which the expected time to reach $\bar{S}$ (and hence a clique of size $m$) is at least $\frac{1}{2\Phi_S}$. This has finished the proof.

\subsubsection{Proof of Lemma \ref{lm: m-gateway prop in hypergraphical plantedclique}}
Let $\cX$ be the set of all pairs $(G, K)$ where $K$ is the clique with size $k$ in $G$ and $\cY$ be the set of pairs $(G,K) \in \cX$ such that $K$ is also a m-gateway.

Let $V = \{0,1, \ldots, n-1\}$ be the set of all nodes and $Q = \{0,1, \ldots, \kappa-1\}$ be the node set of planted clique. Define $f(t) = {\kappa \choose t} {N-\kappa \choose k-t} \left(\frac{1}{2}\right)^{ {k \choose d} - {t \choose d} }$ as the probability that $t$ nodes in $K$ come from $Q$ and rest of $(k-t)$ nodes in $K$ come from $V \setminus Q$ and let $F = \sum_{t=0}^k f(t)$. In the following we define a sampling way of sampling $(G, K)$ from $\cX$ such that it has the same distribution as uniformly sample $(G, K)$ from k-cliques in $G \sim \cG_d(N, \frac{1}{2}, \kappa)$ and we call this sampling strategy ``uniform sampling from $\cX$". The uniform sampling of $(G, K)$ from $\cX$ is the following:
\begin{itemize}
    \item Pick $t \in [0,k]$ with probability $\frac{f(t)}{F}$.
    \item Select $K'$ of size $t$ uniformly at random from $Q$ and select $K''$ of size $(k-t)$ uniformly at random from the subset of $V \setminus Q$, and set $K = K' + K''$. 
    \item Include all edges that have both endpoints in $Q$ or have both endpoints in $K$ in $G$; decide whether to include the remaining potential edges in $G$ with probability $\frac{1}{2}$.
\end{itemize}
We first show that the size of $K'$ is often small when $k = \lceil \left( \frac{d!}{2} (1 + \frac{2}{3}\epsilon ) \log_2 N  \right)^{\frac{1}{d-1}} \rceil $. Using the fact ${\kappa \choose t} \leq \kappa^t$, ${N-\kappa \choose k-t} \leq {N-\kappa \choose k} \left(\frac{2k}{N}\right)^t$ and $\kappa = N^\beta$, $f(t)$ could be upper bounded in the following way
\begin{equation*}
\begin{split}
    f(t) & \leq \kappa^t {N-\kappa \choose k} \left( \frac{2k}{N} \right)^{t} \left(\frac{1}{2}\right)^{ {k \choose d} - {t \choose d} } \\
    & \leq \left( \kappa \frac{2k}{N} 2^{ \frac{(t-1) \ldots (t-d+1) }{d!} }  \right)^t f(0)\\
    & \leq \left( \kappa \frac{2k}{N} 2^{ \frac{k^{d-1}}{d!} }  \right)^t f(0) \leq \left( 2k N^{ -\frac{1}{2} + \beta + \frac{1}{3}\epsilon } \right)^t f(0) \\
    & \leq f(0) N^{-ct},
\end{split}
\end{equation*} where the last inequality is because $\beta - \frac{1}{2} < - \frac{1}{3}\epsilon$. So $f(0) = 1-N^{-c}$ and for any $t^* > 0$, 
\begin{equation}\label{ineq: tail prob of t}
    \bbP(t \geq t^*) \leq N^{-c t^*}.
\end{equation} 

For $(G, K)$ sampled uniformly at random from $\cX$, we show that the probability $(G,K) \in \cY$ is $N^{ - \Omega\left(  (\log_2N)^{\frac{1}{d-1}} \right) }$. When $(G, K) \in \cY$, by definition $K$ is a m-gateway. Consider a path that lead the Metropolis process from $K$ to a $m$-clique and denote $K^*$ as the first clique in this path satisfying $|K^* \setminus K| = m-k$. Set $A = K^* \setminus K$. Then, $a = |A| = m - k$. Since $|K^*| > k$, $|K^* \bigcap K| > k-a $, there exists a set $B \subseteq K$ of cardinality $b = k-a =2k-m$ such that the bipartite subgraph of $G$ induced by $A$ and $B$ is complete. So condition on $t \leq t^*$, the probability that $(G, K) \in \cY$ is less than the probability of the existence of the complete bipartite graph between $A$ and $B \setminus K'$ where $|B \setminus K'| \geq b- t^*$, i.e.,
\begin{equation}\label{ineq: conditional upper bound of prob Y}
\begin{split}
    & \bbP\left((G, K) \in \cY | (G, K) \in \cX, \left| B \bigcap Q \right| \leq t^* \right) \leq {N-k \choose m-k} {k \choose 2k-m} 2^{- (m-k) {b-t^* \choose d-1} }\\
    &\leq {N \choose m-k} {k \choose m-k} 2^{- (m-k) {b-t^* \choose d-1} }\\
    & \leq \left( \frac{eN}{m-k} \frac{ek}{m-k} 2^{- {2k-m -t^* \choose d-1} }    \right)^{m-k}\\
    & \leq \left( \frac{eN}{m-k} \frac{ek}{m-k} 2^{ -\frac{ (2k-m-t^*-d+2)^{d-1} }{(d-1)!} }    \right)^{m-k}\\
    & \leq \left( \frac{eN}{m-k} \frac{ek}{m-k} N^{-1 - \frac{1}{6}\epsilon}   \right)^{m-k} \leq N^{-\Omega\left( (\log_2 N )^{\frac{1}{d-1}} \right)},
\end{split}
\end{equation}
here the first inequality is because all hyper-edges between $|B \setminus K'|$ and $A$ are connected with probability $1/2$ by construction and in the forth inequality, we choose $t^*$ such that $2k-m-t^* = \lceil \left((1+ \frac{1}{3} \epsilon) (d-1)! \log_2 N\right)^{\frac{1}{d-1}} \rceil$ and for large enough $N$, we have $ 2k-m -t^* - d+ 2 \geq \lceil \left((1+ \frac{1}{6} \epsilon) (d-1)! \log_2 N\right)^{\frac{1}{d-1}} \rceil$.

Using Bayesian formula and \eqref{ineq: tail prob of t}, \eqref{ineq: conditional upper bound of prob Y} with the choice of $t^* = 2k-m - \lceil \left((1+ \frac{1}{3} \epsilon) (d-1)! \log_2 N\right)^{\frac{1}{d-1}} \rceil$, and marginalizing $ \left| B \bigcap Q \right| \leq t^*$ we get 
\begin{equation} \label{ineq: transition probability bound}
    \bbP\left((G,K) \in \cY | (G, K) \in \cX \right) \leq N^{-\Omega( (\log_2 N)^{\frac{1}{d-1}} )}.
\end{equation}
Given $G \in \cG_d(N, \frac{1}{2}, \kappa)$ with $ \kappa = N^\beta$, let $X=X(G)$ be the number of $k$-cliques in $G$ and $Y = Y(G)$ be the number of $k$-cliques that are also m-gateways. By the uniform sampling property and \eqref{ineq: transition probability bound}, we have
\begin{equation*}
    \frac{\bbE (Y)}{\bbE (X)} = \bbP\left((G,K) \in \cY | (G, K) \in \cX \right)\leq N^{-\Omega( (\log_2 N)^{\frac{1}{d-1}} )}.
\end{equation*} So to show the result $\frac{Y}{X} = N^{-\Omega( (\log_2 N)^{\frac{1}{d-1}} )}$, we only need to show that $X, Y$ are concentrated around $\bbE (X), \bbE (Y)$. Let $c(N)$ be a sequence that goes to $\infty$ as $N \to \infty$, since $\cY$ is non-negative, by Markov inequality, $\bbP(\cY \geq \bbE(Y) c(N)) \leq \frac{1}{c(N)} \to 0$. Since the subgraph of $G$ induced by $V\setminus Q$ is a random Erd\H{o}s-R{\'e}nyi graph on $N-\kappa$ vertices, by the same argument of \cite{bollobas2001random} p284, the number of $k$-clique in $G$ is concentrated and we have
$$X \geq \frac{3}{4} {N-\kappa \choose k} 2^{-{k \choose d}} + {\kappa \choose k}.$$ 
At the same time, $\bbE(X) = {N-\kappa \choose k} 2^{-{k \choose d}} + {\kappa \choose k}$, and this yields $X \geq \frac{1}{2} \bbE(X)$ for almost every $G$. Finally we have
\begin{equation*}
\frac{Y}{X} \leq \frac{c(N) \bbE(Y)}{\frac{1}{2}\bbE(X)} \leq 2 c(N) N^{-\Omega( (\log_2 N)^{\frac{1}{d-1}} )} \leq O( N^{-\Omega( (\log_2 N)^{\frac{1}{d-1}} )}),
\end{equation*} the last inequality is because we can choose $c(N)$ such that it grows slow enough.

\subsection{Proof of Theorem \ref{prop: fail of low-degree polynomial test}}
We first introduce some preliminary results in literature we will use in the proof. The following Proposition \ref{prop: trun likelihood ration} comes from \cite{hopkins2018statistical}.

\begin{Proposition}[Page 35 of \cite{hopkins2018statistical}]\label{prop: trun likelihood ration}
	Let likelihood ratio be $\LR(x) = \frac{p_{H_1}(x)}{p_{H_0} (x)}: \Omega^n \to \bbR$. For every $D \in \mathbb{Z}_{+}$, we have 
	\begin{equation*}
		 \frac{ \LR^{\leq D} - 1}{\|\LR^{\leq D} - 1\|} = \argmax_{\substack{f \in D-\text{simple}:\\ \bbE_{H_0} f^2(X) = 1,\bbE_{H_0} f(X) = 0}} \bbE_{ H_1} f(X)
	\end{equation*}
	and 
	\begin{equation*}
		\|\LR^{\leq D} - 1\| = \max_{\substack{f: \bbE_{H_0} f^2(X) = 1,\\ \bbE_{H_0} f(X) = 0}} \bbE_{H_1} f(X),
	\end{equation*}
	where $\|f\| = \sqrt{\bbE_{H_0}f^2(X)}$ and $f^{\leq_{H_0} D}$ is the projection of a function $f$ to the span of coordinate-degree-D functions, where the projection is orthonormal with respect to the inner product $\langle \cdot, \cdot \rangle_{H_0}$. When it is clear from the context, we may drop the subscript $H_0$.
\end{Proposition}

By the above Proposition, we know to bound $\bbE_{H_1}f(G)$ for $f$ satisfying $\bbE_{H_0}f(G) = 0$ and $\bbE_{H_0}f^2(G) = 1$, we just need to bound $\|\LR^{\leq D} - 1\|$.

Also suppose $D \geq 1$ is fixed, $f_0, f_1, \ldots, f_m: \Omega^n \to R$ are orthonormal basis for the coordinate-degree $D$ functions (with respect to $\langle \cdot, \cdot \rangle_{H_0}$), and that $f_0(x) = 1$ is a constant function. Then by the property of basis functions, we have
\begin{equation*}
\begin{split}
    \|\LR^{\leq D} - 1\|^2 &= \sum_{i=1}^m \langle f_i, \LR^{\leq D} - 1 \rangle^2 = \sum_{i=1}^m \left(  \bbE_{H_0}(f_i(X) (\LR^{\leq D} - 1) )  \right)^2\\
    & \overset{(a)}= \sum_{i=1}^m \left(  \bbE_{H_0}(f_i(X) \LR(X) )  \right)^2 = \sum_{i=1}^m (\bbE_{H_1} f_i(X))^2,
\end{split}
\end{equation*} here (a) is because $\LR - (\LR^{\leq D} - 1)$ is orthogonal to $f_i$ by assumption.

	Now, we start the proof of the main result. We consider a simple variant of hypergraphic planted clique model with $G \sim \cG_d(N, \frac{1}{2})$ and each vertex is included in the clique set with probability $\frac{\kappa}{N}$ and denote its adjacency tensor as $\bcA$. By concentration result, the clique size of the modified hypergraphic planted clique is of order $\kappa$ and it is easy to see that if we could solve above modified hypergraphic planted clique problem, then we can solve the original hypergraphic planted clique problem with high probability. In the HPC problem, by the Boolean Fourier analysis \cite{o2014analysis}, functions $\{ \chi_{\alpha} (\bcA) = \prod_{(i_1, \ldots, i_d) \in \alpha} (2 \bcA_{[i_1, \ldots, i_d]} - 1)  \}_{\alpha \subseteq \text{set}{N \choose d}, |\alpha| \leq |D| } $ form an orthonormal basis for the degree-D functions, i.e. $f_i$s mentioned above. Here $\text{set}{N \choose d}$ denotes the set of all possible $d$-tuples chosen from $[N]$.
	Fix $\alpha \subseteq \text{set}{N \choose d}$ and the planted clique indices set $S$. Then $$\bbE \chi_\alpha (\bcA) = \bbE_S \prod_{ (i_1, \ldots, i_d) \in \alpha } \bbE \left[ (2\bcA_{[i_1, \ldots, i_d]} - 1)|S \right].$$ This is non-zero if and only if $V(\alpha) \subseteq S$ where $V(\alpha)$ is the vertex set appear in $\alpha$. So we have 
	\begin{equation}\label{eq: expectation on H1 bound}
	    \bbE \chi_\alpha (\bcA) = (\frac{\kappa}{N})^{|V(\alpha)|}.
	\end{equation}
	If $|\alpha| \leq D$, then $|V(\alpha)| \leq d D$ and for every $t \leq dD$, we can compute the number of sets $\alpha $ with $|V(\alpha)| = t$ and $|\alpha| \leq D$ is ${ N \choose t } { {t \choose d} \choose |\alpha|}$. 
	
	Let $D = C\log N, \kappa = N^{\frac{1}{2} - \epsilon}$, then 
	\begin{equation}
		\begin{split}
		 & \sum_{0 < |\alpha| \leq D} (\bbE_{H_1} \chi_\alpha (\bcA) )^2 \overset{(a)}\leq \sum_{t \leq dC \log N} \left(\frac{\kappa}{N}\right)^{2t} { N \choose t } { {t \choose d} \choose |\alpha|}\\
			\leq & \sum_{t \leq dC \log N} N^{-t -2\epsilon t} N^t {t\choose d}^{\min (|\alpha|, {t\choose d}/2 )}\\
			\leq & \sum_{t \leq dC \log N} N^{-2\epsilon t} (t^d)^{ \min (|\alpha|, t^d ) }\\
			\leq & \sum_{t \leq (C \log N)^{\frac{1}{d}}} N^{-2\epsilon t} t^{dt^d} + \sum_{(C \log N)^{\frac{1}{d}} \leq t \leq dC \log N} N^{-2\epsilon t}t^{d C \log N} \\
			= & \sum_{t \leq (C \log N)^{\frac{1}{d}}} \exp(-2\epsilon t \log N + dC \log N \log t) \\
		& \quad 	+ \sum_{(C \log N)^{\frac{1}{d}} \leq t \leq dC \log N} \exp(-2\epsilon t \log N + dC \log N \log t)\\
			\leq & \sum_{t \leq (C \log N)^{\frac{1}{d}}} \exp(-\epsilon t \log N)+ \sum_{(C \log N)^{\frac{1}{d}} \leq t \leq dC \log N} \exp(-\epsilon t \log N) = O(1),
		\end{split}
	\end{equation}
where (a) is due to \eqref{eq: expectation on H1 bound} and last inequality is due to the fact that when $N$ is large enough, we have $\epsilon t \log N \geq C\log t \log N$.

\subsection{Proof of Theorem \ref{thm: low-degree-HPDS-evidence} } \label{sec: proof-conjecture2}
We begin by introducing a few notation and preliminaries from \cite{schramm2020computational} before proving our main results.

Recall $G \in \cG_d(N,\kappa,q_1, q_2)$ and its adjacency tensor is $\bcA$. Let $\bcX = \bbE (\bcA|K)$ where $K$ is a size $\kappa$ subset of $[N]$ drawn uniformly at random and $X = (\bcX_{[i_1,\ldots,i_d]})_{i_1 < \cdots < i_d} \in \{q_1, q_2\}^{\binom{N}{d}}$. First, define the minimum mean squared error of $f$ with degree at most $D$ as $$\text{MMSE}_{\leq D} := \inf_{f: \text{ degree of } f \text{ is at most } D } \bbE (f(\bcA) - v_1)^2.$$
$\text{MMSE}_{\leq D}$ is closely related to the quantity called {\it degree-D maximum correlation}: 
\begin{equation*}
	\text{Corr}_{\leq D} := \sup_{\substack{f: \text{ degree of } f \text{ is at most } D,\\ \bbE (f^2(\bcA)) = 1}} \bbE (f(\bcA) \cdot v_1).
\end{equation*} 
\cite{schramm2020computational} showed that $\text{MMSE}_{\leq D} = \bbE(v_1^2) - \text{Corr}_{\leq D}^2$. Notice $\bbE(v_1^2) = \bbE(v_1) = \frac{\kappa}{N} $, so to lower bound $\text{MMSE}_{\leq D}$, we just need to upper bound $\text{Corr}_{\leq D}^2$. Theorem 2.7 of \cite{schramm2020computational} gives an upper bound for $\text{Corr}_{\leq D}^2$ in the general binary observation model. Specifying it in our context, we have
\begin{Proposition}
	In the HPDS recovery problem,
	\begin{equation}\label{ineq: correlation-upper-bound}
		\text{Corr}^2_{\leq D} \leq \sum_{\alpha \in \{0,1 \}^{ N \choose d }, 0 \leq |\alpha| \leq D} \frac{\kappa_\alpha^2}{(q_2(1-q_1))^{|\alpha|}},
	\end{equation} where $\kappa_\alpha$ for $\alpha \in  \{0,1 \}^{ N \choose d }$ is defined recursively by
	\begin{equation*}
		\kappa_\alpha = \bbE(v_1 X^{\alpha}) - \sum_{0 \leq \beta \lneq \alpha } \kappa_\beta \bbE(X^{\alpha - \beta}).
	\end{equation*} Here for $\alpha,\beta \in  \{0,1 \}^{ N \choose d }$, $|\alpha| := \sum_{i=1}^{ N \choose d } |\alpha_i|$, $X^\alpha := \prod_{i=1}^{ N \choose d } X_i^{\alpha_i}$ and $0 \leq \beta \lneq \alpha$ means $\beta \neq \alpha$ and for $i = 1,\ldots,{ N \choose d }$, we have $0 \leq \beta_i \leq \alpha_i$.
\end{Proposition}
In the hypergraph setting, we can view $\alpha$ as $(\alpha_{[i_1,\cdots,i_d]})_{i_1< \cdots < i_d}$ and it is a realization of hypergraph on the vertex set $[N]$, where $\alpha_{[i_1,\cdots,i_d]}$ encodes whether a hyperedge exists between vertices $i_1,\ldots,i_d$. The quantity $\kappa_\alpha$ can be interpreted as a certain joint cumulant between $v_1$ and $X_i$ for $i = 1,\ldots,{N \choose d}$ and we refer readers to Remark 2.3 and Appendix D of \cite{schramm2020computational} for detailed discussion.

By the proof of Theorem 2.9 and Lemma 3.4 in \cite{schramm2020computational}, we have the following upper bound for $\kappa_\alpha$:
\begin{Lemma} First, $\kappa_0 = \frac{\kappa}{N}$ and for $\alpha$ such that $|\alpha| \geq 1$, we have
\begin{equation} \label{ineq: k-alpha-upper-bound}
	\frac{\kappa_\alpha}{(q_2(1-q_1))^{|\alpha|/2}} \leq (|\alpha| + 1)^{|\alpha|} \theta^{|\alpha|} \kappa_0^{|V(\alpha)|},
\end{equation} here $\theta:= \frac{q_1 - q_2}{\sqrt{q_2(1-q_1)}} $ and $V(\alpha) \subseteq [N]$ denotes the set of vertices spanned by $\alpha$.
\end{Lemma}
We note that Lemma 3.4 of \cite{schramm2020computational} is established when $d = 2$, i.e., graph setting, while it is straightforward to extend it to the hypergraph setting. For brevity, we omit the proof here.

Next we discuss the connectivity of two vertices in hypergraph. We say vertex $i$ and $j$ are connected if there exists a hyperedge path $e_1,\ldots,e_m$ such that $e_1$ contains vertex $i$, $e_m$ contains vertex $j$ and $e_k$ has at least one common vertex with $e_{k-1}$ for $k = 2,\ldots,m$. Following the same proof as Lemma 3.2 of \cite{schramm2020computational}, we have
\begin{Lemma}\label{lm:kappa_alpha}
	If $\alpha$ has a nonempty connected component without vertex $1$, then $\kappa_\alpha = 0$.
\end{Lemma} 
Based on Lemma \ref{lm:kappa_alpha}, we only need to consider the connected hypergraph $\alpha$ containing vertex $1$ when bounding $\text{Corr}^2_{\leq D}$. 

The final ingredient we need is to bound the number of connected hypergraph $\alpha$ containing vertex $1$ in different settings. In the graph case, we can enumerate connected graph $\alpha$ by its number of edges and spanned vertices and it is relative easy to bound the number of connected $\alpha$ when fixing $|\alpha|$ and $|V(\alpha)|$ due to the fact that adding a new edge to a connected graph can introduce at most one new vertex. While in hypergraph, it is difficult to count the number of connected $\alpha$ given $|\alpha|$ and $|V(\alpha)|$ as now adding a new edge to a connected hypergraph can introduce $0$ to $d-1$ different number of new vertices. To overcome this difficulty, we introduce another characteristic of $\alpha$. 
Given any connected hypergraph $\alpha$ containing vertex $1$, we can divide its hyperedges into two groups: in the first group, the edges are sorted such that the first edge contains vertex $1$, the second edge has at least one common vertex with the first edge and at the same time has at least one new vertex that does not appear in previous edges. Similar ideas applies to the rest of the edges in group $1$, i.e., the $i$th edge in group 1 needs to have at least one common vertex with the edges appeared before and also introduce as least one new vertex. The second group contains edges where all their vertices have appeared in the first group. Denote the number of edges in group $1$ as $g_1(\alpha)$ and the number of edges in group $2$ as $g_2(\alpha)$. We note that given $\alpha$, the grouping and values of $g_1(\alpha)$, $g_2(\alpha)$ may not be unique.

The following Lemma plays a key role in proving our main results.
\begin{Lemma} \label{lm: HPDS-connected-hypergraph-counting}
	Given a hypergraph $G \sim \cG_d$ on vertex set $[N]$. For integer $w \geq 1$ and $0 \leq h \leq w-1$, suppose $N$ is sufficiently large and $d,w = o(N)$, the number of connected hypergraph $\alpha$ satisfying (i) $|\alpha| = w$, (ii) $1 \in V(\alpha)$, and (iii) $g_1(\alpha) =w - h$ is at most $(wd(d-1))^{w-h} N^{(w-h)(d-1)} (w(d-1)+1)^{hd}$.
\end{Lemma}

Now, we are ready to bound $\text{Corr}^2_{\leq D}$.
\begin{equation*}
	\begin{split}
		& \text{Corr}^2_{\leq D} \overset{ \eqref{ineq: correlation-upper-bound} }\leq \kappa_0^2 + \sum_{\alpha \in \{0,1 \}^{ N \choose d }, 1 \leq |\alpha| \leq D} \frac{\kappa_\alpha^2}{(q_2(1-q_1))^{|\alpha|}}\\
		\leq & \kappa_0^2 + \sum_{1 \leq |\alpha| \leq D} \frac{\kappa_\alpha^2}{(q_2(1-q_1))^{|\alpha|}} \\
		\overset{(a)}\leq & \kappa_0^2 + \sum_{w=1}^D \sum_{h=0}^{w-1} (wd(d-1))^{w-h} N^{(w-h)(d-1)} (w(d-1)+1)^{hd}  (w + 1)^{2w} \theta^{2w} \kappa_0 ^{2((w-h)(d-1) + 1 )}\\
		\leq & \kappa_0^2 \left( 1+ \sum_{w=1}^D \sum_{h=0}^{w-1} (D(D+1)^2 d(d-1) N^{d-1} \kappa_0^{2(d-1)} \theta^2 )^w \left( \frac{(D(d-1)+1)^d}{Dd(d-1) N^{d-1} \kappa_0^{2(d-1)} } \right)^{h} \right)\\
		\leq & \kappa_0^2 \left( 1+ \sum_{h=0}^{D-1} \sum_{w=h+1}^{D} (D(D+1)^2 d(d-1) N^{d-1} \kappa_0^{2(d-1)} \theta^2 )^{w-h} \left( (D(d-1)+1)^d(D+1)^2 \theta^2 \right)^{h} \right)\\
		\overset{ (b) }\leq & \kappa_0^2 \left( 1+ \sum_{h=0}^{D-1} r^h \sum_{w=h+1}^{D} r^{w-h} \right)\\
		\leq & \kappa_0^2 \left( 1+ r\sum_{h=0}^{\infty} r^h \sum_{w=h+1}^{\infty} r^{w-h-1} \right) = \kappa_0^2 (1 +  \frac{r}{(1-r)^2} ).
	\end{split}
\end{equation*} 
Here, (a) is due to Lemma \ref{lm: HPDS-connected-hypergraph-counting}, \eqref{ineq: k-alpha-upper-bound} and the fact that when $g_1(\alpha) = w-h$, $|V(\alpha)| \leq (w-h)(d-1)+1$; (b) is because of \eqref{ineq: HPDS-low-degree-SNR-lower-bound}, $\theta= \frac{q_1 - q_2}{\sqrt{q_2(1-q_1)}} $ and $\kappa_0 = \frac{\kappa}{N}$.

The first conclusion follows by observing that $\text{MMSE}_{\leq D} = \bbE(v_1^2) - \text{Corr}_{\leq D}^2$, $\bbE(v_1^2) = \kappa_0 = \frac{\kappa}{N} $. 

For the second conclusion, we first notice that the trivial estimator $f_0(\bcA) = \frac{\kappa}{N}$ achieves the mean square error $\frac{\kappa}{N} - (\frac{\kappa}{N})^2$ in estimating $v_1$. Also $q_2 < q_1 < 1- \Omega(1)$ and $$ \liminf_{N \to \infty} \log_N \kappa \geq \frac{1}{2} \quad \text{and}\quad  \limsup_{N \to \infty} \log_N \left(\frac{\kappa^{d-1} (q_1-q_2)}{\sqrt{
			q_2(1-q_2)}}\right) < \frac{d}{2} -\frac{1}{2}$$ implies \eqref{ineq: HPDS-low-degree-SNR-lower-bound} holds. So the second conclusion $\liminf_{N\to \infty}\frac{\mathbb{E}(f(\bcA) - v_1)^2}{\mathbb{E}(f_0(\bcA) - v_1)^2} \geq 1$ follows by the first part of the result and set $r = o(1)$ that decays slowly to zero. This has finished the proof of this Theorem.
\subsubsection{Proof of Lemma \ref{lm: HPDS-connected-hypergraph-counting} }
Given any $\alpha$, we can sort its edges as we mentioned in Section \ref{sec: proof-conjecture2}. Next we bound the number of choices for each edge in group 1 and group 2.

{\bf For Group 1}.
\begin{itemize}
	\item The first edge contains vertex $1$, so there are at most ${N \choose d-1}$ choices for choosing the rest of the $(d-1)$ vertices.
	\item When choose the second edge, the number of its choices is bounded by $\sum_{x=1}^{d-1} {d \choose d-x} {N \choose x} \leq d(d-1) {N \choose d-1} $, where $x$ in the first summand denotes the number of vertices we pick from edge $1$ and the inequality is because ${d \choose d-x} {N \choose x}$ is an increasing function with respect to $x$ when $d,w = o(N)$.
	\item Using the same idea in choose the second edge, the number of choices for the third edge is at most: $ \sum_{x=1}^{d-1} {2(d-1)+1 \choose d-x} {N \choose x} \leq (2(d-1)+1)(d-1) {N \choose d-1} $
	\item $\cdots$
	\item For the $(w-h)$th edge, the number of its choices is bounded by: $ \sum_{x=1}^{d-1} {(w-h-1)(d-1)+1 \choose d-x} {N \choose x} \leq ((w-h-1)(d-1)+1)(d-1) {N \choose d-1} $.
\end{itemize} 
In total, the number of choices in group $1$ is at most $$((w-h-1)(d-1)+1)^{w-h} (d-1)^{w-h} \binom{N}{d-1}^{w-h} \leq (wd(d-1))^{w-h} N^{(w-h)(d-1)}.$$

{\bf For Group 2}. Each edge left has at most $\binom{|V(\alpha)|}{d} \leq |V(\alpha)|^d \leq (w(d-1)+1)^d$ choices. Since there are $h$ edges in group $2$, we have at most $(w(d-1)+1)^{hd}$ choices. This finishes the proof of this lemma. 
\end{appendix}

\end{document}